\documentclass{patmorin}
\usepackage{pat}
\usepackage[T1]{fontenc}
\usepackage[utf8]{inputenc}
\usepackage[inline]{enumitem}

\usepackage{mathtools}  %
\usepackage[longnamesfirst,numbers,sort&compress]{natbib}

\crefformat{equation}{#2(#1)#3}
\let\eqref\cref
\crefformat{subsection}{Subsection #2#1#3}
\crefformat{subsubsection}{Subsubsection #2#1#3} 

\renewcommand{\ge}{\geqslant}
\renewcommand{\le}{\leqslant}
\renewcommand{\geq}{\geqslant}
\renewcommand{\leq}{\leqslant}

\DeclareMathOperator{\id}{id}

\definecolor{brightmaroon}{rgb}{0.76, 0.13, 0.28}
\definecolor{linkblue}{rgb}{0, 0.337, 0.227}
\newcommand{\defin}[1]{\emph{\textcolor{brightmaroon}{#1}}}
\makeatletter
\def\mathcolor#1#{\@mathcolor{#1}}
\def\@mathcolor#1#2#3{%
  \protect\leavevmode
  \begingroup
    \color#1{#2}#3%
  \endgroup
}
\makeatother
\newcommand{\mathdefin}[1]{\mathcolor{brightmaroon}{#1}}

\newcommand{\torso}[2]{{#1}\langle{#2}\rangle}
\DeclareMathOperator{\height}{height}
\DeclareMathOperator{\depth}{depth}
\DeclareMathOperator{\tw}{tw}
\DeclareMathOperator{\lca}{lca}
\DeclareMathOperator{\parent}{p}
\DeclareMathOperator{\lgg}{\hat{lg}}
\DeclareMathOperator{\bin}{bin}
\DeclareMathOperator{\idx}{idx}

\newcommand{\Oh}{O}
\DeclarePairedDelimiter{\set}{\{}{\}}

\DeclarePairedDelimiter{\ceil}{\lceil}{\rceil}

\title{\MakeUppercase{\boldmath  You Can Have It All:  Small sparse induced-universal graphs for minor-closed classes}}

 \author{
 Vida Dujmovi{\'c}\,\footnote{School of Computer Science and Electrical Engineering, University of Ottawa, Ottawa, Canada (\texttt{vida.dujmovic@uottawa.ca}). Research supported by NSERC and a University of Ottawa Research Chair.}
 \qquad
 Cyril Gavoille\footnote{LaBRI, University of Bordeaux, France (\texttt{gavoille@labri.fr}). This work is supported by the French ANR Projets ENEDISC (ANR-24-CE48-7768-01) and TEMPOGRAL (ANR-22-CE48-0001).}
 \qquad
  Gwena\"el Joret\footnote{D\'epartement d'Informatique, Universit\'e libre de Bruxelles, Belgium ({\tt gwenael.joret@ulb.be}). G.\ Joret is supported by the Belgian National Fund for Scientific Research (FNRS) and by the Australian Research Council.}\\[1ex]
Piotr Micek\footnote{Department of Theoretical Computer Science, Faculty of Mathematics and Computer Science, Jagiellonian University, Kraków, Poland (\texttt{piotr.micek@uj.edu.pl}). Research supported by the National Science Center of Poland under grant UMO-2023/05/Y/ST6/00079 within the WEAVE-UNISONO program.}
 \qquad
 Pat Morin\footnote{School of Computer Science, Carleton University, Ottawa, Canada (\texttt{morin@scs.carleton.ca}). Research supported by NSERC.}
 \qquad
 David~R.~Wood\footnote{School of Mathematics, Monash University, Melbourne, Australia (\texttt{david.wood@monash.edu}). Research supported by the Australian Research Council and by NSERC.}
 }

\date{}

\newcommand{\calA}{\mathcal{A}}

\newcommand{\calC}{\mathcal{C}}

\newcommand{\calF}{\mathcal{F}}
\newcommand{\calG}{\mathcal{G}}
\newcommand{\calH}{\mathcal{H}}

\newcommand{\calQ}{\mathcal{Q}}
\newcommand{\calR}{\mathcal{R}}
\newcommand{\calS}{\mathcal{S}}
\newcommand{\calT}{\mathcal{T}}

\begin{document}
\hypersetup{pageanchor=false}
\begin{titlepage}
\maketitle

\begin{abstract}
  We show that, for every fixed proper minor-closed graph class $\calG$ and every positive integer $n$, there exists a graph $U_n$ with $n^{1+o(1)}$ vertices and edges such that every $n$-vertex graph in $\calG$ is isomorphic to an induced subgraph of $U_n$.  

  This result generalizes and unifies the results of Dujmovi\'c et al.\ (\textit{J. ACM} 2021) on induced-universal graphs and the results of Esperet, Joret, and Morin (\textit{J. London Math. Soc.} 2023) on sparse induced-universal graphs from planar graphs to all proper minor-closed graph classes.
  It simultaneously improves the previous bounds for the number of vertices in induced-universal graphs (previously  $O(n^2)$ obtained as a corollary of a result of DeVos et al.\ (\textit{J. Comb. Theory, Ser. {B}} 2004) and a result of Gavoille and Labourel (\textit{ESA} 2007)) and  the number of edges in subgraph-universal graphs (previously  $O(n^{3/2})$ by Babai et al.\ (\emph{Theory and Practice of Combinatorics} 1982)). 
\end{abstract}
\end{titlepage}

\newpage
\tableofcontents
\thispagestyle{empty}
\newpage
\hypersetup{pageanchor=true}
\setcounter{page}{1}

\bigskip
\section{Introduction}

Let $\calG$ be a graph class (i.e., a set of graphs closed under isomorphism). A graph $U_n$ is \defin{$n$-universal} for $\calG$ if it ``contains'' every $n$-vertex graph in $\calG$.  Here, the word ``contains'' can take on at least two meanings.  A graph $U_n$ is \defin{$n$-subgraph-universal} for $\calG$ if every $n$-vertex graph in $\calG$ is isomorphic to some subgraph of $U_n$.  A graph $U_n$ is \defin{$n$-induced-universal} for $\calG$ if every $n$-vertex graph in $\calG$ is isomorphic to some \emph{induced} subgraph of $U_n$.  These two notions of universality have largely been treated as separate concepts evaluated by different quality measures. A notable exception is an early work by \citet{chung:universal} that we discuss below.

Clearly, an $n$-subgraph-universal graph must have at least $n$ vertices, and the complete graph $K_n$ contains every $n$-vertex graph as a subgraph. Thus, for every positive integer $n$, $K_n$ is an $n$-vertex graph that is $n$-subgraph-universal for the set of all graphs. Consequently, research on subgraph-universal graphs has focused on minimizing the number of edges, especially when $\calG$ is a sparse class of graphs, such as trees, planar graphs, or $K_t$-minor-free graphs for some fixed $t$.  This line of research has a long and interesting history, dating back to the work of \citet{chung.graham.ea:on} in 1978.  Notable results include subgraph-universal graphs with $\Oh(n\log n)$ edges for trees and outerplanar graphs \cite{chung.graham.ea:on,chung:universal,kaul.kim.ea:universal}, $\Oh(n)$ edges for bounded-degree trees \cite{bhatt.chung.ea:universal}, $\Oh(kn \log (n/k+1))$ edges for graphs of treewidth at most $k$ \cite{kaul.kim.ea:universal}, and $n^{1+o(1)}$ edges for planar graphs \cite{EJM23}.

In the induced-universal setting, the vertex-disjoint union of all $n$-vertex graphs in $\calG$ is $n$-induced-universal, but this graph generally has an enormous number of vertices.\footnote{More precisely, let $\calG_n$ be the set of all $n$-vertex graphs in $\calG$. The universal graph is constructed by choosing one representative graph $G$ from each isomorphism class of $\calG_n$, and then taking the vertex-disjoint union of all such graphs.}  Thus, research on induced-universal graphs has focused on finding graphs with as few vertices as possible.  The subject of induced-universal graphs also has a long history, especially in computer science where it appears under the notion of adjacency labelling schemes (discussed below).  Notable results include induced-universal graphs with $\Oh(n)$ vertices for trees \cite{alstrup.dahlgaard.ea:optimal}, $n(\log(n/k))^{\Oh(k)}$ vertices for graphs of treewidth at most $k$ \cite{gavoille.labourel:shorter}, and a series of results leading to induced-universal graphs for planar graphs having $n\cdot\exp(\Oh(\sqrt{\log n}))$ vertices \cite{kannan.naor.ea:implicit,gavoille.labourel:shorter,chung:universal,spinrad:efficient,BGP22,dujmovic.esperet.ea:adjacency,gawrychowski.janczewski:simpler}.

A 1990 work by \citet{chung:universal} gave several theorems about induced-universal graphs that simultaneously bound both the number of vertices and the number of edges. This work showed the existence of $n$-induced-universal graphs for forests with $\Oh(n\log n)$ vertices and $\Oh(n^2)$ edges. Using the fact that the edges of a planar graph can be partitioned into three forests, she obtained $n$-induced-universal graphs for planar graphs with $\Oh((n\log n)^3)$ vertices and $\Oh(n^6)$ edges. For the next three decades, however, most research either focused on finding subgraph-universal graphs with few edges or on finding induced-universal graphs with few vertices.

This dichotomy persisted until 2023, when \citet{EJM23} showed that, for the class of planar graphs, there exists an $n$-induced-universal graph with $n^{1+o(1)}$ vertices \emph{and edges}.  Thus, their result simultaneously satisfies the stronger requirement of being induced-universal while (nearly) optimizing both quality measures.  Their work is based on a construction of $n$-induced-universal graphs with $n^{1+o(1)}$ vertices for planar graphs due to \citet{dujmovic.esperet.ea:adjacency}.  In the current paper, we extend these results to all proper minor-closed graph classes.

\subsection{New Results}

As mentioned above, induced-universal graphs are closely related to adjacency labelling schemes.  The results in this paper and the techniques used to obtain them are best described using the language of adjacency labelling schemes, which we now define.
Let $\calG$ be a class of graphs and let $f:\N^+\to\R^+$ be a function (where $\N^+$ and $\R^+$ denote the sets of positive integers and positive real numbers, respectively).
We say that $\calG$ admits an \defin{$f(n)$-bit adjacency labelling scheme} if there exists a function $A:(\set{0,1}^*)^2\to\set{0,1}$ (called an \defin{adjacency tester}) such that for all positive integers $n$, for every $n$-vertex graph $G\in\calG$, there exists a function $\ell: V(G)\to\set{0,1}^*$ such that
$|\ell(v)|\leq f(n)$ for each vertex $v$ in $G$ (where $|\alpha|$ denotes the length of a binary string $\alpha$), and such that
for every two vertices $u$, $v$ in $G$,
\[
  A(\ell(u),\ell(v))=\begin{cases}
    0 & \textrm{if $uv\notin E(G)$}, \\
    1 & \textrm{if $uv\in E(G)$}.
  \end{cases}
\]
(Such a function $\ell$ is called an \defin{adjacency labelling} of $G$ for $A$.) 

A graph $H$ is a \defin{minor} of a graph $G$ if a graph isomorphic to $H$ can
be obtained from a subgraph of $G$ by contracting edges. A class of graphs
$\calG$ is \defin{minor-closed} if every minor of every graph in
$\calG$ is also in $\calG$, and it is \defin{proper} if it is not the class of all graphs. Examples of proper minor-closed graph classes include planar graphs, graphs of treewidth at most $k$, graphs embeddable in a fixed
surface, linklessly embeddable graphs, knotlessly embeddable graphs, and, for every fixed graph $X$, the class of graphs with no $X$-minor. Note that every proper minor-closed class is contained in the class of $K_t$-minor-free graphs for some fixed $t$.  In this paper, we prove the following result. (All logarithms are in base $2$.)

\begin{thm}
  \label{main_result}
  Every proper minor-closed class of graphs admits
  a $(1+o(1))\log n$-bit adjacency labelling scheme.
\end{thm}

Note that all the dependence on the fixed proper minor-closed class in~\cref{main_result} is in the $o(\log n)$ term. Furthermore, \cref{main_result} is optimal up to the lower-order term; in our construction, this $o(\log n)$ term is $\Oh\left((\log n)^{3/4}\right)$. The proof of~\cref{main_result} is constructive: for every fixed proper minor-closed class $\calG$, there is a polynomial-time algorithm that takes a graph $G\in\calG$ as input and constructs the labelling $\ell:V(G)\to\{0,1\}^*$. We prove \cref{main_result} by introducing and studying a strengthening of adjacency labelling schemes called (weighted) mixed labelling schemes.

An adjacency labelling scheme for a class of graphs $\calG$ that uses an adjacency tester $A$ can be used to construct an $n$-induced-universal graph for $\calG$ having $\Oh(2^{f(n)})$ vertices, as follows:  We take $U_n$ to be the graph with vertex set $V(U_n):=\{0,1\}^{\le f(n)}$ (the set of binary strings of length at most $f(n)$) and edge set $E(U_n):=\{\alpha\beta: \alpha, \beta \in V(U_n),\; \alpha\neq\beta, \text{ and } A(\alpha,\beta)=1\}$.  For any $n$-vertex graph $G\in\calG$, there exists an adjacency labelling $\ell_G:V(G)\to\{0,1\}^{\le f(n)}$ of $G$ for $A$.  Then the subgraph $U_n[\{\ell_G(v):v\in V(G)\}]$ of $U_n$ induced by $\{\ell_G(v):v\in V(G)\}$ is isomorphic to $G$.\footnote{A minor detail is that the labelling $\ell_G$ needs to be injective, for each $G\in\calG$.  This is the case for the adjacency labelling schemes developed in this paper. See~\cite[Section~2.1]{spinrad:efficient} for more details on the connection between adjacency labelling schemes and induced-universal graphs.} Thus, \cref{main_result} implies the existence of induced-universal graphs with a near-linear number of vertices for $n$-vertex graphs belonging to a fixed proper minor-closed class of graphs.

\begin{cor}\label{main_result_induced_universal}
  For every proper minor-closed graph class $\calG$, for every positive integer $n$,
  there exists a graph $U_n$ with $n^{1+o(1)}$ vertices such that every $n$-vertex graph in $\calG$ is isomorphic to an induced subgraph of $U_n$.
\end{cor}

Prior to the current work, the best result for adjacency labelling of $X$-minor-free graphs (even for $X=K_6$) required labels of length $(2+o(1))\log n$, leading to $n$-induced-universal graphs with $n^{2+o(1)}$ vertices.  (This result is obtained by combining a $(1+o(1))\log n$-bit adjacency labelling scheme for bounded-treewidth graphs due to \citet{gavoille.labourel:shorter} and a result of \citet{devos.ding.ea:excluding} establishing that the edges of every $X$-minor-free graph can be partitioned into two sets, each of which defines a graph of bounded treewidth.)  

The proof of \cref{main_result} makes use of a variant of the Robertson-Seymour Graph Minor Structure Theorem, which describes any graph in a proper minor-closed family in terms of a tree-decomposition.  To make this work, we resolve the problem of ``combining adjacency-labellings along tree-decompositions'' identified by \citet{BGP22} as the main difficulty in obtaining efficient adjacency-labelling schemes for proper minor-closed graph classes.  We do this by introducing and studying a generalization of adjacency labelling schemes called mixed labelling schemes. 

The adjacency labelling scheme we use to prove \cref{main_result} actually satisfies some additional conditions that not only imply \cref{main_result_induced_universal}, but also ensure that the universal graph $U_n$ is in fact sparse:
\begin{thm}\label{main_result_induced_universal_sparse}
    For every proper minor-closed graph class $\calG$, for every positive integer $n$,
    there exists a graph $U_n$ with $n^{1+o(1)}$ vertices and $n^{1+o(1)}$ edges such that every $n$-vertex graph in $\calG$ is isomorphic to an induced subgraph of $U_n$.
\end{thm}

Prior to the current work, the best result for $n$-subgraph-universal graphs for $H$-minor-free graphs was a construction having $\Oh(n^{3/2})$ vertices and edges due to \citet{babai.chung.ea:on}, obtained using a (by now standard) application of $O(\sqrt{n})$ separators.  \Cref{main_result_induced_universal_sparse} improves this quantitatively by reducing the number of vertices and edges to $n^{1+o(1)}$ and qualitatively by constructing induced-universal graphs instead of subgraph-universal graphs.

We establish \cref{main_result_induced_universal_sparse} by introducing a notion of \emph{orientable} (mixed or adjacency) labelling schemes and showing that any orientable adjacency labelling scheme yields an induced-universal graph that is sparse.

\subsection{Relation to Recent Work}

The current paper is the full version of the conference paper \cite{dujmovic.gavoille.ea:adjacency_minor}, which only establishes \cref{main_result,main_result_induced_universal} using mixed labelling schemes.  In order to establish \cref{main_result_induced_universal_sparse}, we define orientable mixed labelling schemes.  The mixed labelling scheme in \cite{dujmovic.gavoille.ea:adjacency_minor} is not orientable, so the mixed labelling scheme presented here is significantly different. In \cite{dujmovic.gavoille.ea:adjacency_minor}, the tree indexing a tree-decomposition is partitioned into subtrees with at most $b$ vertices in each layer---a partition into \emph{$b$-skinny subtrees} whose weights decrease by a factor of $b$ between subtrees.  This then requires a solution for graphs having $b$-skinny tree-decompositions and a solution for graphs having tree-decompositions of height $\Oh(\log_b n)$ whose torsos themselves have $b$-skinny tree-decompositions. The mixed labelling scheme for graphs having $b$-skinny tree-decompositions in \cite{dujmovic.gavoille.ea:adjacency_minor} is not orientable, and there seems to be no easy way to fix this.

The current paper avoids this two-layer approach and thus avoids having to deal specifically with $b$-skinny tree-decompositions. Instead, the proof of \cref{lem:the_big_one} (the technical core of the current paper) performs a partition into subtrees each having $\Oh(b)$ nodes with the property that, within a single subtree, the nodes have weights that differ by a factor of at most $\Oh(b)$.  

Part of our proof requires that we devise an orientable mixed labelling scheme for a generalization of planar graphs. As in \cite{EJM23}, we do this by modifying the adjacency labelling scheme in \cite{dujmovic.esperet.ea:adjacency}, which is based on the Planar Graph Product Structure Theorem  \cite{dujmovic.joret.ea:planar}.  The modifications we use are different from those in \cite{EJM23}, lead to a simplification of the original adjacency labelling scheme from \cite{dujmovic.esperet.ea:adjacency}, and improve the $n^{o(1)}$ factor obtained in \cite{EJM23}.  Since graph products can be used to describe a variety of minor-closed and non-minor-closed graph classes, we state our general result for graph products as a separate theorem.  Let $\calR_k$ denote the class of graphs $G$ for which there exists a graph $H$ of treewidth at most $k$ and a path $P$ such that $G$ is isomorphic to a subgraph of the strong graph product $H\boxtimes P$.

\begin{thm}\label{product_version}
  For all integers $k,n\ge 1$, there exists a graph $U^k_n$ with $n\cdot \exp(\Oh(k\sqrt{\log n}))$ vertices and edges such that every $n$-vertex graph in $\calR_k$ is isomorphic to an induced subgraph of $U^k_n$.
\end{thm}

The proof of \cref{product_version} combines the adjacency labelling scheme of \citet{dujmovic.esperet.ea:adjacency} and its improvement by \citet{gawrychowski.janczewski:simpler} with logarithmic-height binary-tree decompositions discovered by several groups in the late 1990s and early 2000s in the context of parallel algorithms \cite{bodlaender.hagerup:parallel}, reachability queries \cite{chatterjee.ibsen-jensen.ea:optimal}, and space-efficient algorithms \cite{jain.tewari:space}. 

\paragraph{Organization of the Paper}

The remainder of the paper is organized as follows:  
\begin{itemize}
\item \cref{building_blocks} reviews some basic definitions and supporting results including treewidth, tree- and forest-decompositions, strong products, and prefix-free codes. This section also introduces the Graph Minor Product Structure Theorem (\cref{gmpst}) of \citet{dujmovic.joret.ea:planar} and then finishes by developing some tools for manipulating the tree-decompositions that come out of this theorem.
 \item \cref{weighted_mixed_section} begins by introducing weighted mixed labelling schemes and definitions of orientability for both adjacency labelling and mixed labelling schemes.  Mixed labelling schemes are a generalization of adjacency labelling schemes that also satisfy some additional requirements. In order to solve the problem of combining schemes along tree-decompositions, we need these schemes to support vertex-weighted graphs.  In order to obtain sparse induced-universal graphs, we also need these schemes to be orientable.  As a result of all these requirements, the complete definition of an efficient orientable weighted mixed labelling scheme is somewhat involved.  

After these definitions, this section proves a relatively easy result (\cref{add_apexes}) that applies to any graph family $\calG^{+a}$ whose members can be derived from a graph $G$ in a graph family $\calG$ via the addition of at most $a$ apex vertices.  Specifically, it shows that if $\calG$ admits an efficient orientable weighted mixed labelling scheme then $\calG^{+a}$ also admits an efficient orientable weighted mixed labelling scheme.  Although this result and its proof are not difficult, the proof is not particularly short; verifying that all aspects of the definition are satisfied takes some time.  This is a useful warm-up for what comes next.

Next is the main technical core of the paper, \cref{lem:the_big_one}, which establishes that the class $\calG'$ of graphs having tree-decompositions whose torsos come from a graph family $\calG$ that admits an efficient orientable weighted mixed labelling scheme admits an efficient orientable adjacency labelling scheme. Any such class $\calG'$ therefore has $n$-induced-universal graphs with $n^{1+o(1)}$ vertices and edges.  

\item \cref{proof_main_result} combines the ingredients from the previous two sections with an efficient orientable weighted mixed labelling scheme for the class $\calR_k$ of product-structured graphs (\cref{gk_labelling_precise}) in order to prove the main results of the paper, \cref{main_result,main_result_induced_universal,main_result_induced_universal_sparse}.  Briefly: for any proper minor-closed graph family $\calG$, \cref{gmpst} states that each graph in $\calG$ has a tree-decomposition whose torsos belong to $\calR_k^{+a}$, for some fixed $k$ and $a$.  By \cref{gk_labelling_precise,add_apexes}, $\calR_k^{+a}$ admits an efficient orientable weighted mixed labelling scheme.  By \cref{lem:the_big_one}, $\calG$ admits an efficient orientable adjacency labelling scheme.  The resulting adjacency labelling scheme implies \cref{main_result,main_result_induced_universal} and the fact that this scheme is orientable implies \cref{main_result_induced_universal_sparse}.

\item \cref{weighted_product_proof} establishes \cref{gk_labelling_precise}, the existence of an efficient orientable weighted mixed labelling scheme for $\calR_k$.  In addition to its application in proving \cref{main_result,main_result_induced_universal_sparse}, \cref{gk_labelling_precise} immediately implies \cref{product_version}. This result is left until the appendix because it requires a complete review and simplification of the adjacency labelling scheme in \citet{dujmovic.esperet.ea:adjacency} and its improvement by \citet{gawrychowski.janczewski:simpler}.

At the end of this appendix, we point out that this simplification leads to $n$-induced-universal graphs for the class of graphs of treewidth at most $k$ having $n(\log n)^{\Oh(k)}$ vertices and edges.  This improves the previous best known bound (from \cite{EJM23}) of $n(\log n)^{\Oh(k^2)}$.

\item \cref{conclusion} concludes with a brief discussion.
\end{itemize}

\section{Building blocks}\label{building_blocks}

Let \defin{$\N$} be the set of all nonnegative integers, let \defin{$\N^+$} be the set of all positive integers, let \defin{$\R$} be the set of all real numbers, let \defin{$\R^+$} be the set of all positive real numbers, and let \defin{$\R^{\ge 0}$} be the set of all nonnegative real numbers.  For $p\in\N$, let $\mathdefin{[p]} := \set{1,\dots,p}$ (with $[0]:=\emptyset$).  We denote by \defin{$\set{0,1}^*$} the set of all finite binary strings, by \defin{$|\alpha|$} the length (number of bits) of a binary string $\alpha\in\set{0,1}^*$, by \defin{$\varepsilon$} the empty string (the unique binary string of length $0$), and, for $x\in\R^{\ge 0}$, by \defin{$\set{0,1}^{\le x}$} the set of binary strings of length at most $x$.

In this paper, each graph $G$ we consider is a simple graph with vertex set denoted by $V(G)$ and edge set denoted by $E(G)$. For an edge of $G$ with endpoints $v$ and $w$, we use the shorthand notation $vw$ instead of $\{v,w\}$ for the edge, with occasional exceptions.  For a graph $G$ and a vertex $v$ of $G$, $\mathdefin{N_G(v)}:=\{w\in V(G):vw\in E(G)\}$ is the \defin{(open) neighbourhood} of $v$ in $G$; for a subset $X\subseteq V(G)$, $\mathdefin{N_G(X)}:=(\bigcup_{v\in X}N_G(v))\setminus X$.  For any set $S$, \defin{$G[S]$} is the graph with vertex set $\mathdefin{V(G[S])}:=S\cap V(G)$ and edge set $\mathdefin{E(G[S])}:=\{vw\in E(G):v,w\in S\}$, which we call the subgraph of $G$ \defin{induced} by $V(G)\cap S$.  For any set $S$, $G-S:=G[V(G)\setminus S]$; for a single vertex $v$, we write $G-v$ as shorthand for $G-\{v\}$. For any $R\subseteq E(G)$, $G-R$ is the graph with vertex set $V(G)$ and edge set $E(G)\setminus R$. A non-empty subset $K$ of vertices of $G$ is a \defin{clique} in $G$ if $vw\in E(G)$ for each $\{v,w\}\in\binom{K}{2}$.  

A \defin{graph class} is a set of graphs closed under isomorphism. A graph class $\calG$ is \defin{monotone} if, for every $G\in\calG$, every subgraph of $G$ is also in $\calG$. For a graph class $\calG$ and a nonnegative integer $a$, define \defin{$\calG^{+a}$} as the class of graphs $G$ such that $G-X\in\calG$ for some $X\subseteq V(G)$ with $|X|\leq a$ (the vertices in $X$ are sometimes called \defin{apex vertices}). Note that if $\calG$ is monotone, then $\calG^{+a}$ is also monotone.

\subsection{Trees and Forests}

For a rooted tree $T$, we denote the root of $T$ by \defin{$r(T)$}. For a node $x$ of $T$, define \defin{$P_T(x)$} to be the path in $T$ from $r(T)$ to $x$.  A \defin{forest} $F$ is a graph (possibly disconnected) each of whose components is a tree.  If each component of a forest $F$ is a rooted tree, then $F$ is a \defin{rooted} forest. For a node $x$ in a rooted forest $F$, let $\mathdefin{P_F(x)}:=P_T(x)$, where $T$ is the component of $F$ that contains $x$.

Let $F$ be a rooted forest.  For each $x\in V(F)$, the \defin{$F$-depth} of $x$, denoted \defin{$\depth_F(x)$}, is the number of edges in $P_F(x)$. The \defin{height} of $F$, denoted $\mathdefin{\height(F)}:=\max_{x\in V(F)}\depth_F(x)$, is the maximum $F$-depth of a node in $F$. For each $y\in V(F)$ and each $x\in V(P_F(y))$ (including $y$), the node $x$ is an \defin{$F$-ancestor} of $y$ and $y$ is an \defin{$F$-descendant} of $x$.  For each $x\in V(F)$, \defin{$F_x$} denotes the subtree of $F$ induced by all $F$-descendants of $x$. 
For each edge $xy$ of $F$, where $x$ is an $F$-ancestor of $y$, we say that $x$ is the \defin{$F$-parent} of $y$, denoted \defin{$\parent_F(y)$}, and $y$ is an \defin{$F$-child} of $x$. 
We say that a set $S\subseteq V(F)$ is an \defin{$F$-antichain} if there are no distinct $x,y\in S$ such that $x$ is an $F$-ancestor of $y$.  
(When there is no danger of ambiguity, we may drop the leading $F$- from $F$-depth, $F$-ancestor, $F$-descendant, $F$-parent, and $F$-child, and write $\parent(y)$ for $\parent_F(y)$.) We treat each subgraph $F'\subseteq F$ as a rooted forest in which each component $T$ of $F'$ is rooted at the node $r(T)$ of $T$ having minimum $F$-depth.  

For a rooted tree $T$ and any non-empty set $S\subseteq V(T)$, the \defin{lowest common $T$-ancestor} of $S$, denoted \defin{$\lca_T(S)$}, is the node in $\bigcap_{x\in S} V(P_T(x))$ of maximum $T$-depth (for two nodes $x,y\in V(T)$, we also write $\lca_T(x,y)$ for $\lca_T(\{x,y\})$).

\subsection{Forest-Decompositions and Tree-Decompositions}

It will be convenient to work briefly with forest-decompositions, which are the natural generalization of tree-decompositions, allowing the graph indexing the bags to be a forest.  A \defin{forest-decomposition} of a graph $G$ is a pair $\calF:=(F,(B_x\mid x\in V(F)))$, where $F$ is a forest and $(B_x\mid x\in V(F))$ is a collection of subsets of $V(G)$ indexed by the nodes of $F$ such that  
\begin{enumerate*}[label=(\roman*)]
  \item for each edge $vw$ of $G$, there exists $x\in V(F)$ with $v,w\in B_x$; and
  \item for each vertex $v$ of $G$, $F[\{x\in V(F)\mid v\in B_x\}]$ is non-empty and connected.
\end{enumerate*}  
The sets $B_x$ are called the \defin{bags} of $\calF$. When $F$ is a tree, $\calF$ is a \defin{tree-decomposition}. The \defin{width} of a tree-decomposition $(T,(B_x\mid x\in V(T)))$ is $\max_{x\in V(T)}|B_x|-1$. The \defin{treewidth} of a graph $G$, denoted $\mathdefin{\tw(G)}$, is the minimum width of a tree-decomposition of $G$. 

A forest-decomposition is \defin{rooted} if it is indexed by a rooted forest. For each edge $xy$ of $F$, $B_x\cap B_y$ is called an \defin{adhesion} of $\calF$, and we say that this adhesion is incident to $x$ and $y$.  For each $x\in V(F)$, the graph \defin{$\torso{G}{\calF,B_x}$} obtained from $G[B_x]$ by adding edges so that each adhesion incident to $x$ becomes a clique is called the \defin{torso} of $\calF$ at $x$. For a root node $x\in V(F)$, the \defin{parent adhesion} of $x$ is the empty set. For a non-root node $y\in V(F)$, the \defin{parent adhesion} of $y$ is $B_{\parent_F(y)}\cap B_y$. When $\calF$ is rooted, the \defin{home} (or \defin{home node}) of a vertex $v\in V(G)$ is the unique node $x\in V(F)$ such that $v\in B_x$ and $v$ is not in the parent adhesion of $x$.  The \defin{adhesion-width} of $\calF$ is the size of its largest adhesion (or zero if $\calF$ has no edges).

\subsection{Graph Product Structure}

The \defin{strong product} of graphs~$A$ and~$B$, denoted by~${A \boxtimes B}$, is the graph with vertex set ${V(A) \times V(B)}$, where distinct vertices ${(v,x),(w,y) \in V(A) \times V(B)}$ are adjacent if
${v=w}$ and ${xy \in E(B)}$, or
${x=y}$ and ${vw \in E(A)}$, or
${vw \in E(A)}$ and~${xy \in E(B)}$.
For an integer $k\geq 0$, let \defin{$\calR_{k}$} be the class of graphs isomorphic to a subgraph of $H\boxtimes P$ for some graph $H$ of treewidth at most $k$ and some path $P$.  
Our starting point is the following variant of the Graph Minor Structure Theorem:

\begin{thm}[\protect Graph Minor Product Structure Theorem~\citep{dujmovic.joret.ea:planar}]
  \label{gmpst}
  For every proper minor-closed class $\calG$, there exist integers $k,a\geq 0$ such that every graph in $\calG$ has a tree-decomposition in which every torso is in $\calR_k^{+a}$.
\end{thm}

\subsection{Weight Functions}

For a non-empty set $S$, a \defin{pseudoweight function} over $S$ is a function $\omega:S\to\R^{\ge 0}$.  A \defin{weight function} is a pseudoweight function that assigns strictly positive weights (i.e., $\omega(x)>0$ for each $x\in S$).  For $X\subseteq S$, we use the shorthand $\mathdefin{\omega(X)}:=\sum_{x\in X}\omega(x)$.  When $S:=V(G)$ is the vertex-set of a graph $G$ and $H$ is a subgraph of $G$, we use the shorthand $\mathdefin{\omega(H)}:=\omega(V(H))$.

Let $\omega:S\to\R^+$ be a weight function.  For each $x\in S$, the \defin{(Kraft) ideal codeword length} for $x$ is $\log\omega(S)-\log\omega(x)$.
It will be convenient to bound the total weight as well as the Kraft ideal codeword lengths. A weight function $\omega:S\to\N^+$ is \defin{nice} if $\omega(S)\le 2|S|$. In particular, a nice weight function satisfies
\[
  \log\omega(S)-\log\omega(x)\le\log|S|+1
\]
for every $x\in S$.

\begin{obs}\label{obs:nice}
  For every weight function $\omega:S\to\R^+$, there exists a nice weight function $\omega':S\to\N^+$ such that $\log\omega'(S)-\log\omega'(x)\le\min\{\log|S|,\;\log\omega(S)-\log\omega(x)\}+1$, for each $x\in S$. If $\omega$ is positive integer-valued, then $\omega'(x)\le\omega(x)$ for each $x\in S$, and hence $\omega'(S)\le\omega(S)$.
\end{obs}

\begin{proof}
  Let $m:=|S|$ and $W:=\omega(S)$. For each $x\in S$, define
  \[
    \omega'(x):=\left\lceil\frac{m\omega(x)}{W}\right\rceil\enspace.
  \]
  Then $\omega'(x)\ge 1$ for each $x\in S$ and
  \[
    m\le\omega'(S)<\sum_{x\in S}\left(\frac{m\omega(x)}{W}+1\right)=2m\enspace,
  \]
  so $\omega'$ is nice. Since $\omega'(x)\ge 1$ and $\omega'(x)\ge m\omega(x)/W$, we obtain
  \[
    \frac{\omega'(S)}{\omega'(x)}\le\min\left\{2m,\;\frac{2W}{\omega(x)}\right\}\enspace.
  \]
  Taking logarithms proves the claimed bound. Finally, if $\omega$ is positive integer-valued, then $W\ge m$, so $m\omega(x)/W\le\omega(x)$. Since $\omega(x)$ is an integer, $\omega'(x)=\ceil{m\omega(x)/W}\le\omega(x)$ for each $x\in S$.
\end{proof}

\subsection{On Multipart Labels}
\label{multipart}

In many cases, a label (a bitstring) consists of a sequence of parts (also bitstrings) of varying lengths, $s_1,\ldots,s_p$, where each $s_i$ may be empty.  The naive concatenation of $s_1,\ldots,s_p$ has length
\[
  \textstyle \sum_{i=1}^p|s_i| \enspace .
\]
This ignores the important issue that, in order for a decoder to extract the individual parts $s_1,\ldots,s_p$, it must know the number of parts $p$ and (at least $p-1$ of) the lengths $|s_1|,\ldots,|s_p|$.  We handle this by encoding the value of $p$ and the lengths $|s_1|,\ldots,|s_p|$ and prepending these to the concatenation. This is most easily done using Elias' $\gamma$-code $\gamma:\N^+\to\{0,1\}^*$ for positive integers, in which each $x\in\N^+$ is encoded by a self-delimiting binary string $\gamma(x)$ of length $|\gamma(x)|= 2\lfloor \log x  \rfloor  + 1$ \cite{elias:universal}.\footnote{For each $x\in\N^+$, $\gamma(x)$ consists of $\lfloor\log x\rfloor$ 0 bits, followed by the binary representation of $x$ (which begins with a $1$ bit and has length $\lfloor\log x\rfloor+1$).} In our application, each part $s_i$ may have zero length and $p$ may even be equal to zero, so we encode each $x\in\N$ using $\gamma(x+1)$, which has length $2\floor{\log (x+1)}+1$.
In this way, we obtain a self-delimiting encoding of $s_1,\ldots,s_p$ as
\[  \mathdefin{\langle s_1,\ldots,s_p\rangle}:=\gamma(p+1)\circ\gamma(|s_1|+1)\circ\cdots\circ\gamma(|s_p|+1)\circ s_1\circ\cdots\circ s_p \enspace ,
\]
where $\circ$ denotes string concatenation (and $\langle\rangle:=\gamma(1)$ when $p=0$). This encoding has length
\begin{equation}
  |\langle s_1,\ldots,s_p\rangle| = \textstyle\sum_{i=1}^p|s_i|+\Oh(\lgg p) + \Oh(\sum_{i=1}^p\lgg|s_i|) \enspace , \label{encoding_cost}
\end{equation}
where $\mathdefin{\lgg x} := \log\max\{2,x\}$ (so that $\lgg x$ is defined and $\lgg x\ge 1$ for each $x\in\N$).  The bits of $\langle s_1,\ldots,s_p\rangle$ that come from $s_1,\ldots,s_p$ are called the \defin{payload} bits and the remaining bits are called the \defin{bookkeeping} bits of $\langle s_1,\ldots,s_p\rangle$.
In our setting it will always be the case that $p\in o(\log n)$ and we will use \cref{obs:nice} to ensure that $|s_i|\in \Oh(\log n)$ for each $i\in[p]$.
Thus, the number of bookkeeping bits in $\langle s_1,\ldots,s_p\rangle$ is $\Oh(1 + p\log\log n)$.

In some cases, we include one or more non-negative integers as parts of a multipart label.  When we include $d\in\N$ as part of a label, we include its $\ceil{\log(d+1)}$-bit binary representation \defin{$\bin(d)$} (so $\bin(0)=\varepsilon$).

\subsection{Shannon and Shannon--Fano--Elias Codes}

A rooted tree $T$ is a \defin{binary tree} if each node of $T$ has at most two children among which at most one is the \defin{left child} and at most one is the \defin{right child}.
A binary tree $T$ is \defin{full} if for every node $x$ of $T$,
either $x$ has two children or $x$ has no children.
A full binary tree $T$ is \defin{complete} if all the leaves of $T$ have the same depth in $T$.
Let $T$ be a binary tree. There is a natural ordering, called \defin{left-to-right ordering}, of the leaves of $T$, such that for every node $x$, if $x$ has two children, then all the leaf-descendants of the left child of $x$ precede all the leaf-descendants of the right child of $x$.

The following result is a variant of Alphabetic Binary Trees. Similar results (with the constant $3$ replaced by $2$) are well-known and there are several proofs \cite{knuth:art_iii,mehlhorn:nearly,gilbert.moore:variable_length}. We provide a proof of the weaker statement here because it illustrates a simple but powerful technique---\defin{simulating weight by multiplicity}---that we will use again in \cref{weighted_product_proof}.

\begin{lem}\label{shannon_fano_elias_tree}
  Let $S$ be a non-empty linearly ordered finite set and let $\omega:S\to\N^+$ be a weight function. Then there exists a full binary tree $T$ such that $S$ is the set of leaves of $T$,
  the linear order on $S$ agrees with the left-to-right order of the leaves in $T$, and for each $x\in S$,
  \[
    \depth_{T}(x)\le \log\omega(S) - \log\omega(x) + 3 \enspace.
  \]
\end{lem}

\begin{proof}
  Let $h:=\lceil\log\omega(S)\rceil$.
  Let $Q$ be the complete binary tree of height $h$, so with $2^h$ leaves.
  For each $x\in S$, define the $\omega(x)$-element set
  \[  S_{\omega}(x):=\set{(x,1),\,(x,2),\,\ldots\,,\,(x,\omega(x))}\]
  and let $S_\omega:=\bigcup_{x\in S}S_{\omega}(x)$.  Then $|S_\omega|=\omega(S)$.  Treat $S_{\omega}$ as a linearly ordered set by using lexicographic order and rename the leftmost $\omega(S)$ leaves of $Q$, in order, with elements of $S_{\omega}$.

  Now for each $x\in S$, we will fix a node $r(x)$ in $Q$ such that
  (1) all the leaves of $Q_{r(x)}$ lie in $S_{\omega}(x)$; and
  (2) $Q_{r(x)}$ contains at least $\omega(x)/4$ leaves.

  Let $x\in S$.
  If $\omega(x)=1$ then define $r(x):=(x,1)$, which is a leaf of $Q$.
  Clearly, (1) and (2) hold.
  Now suppose that $\omega(x)\geq2$.
  Consider the lowest common $Q$-ancestor $y$ of $S_{\omega}(x)$.
  Thus, $y$ is an internal node of $Q$ and each child of $y$ has a leaf-descendant in $S_{\omega}(x)$.
  Since $S_{\omega}(x)$ is a set of consecutive leaves in $Q$,
  all leaves of the left child $y'$ of $y$ in $Q$ that are in $S_{\omega}(x)$ are rightmost in $Q_{y'}$ and
  all leaves of the right child $y''$ of $y$ in $Q$ that are in $S_{\omega}(x)$ are leftmost in $Q_{y''}$.
  Let $z$ be the child of $y$ in $Q$ such that $Q_z$ has at least $\omega(x)/2$ leaves in $S_{\omega}(x)$ and let $S_z$ be the set of all such leaves.
  Thus, $|S_z|\geq \omega(x)/2$.
  Since the argument is completely symmetric in these cases, assume that $z$ is the left child of $y$.
  Consider a path in $Q_z$ starting from $z$ and always taking an edge to the right child of the current node until we finish at a leaf.
  Let $v$ be the first vertex of that path such that $Q_v$ has all leaves in $S_z$.
  If $v=z$ then we set $r(x):=v$. Clearly (1) and (2) hold.
  Suppose that $v\neq z$.
  Let $u$ be the parent of $v$.
  Since the set of leaves of $Q_u$ form an interval of rightmost leaves of $Q_z$ and since $Q_u$ contains a leaf not in $S_z$, we conclude that $Q_u$ contains all nodes in $S_z$ as leaves.
  Therefore, $Q_v$ contains at least $|S_z|/2$ vertices of $S_z$ as leaves.
  We set $r(x):=v$ and again (1) and (2) hold.

  Note that the set $\set{r(x) \mid x\in S}$ is a $Q$-antichain.
  Indeed, if $x,y\in S$, $x\neq y$ and $r(x)$ is a $Q$-ancestor of  $r(y)$, then the leaves of $Q_{r(y)}$ are contained in the leaves of $Q_{r(x)}$ which by (1) implies that $Q_{r(x)}$ contains all nodes in $S_{\omega}(y)$ as leaves.
  However, again by (1) all the leaves of $Q_{r(x)}$ lie in $S_{\omega}(x)$ which is disjoint from $S_{\omega}(y)$, a contradiction.

  Let $x\in S$. Since $Q_{r(x)}$ contains at least $\omega(x)/4$ leaves, $\height(Q_{r(x)}) \geq \log(\omega(x)/4) = \log\omega(x) - 2$. Finally,
  \begin{align*}
    \depth_Q(r(x))
     & \le \height(Q)-\height(Q_{r(x)})                    \\
     & \le \lceil\log\omega(S)\rceil - (\log\omega(x) - 2) \\
     & \le \log\omega(S) - \log\omega(x) +3 \enspace.
  \end{align*}

  Let $T_0:=\bigcup_{x\in S} P_Q(r(x))$.
  Since $\set{r(x) \mid x\in S}$ is a $Q$-antichain,
  we conclude that $\set{r(x)\mid x\in S}$ is the set of all leaves of $T_0$.
  Obtain $T$ from $T_0$ by repeatedly contracting an edge incident to any internal node that has only one child.
  Then $T$ is a full binary tree with leaf set $\set{r(x)\mid x\in S}$, and $\depth_T(r(x))\le\depth_{T_0}(r(x))=\depth_Q(r(x))$ for each $x\in S$.
  Identifying each leaf $r(x)$ of $T$ with $x$ for each $x\in S$ gives the desired tree.
\end{proof}

A \defin{code} for a set $S$ is an injective function $\rho:S\to\{0,1\}^*$. A code $\rho:S\to\{0,1\}^*$ is \defin{prefix-free} if for all distinct $x,y\in S$, $\rho(x)$ is not a prefix of $\rho(y)$.
In coding theory, \cref{shannon_fano_elias_tree} is the basis of the Shannon--Fano--Elias coding scheme, as we now explain.  By assigning each edge $e$ of a binary tree $T$ a $0$ or $1$ depending on whether $e$ joins a parent to its left or right child, we obtain an encoding of each root-to-leaf path in $T$ as a binary string, which gives a prefix-free code for the leaves of $T$. Doing this with the tree from \cref{shannon_fano_elias_tree}, we obtain the following.

\begin{cor}\label{shannon_fano_elias_code}
  There exists a function $F:(\{0,1\}^*)^2\to\{0,1\}^*$ such that, for every $n\in\N^+$, every set $S:=\{x_1,\ldots,x_n\}$, and every weight function $\omega:S\to\R^+$, there exists a prefix-free code $\rho:S\to \{0,1\}^*$ for $S$ and functions $\overleftarrow{\pi}:S\to\{0,1\}^*$, $\overrightarrow{\pi}:S\to\{0,1\}^*$ such that:
  \begin{enumerate}[nosep,nolistsep,label=\rm(\alph*)]
    \item\label{sfec_len} $|\rho(x_i)|\le \log\omega(S) - \log\omega(x_i)+\Oh(1)$ for each $i\in[n]$;
    \item\label{sfec_size} $|\overleftarrow{\pi}(x_i)|+|\overrightarrow{\pi}(x_i)|\in O(\log\log n)$ for each $i\in[n]$;
    \item\label{sfec_prev} $F(\rho(x_i),\overleftarrow{\pi}(x_i))=\rho(x_{i-1})$ for each $i\in\{2,\ldots,n\}$; and    
    \item\label{sfec_next} $F(\rho(x_i),\overrightarrow{\pi}(x_i))=\rho(x_{i+1})$ for each $i\in\{1,\ldots,n-1\}$.    
  \end{enumerate}
\end{cor}

\begin{proof}
  Let $\omega':S\to\N^+$ be the nice weight function obtained by applying \cref{obs:nice} to $\omega$.  Then $\log\omega'(S)\le\log n+ O(1)$ and $\log\omega'(x_i)\ge 0$ for each $i\in[n]$.
  The prefix-free code $\rho$ is obtained by applying \cref{shannon_fano_elias_tree} to $S$ (ordered by $x_1<\cdots<x_n$) with weight function $\omega'$ to obtain a full binary tree $T$ and then encoding each root-to-leaf path in $T$ with a bitstring as described above.  This gives a code $\rho$ that satisfies \ref{sfec_len}.
  
  Set $\overleftarrow{\pi}(x_1):=\varepsilon$. Fix $i\in\{2,\ldots,n\}$ and let $y_0,y_1,\ldots,y_\ell$ be the path from $y_0:=x_i$ to $y_\ell:=x_{i-1}$ in $T$.  Let $j$ be such that $y_j=\lca_T(x_i,x_{i-1})$.  Since $T$ is a full binary tree, $y_{j+1}$ is the left child of $y_j$ and $y_k$ is the right child of $y_{k-1}$ for each $k\in\{j+2,\ldots,\ell\}$.  Thus, $\rho(x_{i-1})$ is the prefix of $\rho(x_i)$ of length $\depth_T(y_j)$, followed by a $0$ bit, followed by a string of $1$ bits of length $\ell-j-1$.  We set $\overleftarrow{\pi}(x_i):=\langle 0, \bin(j), \bin(\ell-j-1)\rangle$.  Since $j <\ell$ and $\ell\le 2\height(T)\in O(\log n)$, $\overleftarrow{\pi}(x_i)$ has length $O(\log\log n)$.  Symmetrically, set $\overrightarrow{\pi}(x_n):=\varepsilon$ and, for each $i\in\{1,\ldots,n-1\}$, define $\overrightarrow{\pi}(x_i):=\langle 1, \bin(j), \bin(\ell-j-1)\rangle$ using the path $y_0,y_1,\ldots,y_\ell$ from $y_0:=x_i$ to $y_\ell:=x_{i+1}$ in $T$, where $y_j=\lca_T(x_i,x_{i+1})$.  

  The operation of $F$ is now obvious.  Given the arguments $\rho(x_i)$ and $\langle a,\bin(b),\bin(c)\rangle$, $F$ checks the value of $a\in\{0,1\}$.  If $a=0$, then it removes the suffix of length $b$ from $\rho(x_i)$ and appends the string $01^c$ to obtain $\rho(x_{i-1})$.  If $a=1$, then $F$ removes the suffix of length $b$ from $\rho(x_i)$ and appends the string $10^c$ to obtain $\rho(x_{i+1})$. 
\end{proof}

\subsection{Working with Tree Decompositions}
\label{working_with_tree_decomps}

In this section, we develop some tools for manipulating tree-decompositions like those obtained from \cref{gmpst}.

\begin{lem}\label{skinny_edges}
  Let $b>1$, let $T$ be a rooted tree, and let $\omega:V(T)\to\N^+$ be a weight function.  Then there exists $Y\subseteq E(T)$ such that
  \begin{enumerate}[nosep,nolistsep,label=\rm(\alph*)] 
    \item each root-to-leaf path in $T$ contains at most $\log_b \omega(T)$ edges in $Y$; and
    \item $\omega(T_y)\le b\cdot\omega(T_x)$ for any two vertices $x,y$ in the same component of $T-Y$.
  \end{enumerate}
\end{lem}

\begin{proof}
  The proof is by induction on $|V(T)|$.  Let $X_0:=\{x\in V(T):\omega(T_x)\ge \omega(T)/b\}$.  Since the root $r$ of $T$ satisfies $\omega(T_r)=\omega(T)>\omega(T)/b$ (as $b>1$) and every ancestor of a vertex in $X_0$ is also in $X_0$, the induced subgraph $T_0:=T[X_0]$ is a connected subtree of $T$ containing $r$.  If $X_0=V(T)$, then setting $Y:=\emptyset$ satisfies the lemma: each root-to-leaf path contains $0 \le \log_b\omega(T)$ edges of $Y$ (since $\omega(T)\ge 1$), and for any $x,y\in V(T)$, $\omega(T_y)\le \omega(T)\le b\cdot\omega(T_x)$ by the definition of $X_0$.

  Otherwise, let $F_0:=T-X_0$ and let $T_1,\ldots,T_s$ be the connected components of $F_0$.  By the definition of $X_0$, for any $x,y\in V(T_0)$ we have $\omega(T_y)\le\omega(T)\le b\cdot\omega(T_x)$.  Furthermore, for each $i\in[s]$, the component $T_i$ is the subtree $T_{y_i}$ rooted at the child $y_i\notin X_0$ of some vertex $x_i\in V(T_0)$, so $\omega(T_i)<\omega(T)/b$.  Since $r\in X_0$, each component $T_i$ satisfies $|V(T_i)|<|V(T)|$.  Applying induction to each $T_i$ with weight function $\omega|_{V(T_i)}$ yields an edge set $Y_i\subseteq E(T_i)$.

  Let $Y_0:=\{x_iy_i\mid i\in[s]\}$ and let $Y:=Y_0\cup\bigcup_{i=1}^s Y_i$.  Each component $C$ of $T-Y$ is either $T_0$ or a component of $T_i-Y_i$ for some $i\in[s]$.  In the former case, $\omega(T_y)\le b\cdot\omega(T_x)$ for all $x,y\in V(T_0)$ as established above.  In the latter case, since $(T_i)_z = T_z$ for every $z\in V(T_i)$, the inductive hypothesis implies $\omega(T_y)\le b\cdot\omega(T_x)$ for all $x,y\in V(C)$.  This establishes condition~(b).

  Finally, let $P$ be a root-to-leaf path in $T$.  If $P$ is contained entirely in $T_0$, then $E(P)\cap Y = \emptyset$, so $|E(P)\cap Y| = 0 \le \log_b\omega(T)$.  Otherwise, $P$ leaves $T_0$ via a unique edge $x_iy_i\in Y_0$ into some component $T_i$, and the suffix $P_i := P\cap T_i$ is a root-to-leaf path in $T_i$.  By induction, $|E(P_i)\cap Y_i|\le\log_b\omega(T_i)$.  Therefore,
  \[  
    |E(P)\cap Y|
      = |E(P)\cap Y_0| + |E(P_i)\cap Y_i|
      \le 1+\log_b\omega(T_i)
      < 1+\log_b(\omega(T)/b)
      = \log_b\omega(T) \enspace . \qedhere 
  \]
\end{proof}

The following lemma is a standard weighted separator theorem for trees; we include its short proof for completeness.

\begin{lem}\label{bag_separator}
  Let $T$ be a tree and let $\omega:V(T)\to\R^{\ge 0}$ be a pseudoweight function.  Then, for any $b\in\N^+$, there exists $S\subseteq V(T)$ with $|S|\le b-1$ such that $\omega(C)\le \omega(T)/b$ for each component $C$ of $T-S$.
\end{lem}

\begin{proof}
  The proof is by induction on $b$.  If $b=1$, then $S:=\emptyset$ satisfies the requirements of the lemma since $T$ is the only component of $T-S$.  Now assume $b\ge 2$. Root $T$ at an arbitrary vertex. The set $\{v\in V(T):\omega(T_v)\ge \omega(T)/b\}$ contains the root of $T$ and is therefore non-empty. Let $x$ be a vertex in this set of maximum $T$-depth. Then $\omega(T_x)\ge\omega(T)/b$, while $\omega(T_y)<\omega(T)/b$ for each child $y$ of $x$ (if $\omega(T)=0$, $x$ is a leaf and this holds vacuously).

  The components of $T-x$ consist of the subtrees $T_y$ for each child $y$ of $x$ (each having weight $\omega(T_y)<\omega(T)/b$), and possibly the component $T_0 := T-V(T_x)$ containing the parent $\parent(x)$ of $x$. If $\omega(T_0)\le \omega(T)/b$ (or if $T_0$ does not exist), then $S:=\{x\}$ satisfies the lemma since $|S|=1\le b-1$ and every component of $T-x$ has weight at most $\omega(T)/b$.

  Otherwise, $\omega(T_0)>\omega(T)/b$. Since $\omega(T_x)\ge\omega(T)/b$, we have $\omega(T_0)=\omega(T)-\omega(T_x)\le (1-1/b)\omega(T) = \frac{b-1}{b}\omega(T)$. Applying induction to the tree $T_0$ with pseudoweight function $\omega|_{V(T_0)}$ and parameter $b-1$ yields a set $S_0\subseteq V(T_0)$ of size at most $b-2$ such that each component $C$ of $T_0-S_0$ satisfies $\omega(C)\le \frac{\omega(T_0)}{b-1}\le \frac{\omega(T)}{b}$.

  We define $S:=S_0\cup\{x\}$, so $|S|\le (b-2)+1 = b-1$. Each component of $T-S$ is either a component of $T_0-S_0$ (which has weight at most $\omega(T)/b$ by the inductive hypothesis) or one of the subtrees $T_y$ for a child $y$ of $x$ (which has weight $\omega(T_y)<\omega(T)/b$). Thus, $S$ satisfies the requirements of the lemma.
\end{proof}

\begin{lem}\label{log_b_td}
  Let $G$ be an $n$-vertex graph, with $n\in\N^+$, and let $\calT$ be a tree-decomposition of $G$.  Then, for any integer $b\ge 2$, there exists a rooted tree-decomposition $\calQ:=(Q,(C_x\mid x\in V(Q)))$ of $G$ and $Y\subseteq E(Q)$ with the following properties:
  \begin{enumerate}[nosep,nolistsep,label=\rm(\alph*)]
    \item\label{req:height_0} each root-to-leaf path in $Q$ contains at most $\log_b n$ edges of $Y$;
    \item\label{req:subtree_size_0} each connected component of $Q-Y$ has at most $4b-1$ vertices;
    \item\label{req:adhesion_width_0} each adhesion of $\calQ$ is contained in the union of at most two adhesions of $\calT$; and
    \item\label{req:subtree_torso_apex_0} for each $x\in V(Q)$, at least one of the following holds:
    \begin{itemize}[nosep,nolistsep]
      \item\label{req:three_adhesions_0} $C_x$ is contained in the union of at most three adhesions of $\calT$; or
      \item\label{req:adhesion_torso_0} $\torso{G}{\calQ,C_x}-A$ is a subgraph of some torso of $\calT$, where $A=\emptyset$ or $A$ is an adhesion of $\calT$.
    \end{itemize}
  \end{enumerate}
\end{lem}

\begin{proof}
  The proof is by induction, with a generalized inductive hypothesis. An inductive instance $(G,\calT, \alpha, \beta)$ involves a graph $G$, a tree-decomposition $\calT:=(T,(B_x\mid x\in V(T)))$ of $G$, and two (not necessarily distinct) leaves $\alpha,\beta$ of $T$ such that each of $B_\alpha$ and $B_\beta$ equals the adhesion of $\calT$ incident to $\alpha$ and $\beta$ respectively, and $|V(G)\setminus (B_\alpha\cup B_\beta)|\ge 1$.  We will produce a rooted tree-decomposition $\calQ:=(Q,(C_x\mid x\in V(Q)))$ of $G$ and $Y\subseteq E(Q)$ that satisfies all the requirements of the lemma as well as the following conditions, the first of which is a strengthening of \ref{req:height_0}:
  \begin{enumerate}[nosep,nolistsep,label=\rm(\alph*)]\setcounter{enumi}{4}
    \item\label{req:height_prime_0} each root-to-leaf path in $Q$ contains at most $\log_b(n-|B_\alpha\cup B_\beta|)$ edges of $Y$;
    \item\label{req:root_0} $\alpha\in V(Q)$ and $\alpha$ is the root of $Q$; and
    \item\label{req:adhesion_pairs_0} $C_\alpha=B_\alpha\cup B_\beta$.
  \end{enumerate}
  The induction is on the quantity $n-|B_\alpha\cup B_\beta| = |V(G)\setminus (B_\alpha\cup B_\beta)|$.  To see that this suffices to establish the lemma for an $n$-vertex graph $G$, we can attach two new leaves $\alpha$ and $\beta$ to $T$ to obtain $T'$, set $B_\alpha:=B_\beta:=\emptyset$ (so their incident adhesions are empty), and apply the refined result to $(G,(T',(B_x\mid x\in V(T'))),\alpha,\beta)$ to obtain the desired rooted tree-decomposition $\calQ$ and edge set $Y$.

  We now proceed with the inductive proof.  Root $T$ at $\alpha$ so that $\calT$ is a rooted tree-decomposition of $G$.  For each $x\in V(T)\setminus\{\alpha\}$, let $A_x$ be the parent adhesion of $x$ in $\calT$, and set $A_\alpha:=\emptyset$.  Define the pseudoweight function $\omega: V(T)\to\N$ by $\omega(x):=|B_x\setminus (A_x\cup B_\alpha\cup B_\beta)|$.  Since each vertex in $V(G)\setminus (B_\alpha\cup B_\beta)$ belongs to $B_x\setminus A_x$ for exactly one node $x\in V(T)$, we have $\omega(T) = n-|B_\alpha\cup B_\beta|$.

  Apply \cref{bag_separator} to the tree $T$ with pseudoweight function $\omega$ to obtain a set $S_0$ of at most $b-1$ nodes of $T$ such that $\omega(C)\le\omega(T)/b=(n-|B_\alpha\cup B_\beta|)/b$ for each component $C$ of $T-S_0$.
  Let $T'$ be the vertex-minimal subtree of $T$ that spans $S_0\cup\{\alpha,\beta\}$.  Let $S_1$ be the set of nodes of $T'$ that have degree at least $3$ in $T'$.  Then $|S_1|\le |S_0\cup\{\alpha,\beta\}|-2\le b-1$.  Let $S:=S_0\cup S_1\cup\{\alpha,\beta\}$, so $|S|\le 2b$.  Let $C_1,\ldots,C_s$ be the connected components of $T-S$.  For each $i\in[s]$, let $N_i:=N_T(V(C_i))$.  Since every branching node of $T'$ belongs to $S_1\subseteq S$, no component $C_i$ can have three or more neighbors in $S$, so $|N_i|\in\{1,2\}$ for each $i\in[s]$.  For each $i\in[s]$, write $N_i:=\{a_i,b_i\}$ (where $a_i=b_i$ if $|N_i|=1$).

  For each $i\in[s]$, let $G_i:=G[\bigcup_{x\in V(C_i)}B_x]$ and let $m_i:=|V(G_i)\setminus (B_{a_i}\cup B_{b_i})|$.  Observe that  $m_i\le\omega(C_i)\le\omega(T)/b$ for each $i\in [s]$. For each $i\in[s]$, let $\{x_i\}:=N_T(a_i)\cap V(C_i)$ and  $\{y_i\}:=N_T(b_i)\cap V(C_i)$ so that $x_i$ and $y_i$ are the unique neighbours in $C_i$ of $a_i$ and $b_i$, respectively (possibly $x_i=y_i$). Since $\calT$ is a tree-decomposition,
    $B_{a_i}\cap V(G_i)=B_{a_i}\cap B_{x_i}$ and $B_{b_i}\cap V(G_i)=B_{b_i}\cap B_{y_i}$.

  We construct a tree $Q_0$ as follows: $V(Q_0)$ contains all nodes of $S$.  For each $i\in[s]$ such that $|N_i|=2$, we add a node $v_i$ and edges $a_iv_i$ and $b_iv_i$ to $Q_0$. 
  For each edge $xy\in E(T)$ with $x,y\in S$, we add the edge $xy$ to $Q_0$.  Since each component $C_i$ with $|N_i|=2$ corresponds to a path in $T'$ between two nodes in $S$, $Q_0$ is a tree with $V(Q_0)\supseteq S$ and $|V(Q_0)|\le 2|S|-1\le 4b-1$.  We root $Q_0$ at $\alpha$.

  For each $x\in S$, set $C_x:=B_x$.  For each $i\in[s]$ with $|N_i|=2$, set $C_{v_i}:=(B_{a_i}\cup B_{b_i})\cap V(G_i)$.  To ensure property \ref{req:adhesion_pairs_0}, we add the vertices in $B_\beta$ to $C_x$ for each $x$ on the path from $\beta$ to $\alpha$ in $Q_0$.  This completes the definition of the bags $C_x$ for all $x\in V(Q_0)$.  Then $\calQ_0:=(Q_0,(C_x\mid x\in V(Q_0)))$ is a tree-decomposition of $G_S:=G[\bigcup_{x\in S}B_x]$.

  If $m_i=0$, then $V(G_i)\subseteq B_{a_i}$ when $|N_i|=1$, and $V(G_i)\subseteq C_{v_i}$ when $|N_i|=2$. Thus all vertices and edges of $G_i$ are covered by a single bag of $\calQ_0$, and no recursive call is needed for this component.  For each $i\in[s]$ such that $m_i>0$, we construct a tree-decomposition of $G_i$ as follows. 

  \begin{itemize}
    \item If $|N_i|=1$ (so $N_i=\{a_i\}$) then $x_i\in V(C_i)$ is adjacent to $a_i$ in $T$. Let $T_i$ be the tree obtained from $C_i$ by adding a new leaf $\alpha_i=\beta_i$ adjacent to $x_i$, and set $B_{\alpha_i}:=B_{a_i}\cap B_{x_i}$.  Then $\calT_i:=(T_i, (B_x\mid x\in V(T_i)))$ is a tree-decomposition of $G_i$.  Apply induction to the instance $(G_i,\calT_i,\alpha_i,\beta_i)$.
    \item If $|N_i|=2$ then $x_i,y_i\in V(C_i)$ are the unique neighbours of $a_i$ and $b_i$ in $C_i$ (possibly $x_i=y_i$). Let $T_i$ be the tree obtained from $C_i$ by adding new leaves $\alpha_i$ adjacent to $x_i$ and $\beta_i$ adjacent to $y_i$, and set $B_{\alpha_i}:=B_{a_i}\cap B_{x_i}$ and $B_{\beta_i}:=B_{b_i}\cap B_{y_i}$.  Then $\calT_i:=(T_i, (B_x\mid x\in V(T_i)))$ is a tree-decomposition of $G_i$.  Apply induction to the instance $(G_i,\calT_i,\alpha_i,\beta_i)$.
  \end{itemize}
  In either case, induction yields a rooted tree-decomposition $\calQ_i:=(Q_i,(C_x\mid x\in V(Q_i)))$ of $G_i$ with root $\alpha_i$, and an edge set $Y_i\subseteq E(Q_i)$.  By \ref{req:adhesion_pairs_0}, the root bag satisfies $C_{\alpha_i}=B_{\alpha_i}\cup B_{\beta_i}$.

  We now construct the tree $Q$ from $Q_0$ and $\{Q_i\mid m_i\ge 1\}$.   
  For each $i\in[s]$ such that $m_i>0$ and  $|N_i|=1$, let $e_i:=a_i\alpha_i$.
  For each $i\in[s]$ such that $m_i>0$ and $|N_i|=2$, let $e_i:=v_i\alpha_i$.
  Then $V(Q):=V(Q_0)\cup \bigcup_{i: m_i\ge 1} V(Q_i)$ and $E(Q):=E(Q_0)\cup \bigcup_{i: m_i\ge 1} (E(Q_i) \cup \{e_i\})$.  Define $Y := \bigcup_{i: m_i\ge 1} (Y_i \cup \{e_i\})$.  Since $V(G_S)\cap V(G_i)=B_{\alpha_i}\cup B_{\beta_i}=C_{\alpha_i}$ and $C_{\alpha_i}\subseteq C_{a_i}$ (when $|N_i|=1$) or $C_{\alpha_i}\subseteq C_{v_i}$ (when $|N_i|=2$), $\calQ:=(Q,(C_x\mid x\in V(Q)))$ is a valid rooted tree-decomposition of $G$ with root $\alpha$, satisfying \ref{req:root_0} and \ref{req:adhesion_pairs_0}.

  Because each edge $e_i$ belongs to $Y$, the connected components of $Q-Y$ are precisely $Q_0$ and the connected components of $Q_i-Y_i$ for each $i$ with $m_i\ge 1$.  Since $|V(Q_0)|\le 4b-1$ and each component of $Q_i-Y_i$ has at most $4b-1$ vertices by induction, property \ref{req:subtree_size_0} holds.

  Any root-to-leaf path in $Q$ consists of a path in $Q_0$ followed by at most one edge $e_i$ and a root-to-leaf path in $Q_i$ for some $i$ with $m_i\ge 1$.  A path contained in $Q_0$ uses no edges of $Y$. Otherwise, since $E(Q_0)\cap Y=\emptyset$, the number of edges of $Y$ on the path is at most:
  \[
    1 + \log_b(|V(G_i)| - |B_{\alpha_i}\cup B_{\beta_i}|) \le 1 + \log_b\omega(C_i) \le 1 + \log_b((n-|B_\alpha\cup B_\beta|)/b) = \log_b(n-|B_\alpha\cup B_\beta|) \enspace .
  \]
  Thus, property \ref{req:height_prime_0} holds.

  Next, we verify property \ref{req:adhesion_width_0}. For an edge $xy\in E(Q_0)$ with $x,y\in S$, the adhesion $C_x\cap C_y$ is contained in $(B_x\cap B_y)\cup B_\beta$, which is the union of at most two adhesions of $\calT$. For an edge $a_i v_i$ (or $b_i v_i$), the adhesion is contained in $(B_{a_i}\cap B_{x_i})\cup B_\beta$, also the union of at most two adhesions of $\calT$. For each edge $e_i$, the adhesion is contained in $C_{\alpha_i} = B_{\alpha_i}\cup B_{\beta_i}$, which is the union of two adhesions of $\calT$. Finally, for an edge of $Q_i$, the adhesion is contained in the union of at most two adhesions of $\calT_i$ by induction, and each adhesion of $\calT_i$ is an adhesion of $\calT$. Thus, property~\ref{req:adhesion_width_0} holds.

  Finally, we verify property \ref{req:subtree_torso_apex_0}.  First we consider the torsos of $\calQ$ at vertices $x\in S$.  Let $x$ be a vertex in $S$, so $x\in V(T)\cap V(Q_0)$ and $C_x\subseteq B_x\cup B_\beta$.  We claim that $\torso{G}{\calQ,C_x}-B_\beta$ is a subgraph of $\torso{G}{\calT,B_x}-B_\beta$.  Let $xy$ be an edge of $Q$ incident to $x$. There are three possibilities:
  \begin{itemize}[nosep,nolistsep]
    \item If $y\in S$ then $C_x\subseteq B_x\cup B_\beta$ and $C_y\subseteq B_y\cup B_\beta$.  Then $C_x\cap C_y\setminus B_\beta=B_x\cap B_y\setminus B_\beta$ is a clique in both $\torso{G}{\calQ,C_x}-B_\beta$ and $\torso{G}{\calT,B_x}-B_\beta$.
    \item If $y=v_i$ for some $i\in[s]$ with $|N_i|=2$, then $x\in\{a_i,b_i\}$. Without loss of generality, suppose $x=a_i$.  Then $C_x\cap C_y\setminus B_\beta=B_{a_i}\cap V(G_i)\setminus B_\beta$ is contained in the original adhesion $B_{a_i}\cap B_{x_i}$ incident to $x$.  Thus $C_x\cap C_y\setminus B_\beta$ is a clique in both $\torso{G}{\calQ,C_x}-B_\beta$ and $\torso{G}{\calT,B_x}-B_\beta$.
    \item If $xy=e_i=a_i\alpha_i$ for some $i\in[s]$ with $|N_i|=1$ and $m_i>0$, then $C_x\cap C_y\setminus B_\beta = B_{a_i}\cap V(G_i)\setminus B_\beta$ is contained in the original adhesion $B_{a_i}\cap B_{x_i}$ incident to $x$.  Thus $C_x\cap C_y\setminus B_\beta$ is a clique in both $\torso{G}{\calQ,C_x}-B_\beta$ and $\torso{G}{\calT,B_x}-B_\beta$.
  \end{itemize}
  Thus, for each $x\in S$, $\torso{G}{\calQ,C_x}-B_\beta$ is a subgraph of $\torso{G}{\calT,B_x}-B_\beta$, so $x$ satisfies the second alternative of property~\ref{req:adhesion_torso_0}, since $B_\beta$ is an adhesion of $\calT$.

  For each $x\in V(Q_0)\setminus S$, $x=v_i$ for some $i\in[s]$ with $|N_i|=2$.  Then $C_x=C_{v_i}\subseteq (B_{a_i}\cap B_{x_i})\cup (B_{b_i}\cap B_{y_i})\cup B_\beta$ is contained in the union of at most three adhesions of $\calT$ so $x$ satisfies the first alternative of property~\ref{req:three_adhesions_0}.

  Lastly, for each $x\in V(Q)\setminus V(Q_0)$, there are two possibilities:
  \begin{itemize}[nosep,nolistsep]
    \item $x=\alpha_i$ for some $i\in[s]$ with $m_i>0$. By \ref{req:adhesion_pairs_0}, $C_{\alpha_i}=B_{\alpha_i}\cup B_{\beta_i}$. Each of $B_{\alpha_i}$ and $B_{\beta_i}$ is an adhesion of $\calT$, so $x$ satisfies the first alternative of property~\ref{req:three_adhesions_0}.
    \item $x$ is a non-root node of $Q_i$ for some $i\in[s]$ with $m_i>0$.  By induction, $x$ satisfies \ref{req:subtree_torso_apex_0} with respect to $G_i$, $\calT_i$, and $\calQ_i$. 
    By definition, $\torso{G_i}{\calT_i,B_y}=\torso{G}{\calT,B_y}$ for each $y\in V(C_i)$, and every adhesion of $\calT_i$ is an adhesion of $\calT$, and $\torso{G_i}{\calQ_i,C_x}=\torso{G}{\calQ,C_x}$. Therefore, $x$ also satisfies \ref{req:subtree_torso_apex_0} with respect to $G$, $\calT$, and $\calQ$. \qedhere
  \end{itemize}
\end{proof}

Let $\calF:=(F,(B_x\mid x\in V(F)))$ be a rooted forest-decomposition of a graph $G$ and, for each $x\in V(F)$, let $A_x$ denote the parent adhesion of $x$. We say that $\calF$ is \defin{tidy} if
\begin{enumerate*}[label=(\arabic*)]
  \item $B_x\setminus A_x\neq\emptyset$ for each $x\in V(F)$; and
  \item $A_y\setminus A_x\neq\emptyset$ for each non-root $y\in V(F)$ with $x=\parent_F(y)$.
\end{enumerate*}

\begin{lem}\label{all_singing_all_dancing}
  Let $G_0$ be an $n$-vertex graph, with $n\in\N^+$, and let $\calT$ be a tree-decomposition of $G_0$.  Then, for any integer $b\ge 2$, there exists a graph $G$ with an isolated vertex $v_0$ such that $G-v_0=G_0$, a tidy rooted tree-decomposition $\calQ:=(Q,(C_x\mid x\in V(Q)))$ of $G$ rooted at $x_0\in V(Q)$, and a subset of edges $Y\subseteq E(Q)$ with the following properties:
  \begin{enumerate}[nosep,nolistsep,label=\rm(\alph*)]
    \item\label{req:root_bag} $C_{x_0}=\{v_0\}$;
    \item\label{req:height} each root-to-leaf path in $Q$ contains at most $2\log_b(n+1) + 1$ edges of $Y$;
    \item\label{req:subtree_size} each connected component of $Q-Y$ has at most $4b-1$ vertices;
    \item\label{req:same_size} letting $A_x$ denote the parent adhesion of $x$ in $\calQ$ (with $A_{x_0}:=\emptyset$), for each connected component $C$ of $Q-Y$ and each $x,y\in V(C)$,
    \[
      \textstyle \left|\bigcup_{z\in V(Q_x)} (C_z\setminus A_z)\right| \le b \left|\bigcup_{z\in V(Q_y)} (C_z\setminus A_z)\right| \enspace ;
    \]
    \item\label{req:adhesion_width} each adhesion of $\calQ$ is either $\{v_0\}$ or is contained in the union of at most two adhesions of $\calT$; in particular, if $\calT$ has adhesion-width at most $k$, then $\calQ$ has adhesion-width at most $\max\{2k,1\}$; and
    \item\label{req:subtree_torso_apex} for each $x\in V(Q)$, at least one of the following holds: 
    \begin{itemize}[nosep,nolistsep]
      \item\label{req:three_adhesions} $C_x\setminus\{v_0\}$ is contained in the union of at most three adhesions of $\calT$; or
      \item\label{req:adhesion_torso} $\torso{G}{\calQ,C_x}-(A\cup\{v_0\})$ is a subgraph of some torso of $\calT$, where $A=\emptyset$ or $A$ is an adhesion of $\calT$.
    \end{itemize}
  \end{enumerate}
\end{lem}

\begin{proof}
  Let $\calT_0:=(T_0,(B_x\mid x\in V(T_0)))$ and $Y_0\subseteq E(T_0)$ be the tree-decomposition of $G_0$ obtained by applying \cref{log_b_td} to the tree-decomposition $\calT$.  Thus, each torso of $\calT_0$ satisfies \cref{log_b_td}\ref{req:subtree_torso_apex_0}.    
  First, we transform $\calT_0$ and $Y_0$ into a tidy rooted forest-decomposition $\calF=(F,(B_x\mid x\in V(F)))$ and $Y_1\subseteq E(F)$ in five steps:
  \begin{enumerate}[nosep,nolistsep,label=(\arabic*)]
    \item Initially, set $F:=T_0$ and $Y_1:=Y_0$.
    \item Remove every node $x$ of $F$ with $B_x=\emptyset$ (along with any incident edges in $F$ and $Y_1$).
    \item Remove every edge $xy$ of $F$ (and $Y_1$) with $B_x\cap B_y=\emptyset$.
    \item Repeatedly replace an edge $yz\in E(F)$ with the edge $xz$ if $x$ is the $F$-parent of $y$, $y$ is the $F$-parent of $z$, and $B_z\cap B_y\subseteq B_x$. 
    
    If the edge $yz$ is in $Y_1$, replace $yz$ with $xz$ in $Y_1$.  
    
    If the edge $xy$ is in $Y_1$, add the edge $xz$ to $Y_1$ (and leave $xy$ in $Y_1$).
    \item Repeatedly remove every leaf $y$ of $F$ with $B_y\subseteq B_{\parent(y)}$ (along with the edge $\parent(y)y$ from $F$ and $Y_1$). (Note that any non-root node $y$ with $B_y\subseteq B_{\parent(y)}$ has $B_z\cap B_y\subseteq B_y\subseteq B_{\parent(y)}$ for each child $z$ of $y$, so $y$ loses all its children in step~(4) and is removed in step~(5).)
  \end{enumerate} 
  Let $\calF=(F,(B_x\mid x\in V(F)))$ be the resulting tidy forest-decomposition.  (For each $x\in V(F)$, the bag $B_x$ is well-defined since $V(F)\subseteq V(T_0)$.) 
  Note that each torso of $\calF$ is a subgraph of a torso of $\calT_0$.
  In step (4), the new edge $xz$ satisfies $B_x\cap B_z = B_x\cap (B_y\cap B_z) = B_y\cap B_z$, so the adhesion on $xz$ is identical to the adhesion on $yz$. Thus every adhesion of $\calF$ is an adhesion of $\calT_0$.
  Furthermore, for any node $z\in V(F)$, the $F$-ancestors of $z$ form a subset of the $T_0$-ancestors of $z$.  Furthermore, if $x$ is both an $F$-ancestor and $T_0$-ancestor of $z$ then the path from $x$ to $z$ in $F$ contains an edge in $Y_1$ if and only if the path from $x$ to $z$ in $T_0$ contains an edge in $Y_0$. Therefore, for any connected component $C$ of $F-Y_1$, $V(C)\subseteq V(C_0)$ for some component $C_0$ of $T_0-Y_0$.  Thus $|V(C)|\le |V(C_0)|\le 4b-1$ for each component $C$ of $F-Y_1$.  Let $Q$ be the tree obtained from $F$ by introducing a vertex $x_0$ and adding an edge $x_0x'$ for each root $x'$ of $F$. Let $Y_2:=Y_1\cup \{x_0x'\mid x'\text{ is a root of }F\}$. Then each component of $Q-Y_2$ has at most $4b-1$ nodes.
  
  Let $G$ be the graph obtained by adding an isolated vertex $v_0$ to $G_0$.  Define $C_{x_0}:=\{v_0\}$. For each child $x'$ of $x_0$, define $C_{x'}:=B_{x'}\cup\{v_0\}$.  For each other vertex $x$ in $Q$, define $C_x:=B_x$.  Then $\calQ:=(Q,(C_x\mid x\in V(Q)))$ is a tidy rooted tree-decomposition of $G$ that satisfies \ref{req:root_bag} and \ref{req:subtree_torso_apex}. Moreover, for each root $x'$ of $F$, the adhesion on $x_0x'$ is $C_{x_0}\cap C_{x'}=\{v_0\}$, while for every other edge $xy\in E(Q)$, $C_x\cap C_y = B_x\cap B_y$ is an adhesion of $\calF$ (and hence of $\calT_0$). By \cref{log_b_td}\ref{req:adhesion_width_0}, each adhesion of $\calT_0$ is contained in the union of at most two adhesions of $\calT$. Therefore, $\calQ$ satisfies \ref{req:adhesion_width}.

  Let $A_{x_0}:=\emptyset$ and, for each vertex $x\in V(Q)\setminus\{x_0\}$, let $A_x$ be the parent adhesion of $x$ in $\calQ$.  For each $x\in V(Q)$, define $\omega(x):=|C_x\setminus A_x|$.  Since $\calQ$ is tidy, $\omega(x)\ge 1$ for each $x\in V(Q)$, so $\omega:V(Q)\to\N^+$ is a weight function with total weight $\omega(Q)=\sum_{x\in V(Q)}|C_x\setminus A_x|=|V(G)|=n+1$.  Let $Y_3\subseteq E(Q)$ be the set of edges obtained by applying \cref{skinny_edges} to $Q$ with the weight function $\omega$.  Then, for each component $C$ of $Q-Y_3$ and each $x,y\in V(C)$, we have   
  \[
      \textstyle \left|\bigcup_{z\in V(Q_x)}(C_z\setminus A_z)\right|=\omega(Q_x)\le b\cdot\omega(Q_y)=b\left|\bigcup_{z\in V(Q_y)}(C_z\setminus A_z)\right| \enspace .
  \]
  Let $Y:=Y_2\cup Y_3$.  Since each component of $Q-Y$ is contained in a component of $Q-Y_2$ and in a component of $Q-Y_3$, the pair $(\calQ,Y)$ satisfies \ref{req:subtree_size} and \ref{req:same_size}.  Each root-to-leaf path in $Q$ contains at most $1+\log_b n$ edges in $Y_2$ and at most $\log_b(n+1)$ edges in $Y_3$.  Thus, each root-to-leaf path in $Q$ contains at most $1 + \log_b n + \log_b(n+1) \le 2\log_b(n+1)+1$ edges in $Y$, so $(\calQ,Y)$ satisfies \ref{req:height}.
\end{proof}

\section{Weighted Mixed Labelling Schemes}
\label{weighted_mixed_section}
 
Let $\calG$ be a class of graphs. We say that $\tau:=(G^+,G,\omega)$ is a \defin{$\calG$-admissible triple} if $G^+$ is a graph in $\calG$ with at least one vertex, $G$ is a spanning subgraph of $G^+$, and $\omega:V(G)\to\R^+$ is a weight function.  Let \defin{$\calT_\calG$} be the set of all $\calG$-admissible triples. For each $n\in\N$, let $\mathdefin{\calT_{\calG,\le n}}:=\{(G^+,G,\omega)\in\calT_\calG:|V(G)|\le n\}$.

Let $g_i:\N\to\R$ for each $i\in[3]$.
Let $A:(\{0,1\}^*)^2\to\{0,1\}$ (the \defin{adjacency tester}) and $I:(\{0,1\}^*)^3\to\{0,1\}$ (the \defin{identity tester}).
For a $\calG$-admissible triple $\tau:=(G^+,G,\omega)$, let $n:=|V(G)|$.
A pair of functions
\begin{itemize}[nosep,nolistsep]
  \item $\mu: V(G)\cup \{K\mid \text{$K$ is a clique in $G^+$}\} \to \set{0,1}^*$ (the \defin{mixed labelling function}); and 
  \item $\kappa:\{(K, u)\mid \text{$K$ is a clique in $G^+$}, u\in K\}\to \set{0,1}^*$ (the \defin{local identifier function})
\end{itemize} 
is a \defin{mixed labelling} of $(G^+,G)$ for $(A,I)$ if it satisfies:
\begin{enumerate}[nosep,nolistsep,label=(M\arabic*)]
  \item\label{mls:injective_on_cliques} for all pairs of distinct cliques $K$ and $L$ in $G^+$, $\mu(K)\neq\mu(L)$;
  \item\label{mls:adjacency_tester} for all $u,v\in V(G)$, $A(\mu(u),\mu(v))=1$ if and only if $uv\in E(G)$;
  \item\label{mls:identity_tester} for every clique $K$ in $G^+$, every $u\in K$, and every $v\in V(G)$, $I(\mu(K),\kappa(K, u),\mu(v))=1$ if and only if $u=v$.
\end{enumerate}
If $(\mu,\kappa)$ additionally satisfies the length bounds
\begin{enumerate}[resume,nosep,nolistsep,label=(M\arabic*)]
  \item\label{mls:vertex_label_length} for every $v\in V(G)$, $|\mu(v)| \le \log \omega(G) - \log\omega(v) + g_1(n)$;
  \item\label{mls:clique_label_length} for every clique $K$ in $G^+$, $|\mu(K)|\le \log \omega(G) - \log\min_{v\in K}\omega(v) + g_3(n)$;
  \item\label{mls:local_identifier_length}  for every clique $K$ in $G^+$ and every $u\in K$, $|\kappa(K, u)| \le g_2(n)$,
\end{enumerate}
then $(\mu,\kappa)$ is called a \defin{weighted $(g_1,g_2,g_3)$ mixed labelling} of $(G^+,G,\omega)$ for $(A,I)$.  

Here is a little intuition about the last three conditions. The function $g_1(n)$ bounds how much the length of a vertex label exceeds the ideal value $\log \omega(G) - \log\omega(v)$ that appears in Kraft's inequality \citep{Kraft1949}. The function $g_2(n)$ is the length of the longest local identifier. The function $g_3(n)$ bounds how much the length of a clique label exceeds the Kraft-like quantity $\log\omega(G)-\log\min_{v\in K}\omega(v)$.

To prove \cref{main_result_induced_universal_sparse}, it will be helpful to fix a particular mixed labelling of each $\calG$-admissible triple.  
A \defin{weighted $(g_1,g_2,g_3)$ mixed labelling} of $\calG$ for $(A,I)$ is an indexed collection $\mathdefin{\calS}:=((\mu_\tau,\kappa_\tau):\tau\in\calT_\calG)$ where $(\mu_\tau,\kappa_\tau)$ is a weighted $(g_1,g_2,g_3)$ mixed labelling of $\tau$ for $(A,I)$, for each $\tau\in\calT_\calG$.  The triple $(\calS,A,I)$ is called a \defin{weighted $(g_1,g_2,g_3)$ mixed labelling scheme} for $\calG$. We say that $\calG$ \defin{admits} a weighted $(g_1,g_2,g_3)$ mixed labelling scheme if such a triple $(\calS,A,I)$ exists.  Similar definitions apply to adjacency labelling schemes. An \defin{$f(n)$-bit adjacency labelling} of $\calG$ for $A$ is an indexed collection $\mathdefin{\calA}:=(\ell_G:G\in\calG)$ where, for each $G\in\calG$,  $\ell_G$ is an adjacency labelling of $G$ for $A$ and $|\ell_G(v)|\le f(|V(G)|)$ for each $v\in V(G)$.  The pair $(\calA,A)$ is an \defin{$f(n)$-bit adjacency labelling scheme} for $\calG$.

\subsection{Basic Properties}
\label{basic_properties}

Let $(\calS,A,I)$ be a weighted $(g_1,g_2,g_3)$ mixed labelling scheme for a graph class $\calG$, where $\calS = ((\mu_\tau,\kappa_\tau):\tau\in\calT_\calG)$. For any graph $G$, let $\omega_1:V(G)\to\{1\}$ be the \defin{uniform weight function} defined by $\omega_1(v):=1$ for all $v\in V(G)$.  Let $\calG^-$ be the class of all spanning subgraphs of graphs in $\calG$. For each non-empty $G\in\calG^-$, choose a graph $G^+\in\calG$ such that $G$ is a spanning subgraph of $G^+$, and let $\mathdefin{\calA}:=(\mu_{(G^+,G,\omega_1)}|_{V(G)}:G\in\calG^-)$. Then $(\calA,A)$ is a $(\log n+g_1(n))$-bit adjacency labelling scheme for $\calG^-$.  Thus, we are interested in the case where $g_1(n)\in o(\log n)$. We say that a graph class $\calG$ admits an \defin{efficient} weighted mixed labelling scheme if $\calG$ admits a weighted $(g_1,g_2,g_3)$ mixed labelling scheme with $g_i(n)\in o(\log n)$ for each $i\in [3]$. So any graph class that admits an efficient weighted mixed labelling scheme admits a $(1+o(1))\log n$-bit adjacency labelling scheme. 

An adjacency labelling $\ell$ of a graph $G$ is \defin{injective} if $\ell(u)\neq\ell(v)$ for all distinct $u,v\in V(G)$. An adjacency labelling scheme $(\calA,A)$ for $\calG$ is \defin{injective} if $\ell_G$ is injective for each $G\in\calG$. The following lemma shows that any adjacency labelling obtained from a mixed labelling scheme is injective, which justifies \cref{main_result_induced_universal}.

\begin{lem}\label{injective_vertex_labels}
  Let $G^+$ be a graph with at least one vertex, let $G$ be a spanning subgraph of $G^+$, and let $(\mu,\kappa)$ be a mixed labelling of $(G^+,G)$ for $(A,I)$ where $A$ is an adjacency tester and $I$ is an identity tester. Then the restriction $\mu|_{V(G)}$ of $\mu$ to $V(G)$ is injective.
\end{lem}

\begin{proof}
  Suppose, for the sake of contradiction, that there exist distinct $v,w\in V(G)$ such that $\mu(v)=\mu(w)$. Since $K:=\{v\}$ is a clique in $G^+$, \ref{mls:identity_tester} implies that $I(\mu(K),\kappa(K,v),\mu(v))=1$ and $I(\mu(K),\kappa(K,v),\mu(w))=0$. Since $\mu(v)=\mu(w)$, this immediately leads to the contradiction
  \[
    1 = I(\mu(K),\kappa(K,v),\mu(v)) = I(\mu(K),\kappa(K,v),\mu(w)) = 0 \enspace. \qedhere
  \]
\end{proof}

It will actually be convenient to consider a slightly weaker notion of weighted mixed labelling schemes.  A $\calG$-admissible triple $(G^+,G,\omega)$ is \defin{nice} if $\omega:V(G)\to\N^+$ is a nice weight function.  A \defin{weak weighted $(g_1,g_2,g_3)$ mixed labelling scheme} $(\calS,A,I)$ for a graph class $\calG$ satisfies all the conditions of a weighted $(g_1,g_2,g_3)$ mixed labelling scheme for $\calG$, with the exception that $\calS$ only contains mixed labellings for nice triples $(G^+,G,\omega)\in\calT_\calG$. 
We say that $\calG$ \defin{admits} a weak weighted $(g_1,g_2,g_3)$ mixed labelling scheme if such a scheme exists. Weak mixed labelling schemes and weighted mixed labelling schemes are essentially equivalent, as shown by the following lemma. This will allow us to only consider nice weight functions in the proofs.

\begin{lem}
  \label{lem:nice_weights}
  Let $\calG$ be a class of graphs.
  Let $g_i:\N\to\R$ for each $i\in[3]$.
  If $\calG$ admits a weak weighted $(g_1,g_2,g_3)$ mixed labelling scheme, then $\calG$ admits a weighted $(g_1+1,g_2,g_3+1)$ mixed labelling scheme.
\end{lem}
\begin{proof}
  Let $(\calS,A,I)$ be a weak weighted $(g_1,g_2,g_3)$ mixed labelling scheme for $\calG$. We will use the same pair $(A, I)$ for our weighted $(g_1+1,g_2,g_3+1)$ mixed labelling scheme. Let $\tau=(G^+,G,\omega)\in\calT_\calG$ and let $n:=|V(G)|$. Let $\omega'$ be the nice weight function defined by $\omega'(v):=\ceil{n\omega(v)/\omega(G)}$ for each $v\in V(G)$, as in the proof of \cref{obs:nice}. Thus, for every $v\in V(G)$,
  \[
    \log\omega'(G)-\log\omega'(v)\le\min\{\log n,\;\log\omega(G)-\log\omega(v)\}+1 \enspace.
  \]
  Therefore, $\calS$ contains a weighted $(g_1,g_2,g_3)$ mixed labelling $(\mu,\kappa)$ of $(G^+,G,\omega')$ for $(A,I)$.
  Conditions \ref{mls:injective_on_cliques}--\ref{mls:identity_tester} depend only on $G^+$ and $G$, so they hold for $(\mu,\kappa)$. For every $v\in V(G)$, \ref{mls:vertex_label_length} gives
  \[
    |\mu(v)| \leq \log \omega'(G) - \log\omega'(v) + g_1(n) \leq \log \omega(G) - \log\omega(v) + g_1(n) + 1\enspace.
  \]
  For every clique $K$ in $G^+$, if $x\in K$ minimizes $\omega'(x)$, then $\log\omega'(G)-\log\omega'(x) \le \log\omega(G)-\log\omega(x)+1 \le \log\omega(G)-\log\min_{v\in K}\omega(v)+1$, and thus by \ref{mls:clique_label_length},
  \[
    |\mu(K)|\le \log \omega'(G) - \log\min_{v\in K}\omega'(v) + g_3(n)
    \leq \log \omega(G) - \log\min_{v\in K}\omega(v) + g_3(n) + 1\enspace.
  \]
  Finally, \ref{mls:local_identifier_length} gives $|\kappa(K,u)|\le g_2(n)$ for every clique $K$ in $G^+$ and every $u\in K$. It follows that $(\mu,\kappa)$ is a weighted $(g_1+1,g_2,g_3+1)$ mixed labelling of $(G^+,G,\omega)$ for $(A,I)$. Setting $(\mu_\tau,\kappa_\tau) := (\mu,\kappa)$ for each $\tau\in\calT_\calG$, the indexed collection $\mathdefin{\calS'}:=((\mu_\tau,\kappa_\tau):\tau\in\calT_\calG)$ yields the desired weighted $(g_1+1,g_2,g_3+1)$ mixed labelling scheme $(\calS',A,I)$ for $\calG$.
\end{proof}

If the ultimate goal is to obtain an adjacency labelling scheme, then an even weaker notion suffices. Recall that $\omega_1:V(G)\to\{1\}$ denotes the uniform weight function. We say that a class of graphs $\calG$ \defin{admits an unweighted $(g_1,g_2,g_3)$ mixed labelling scheme} $(\calS,A,I)$ if $\calS$ satisfies all the conditions of a weighted $(g_1,g_2,g_3)$ mixed labelling scheme for $\calG$, with the exception that $\calS$ only contains mixed labellings for triples $(G^+,G,\omega_1)\in\calT_\calG$ in which $\omega_1$ is the uniform weight function.

\begin{obs}\label{adj_labelling_from_mixed_labelling}
  If a graph class $\calG$ admits an unweighted $(g_1,g_2,g_3)$ mixed labelling scheme, then it admits an injective $f(n)$-bit adjacency labelling scheme where $f(n)=\log n+g_1(n)$. Furthermore, if the mixed labelling scheme is $g_0$-orientable (defined in \cref{sec:orientable}, below), then the resulting adjacency labelling scheme is also $g_0$-orientable.
\end{obs}
\begin{proof}
  Let $(\calS,A,I)$ be an unweighted $(g_1,g_2,g_3)$ mixed labelling scheme for $\calG$. For each non-empty graph $G\in\calG$, let $\tau := (G, G, \omega_1)$, where $\omega_1:V(G)\to\{1\}$ is the uniform weight function, and let $(\mu_\tau,\kappa_\tau)$ be the mixed labelling of $\tau$ in $\calS$. By \ref{mls:adjacency_tester}, $\ell_G := \mu_{\tau}|_{V(G)}$ is an adjacency labelling of $G$ for $A$. For $n:=|V(G)|$, \ref{mls:vertex_label_length} implies that for each $v\in V(G)$,
  \[
    |\ell_G(v)| \le \log\omega_1(G) - \log\omega_1(v) + g_1(n) = \log n + g_1(n) \enspace.
  \]
  By \cref{injective_vertex_labels}, $\ell_G$ is injective. Thus, $(\calA, A)$ is an injective $(\log n + g_1(n))$-bit adjacency labelling scheme for $\calG$, where $\mathdefin{\calA} := (\ell_G : G\in\calG)$. The orientability claim is verified in \cref{sec:orientable}.
\end{proof}

\subsection{Orientable Labelling Schemes}\label{sec:orientable}

In order to establish \cref{main_result_induced_universal_sparse}, which promises sparse universal graphs, we need two additional definitions.
Let $(\calS,A,I)$ be a weighted $(g_1,g_2,g_3)$ mixed labelling scheme for a graph class $\calG$, where $\calS = ((\mu_\tau,\kappa_\tau):\tau\in\calT_\calG)$.  For each $n\in\N^+$, let
\[
  \mathdefin{\calT'_{\calG,\le n}}:=\left\{(G^+,G,\omega)\in\calT_{\calG,\le n}: \omega:V(G)\to\N^+,\; \omega(G)\le 2n\right\}.
\]
The triples in $\calT'_{\calG,\le n}$ need not be nice: their total weight is bounded by $2n$, which may be much larger than their number of vertices. The \defin{$n$-universal graph} $U_n$ defined by $(\calS,A)$ is the graph with vertex set
\[
V(U_n):=\left\{ \mu_\tau(v): \tau = (G^+,G,\omega)\in\calT'_{\calG,\le n},\; v\in V(G) \right\}
\]
and edge set
\[
  E(U_n):=\left\{ \{\mu_\tau(v),\mu_\tau(v')\}: \tau = (G^+,G,\omega)\in\calT'_{\calG,\le n},\; vv'\in E(G) \right\}
  \enspace.
\]
These definitions imply that, for each $\tau = (G^+,G,\omega)\in\calT'_{\calG,\le n}$, the graph $G$ is isomorphic to the induced subgraph $U_n[\{\mu_\tau(v):v\in V(G)\}]$.  Taking uniform weights, we see that $U_n$ contains every graph $G\in\calG$ with at most $n$ vertices as an induced subgraph. (In fact, $U_n$ contains every subgraph of every at-most-$n$-vertex graph in $\calG$, even if $\calG$ itself is not monotone.)  Note that $U_n$ is finite. Indeed, for each $\tau=(G^+,G,\omega)\in\calT'_{\calG,\le n}$ and each $v\in V(G)$, we have $\omega(v)\ge 1$ and $\omega(G)\le 2n$, so
\[
  |\mu_\tau(v)|\le \log\omega(G)-\log\omega(v)+g_1(|V(G)|)
  \le \log n+1+\max_{1\le m\le n}g_1(m)\enspace.
\]
If $g_1$ is non-decreasing, this bound is at most $\log n+1+g_1(n)$.

Since the graph $U_n$ defined above depends only on the vertex labellings in $\calS$ (and not on clique labels or local identifiers), the same definition applies directly to any injective adjacency labelling $\calA:=(\ell_G:G\in\calG)$ of $\calG$ for $A$. In this case, for each $n\in\N^+$ and each $G\in\calG$ with at most $n$ vertices, the $n$-universal graph $U_n$ defined by $(\calA,A)$ has an induced subgraph $U_n[\{\ell_G(v):v\in V(G)\}]$ that is isomorphic to $G$.

Let $g_0:\N\to\R^{\ge 0}$.
We say that $(\calS,A)$ (or $(\calA,A)$) is \defin{$g_0$-orientable} if, for each $n\in\N^+$, the $n$-universal graph $U_n$ defined by $(\calS,A)$ (or $(\calA,A)$) has an orientation in which each vertex has outdegree at most $g_0(n)$. A weighted mixed labelling scheme $(\calS,A,I)$ is \defin{$g_0$-orientable} if $(\calS,A)$ is $g_0$-orientable. This definition immediately implies that the $n$-universal graph $U_n$ has at most $|V(U_n)|\cdot g_0(n)$ edges. We say that a graph class $\calG$ \defin{admits a $g_0$-orientable weighted mixed labelling scheme} if such a scheme exists (and similarly for adjacency labelling schemes).
Observe that if $(\calS,A,I)$ is $g_0$-orientable, then the injective adjacency labelling scheme $(\calA,A)$ obtained from $\calS$ in \cref{adj_labelling_from_mixed_labelling} is also $g_0$-orientable. Indeed, for each $G\in\calG$ with $|V(G)|\le n$, the triple $\tau=(G,G,\omega_1)$ belongs to $\calT'_{\calG,\le n}$, so the $n$-universal graph of $(\calA,A)$ is an induced subgraph of the $n$-universal graph of $(\calS,A)$, and therefore inherits an orientation of maximum outdegree at most $g_0(n)$.

To control the number of edges when combining mixed labellings, we also need to bound the number of vertex labels associated with a fixed clique label and local identifier. This motivates the following strengthening of orientability.
We say that $(\calS,A,I)$ is \defin{strongly $g_0$-orientable} if $(\calS,A)$ is $g_0$-orientable and, for each $n\in\N^+$, each triple $\tau_0:=(G_0^+,G_0,\omega_0)\in\calT'_{\calG,\le n}$, each clique $K_0$ in $G_0^+$, and each $u_0\in K_0$, there are at most $g_0(n)$ binary strings $\beta$ for which there exists a triple $\tau:=(G^+,G,\omega)\in\calT'_{\calG,\le n}$, a clique $K$ in $G^+$, and a vertex $u\in K$ such that
\begin{enumerate}[nosep,nolistsep,label=(\alph*)]
  \item $\mu_\tau(K)=\mu_{\tau_0}(K_0)$;
  \item $\kappa_\tau(K,u)=\kappa_{\tau_0}(K_0,u_0)$; and
  \item $\mu_\tau(u)=\beta$.
\end{enumerate}
We say that $\calG$ \defin{admits a strongly $g_0$-orientable weighted $(g_1,g_2,g_3)$ mixed labelling scheme} if such a scheme exists.

Less formally, fix a clique label $\zeta$ and a local identifier $\gamma$ that occur together in a triple in $\calT'_{\calG,\le n}$. Strong $g_0$-orientability guarantees that at most $g_0(n)$ distinct vertex labels $\beta$ occur with this same pair: there must be a triple $\tau=(G^+,G,\omega)\in\calT'_{\calG,\le n}$, a clique $K$ in $G^+$, and a vertex $u\in K$ such that $\mu_\tau(K)=\zeta$, $\kappa_\tau(K,u)=\gamma$, and $\mu_\tau(u)=\beta$. Any such occurrence implies that $I(\zeta,\gamma,\beta)=1$.

Our construction will construct an adjacency labelling scheme $(\calA',A')$ for a graph class $\calG'$. The graphs in $\calG'$ will contain subgraphs that are in some class $\calG$ that admits a weighted mixed labelling scheme $(\calS,A,I)$.  In our constructions, a vertex label $\ell_{G'}(v)$ may contain $\zeta$ and $\gamma$, while the label of a neighbour $w$ contains a corresponding vertex label $\beta$. If every edge handled in this way has a witness as above, then fixing $\ell_{G'}(v)$ leaves at most $g_0(n)$ choices for $\beta$. If, for each such $\beta$, there are at most $M(n)$ possible complete labels $\ell_{G'}(w)$, then orienting these edges from $\ell_{G'}(v)$ towards $\ell_{G'}(w)$ contributes at most $g_0(n)\cdot M(n)$ to the outdegree of $\ell_{G'}(v)$. Here $M(n)$ accounts for the fact that $\beta$ is only \emph{part} of the vertex label $\ell_{G'}(w)$. There may be up to $M(n)$ vertex labels $\ell_{G'}(w)$ that include $\beta$.

The definitions of $U_n$, $g_0$-orientability, and strong $g_0$-orientability also apply to weak weighted mixed labelling schemes, with all triples restricted to nice triples in $\calT'_{\calG,\le n}$.

\begin{obs}\label{obs:normalization_orientability}
  The construction in \cref{lem:nice_weights} preserves $g_0$-orientability and strong $g_0$-orientability.
\end{obs}
\begin{proof}
  Let $(\calS,A,I)$ be a weak weighted mixed labelling scheme and let $(\calS',A,I)$ be the weighted mixed labelling scheme constructed in \cref{lem:nice_weights}. Fix $n\in\N^+$. For every $\tau=(G^+,G,\omega)\in\calT'_{\calG,\le n}$, the construction assigns to $\tau$ the labelling in $\calS$ of $\tau'=(G^+,G,\omega')$, where $\omega'(v)=\ceil{|V(G)|\omega(v)/\omega(G)}$. By \cref{obs:nice}, $\tau'$ is nice and $\omega'(G)\le\omega(G)\le 2n$, so $\tau'\in\calT'_{\calG,\le n}$. Consequently, the $n$-universal graph defined by $(\calS',A)$ is a subgraph of the one defined by $(\calS,A)$, and any orientation of the latter restricts to an orientation of the former with no larger outdegrees. Similarly, every occurrence of a clique label, local identifier, and vertex label in $\calS'$ on a triple in $\calT'_{\calG,\le n}$ also occurs in $\calS$ on a nice triple in the same set. Thus, the bound required for strong $g_0$-orientability is preserved as well.
\end{proof}

\subsection{Adding Apex Vertices}
\label{sec:adding_apexes}

The next lemma shows how a mixed labelling scheme can accommodate the addition of apex vertices. Recall that $\calG^{+a}$ is the class of graphs $G$ such that $G-X\in\calG$ for some $X\subseteq V(G)$ with $|X|\leq a$.

\begin{lem}\label{add_apexes}
  Let $\calG$ be a class of graphs that admits a strongly $g_0$-orientable weighted $(g_1,g_2,g_3)$ mixed labelling scheme, where $g_0:\N\to\R^{\ge 0}$ is a function, and $g_1,g_2,g_3:\N\to\R^{\ge 0}$ are non-decreasing functions with $g_1(n),g_3(n)\in \Oh(\log n)$.
  Let $a\in\N^+$. Then $\calG^{+a}$ admits a strongly $g_0'$-orientable weighted $(g'_1,g'_2,g'_3)$ mixed labelling scheme, where
  \begin{align*}
    g_0'(n)
     & = (2g_0(n) + a)2^{a} \enspace,  \\
    g'_1(n)
     & = g_1(n)+\Oh(a)+\Oh(\log\log n) \enspace,  \\
    g'_2(n)
     & = g_2(n)+\Oh(\log a) \enspace,             \\
    g'_3(n)
     & = g_3(n)+\Oh(a)+\Oh(\log\log n) \enspace .
  \end{align*}
\end{lem}

\begin{proof} 
  Let $(\calS,A,I)$ be a strongly $g_0$-orientable weighted $(g_1,g_2,g_3)$ mixed labelling scheme for $\calG$, where $\calS:=((\mu^-_\tau,\kappa^-_\tau):\tau\in\calT_\calG)$.  For each $n\in\N^+$, let $U_n$ be the $n$-universal graph defined by $(\calS,A)$.

  By \cref{lem:nice_weights,obs:normalization_orientability}, it suffices to construct a strongly $g'_0$-orientable weak weighted mixed labelling scheme with the claimed asymptotic label-length bounds. The additive constants needed to extend it to arbitrary weights are absorbed into the $\Oh(a)$ terms in $g'_1$ and $g'_3$. We will describe testers $A'$ and $I'$ and construct a mixed labelling $(\mu_\tau,\kappa_\tau)$ for each nice triple $\tau\in\calT_{\calG^{+a}}$. Let $\calS'$ be the indexed collection of these labellings over all nice triples.   To fix a specific labelling $(\mu_\tau,\kappa_\tau)$ for each nice $\calG^{+a}$-admissible triple $\tau=(G^+,G,\omega)$, we first choose, for each $G^+\in\calG^{+a}$, a specific set $X(G^+)\subseteq V(G^+)$ with $|X(G^+)|\le a$ such that $G^+-X(G^+)\in\calG$.

  Let $\tau:=(G^+,G,\omega)\in\calT_{\calG^{+a}}$ be a triple with a nice weight function $\omega$, let $X:=X(G^+)$, and let $n:=|V(G)|$. Since $\omega$ is nice, $\omega(G)\le 2n$ and $\log\omega(G)-\log\omega(v)\le\log n + 1$ for each $v\in V(G)$. We will describe the weighted mixed labelling $(\mu,\kappa):=(\mu_\tau,\kappa_\tau)$ of $\tau$ for $(A',I')$.

  \paragraph{\boldmath The subgraph label $\mu^-$:}
  If $G-X$ is non-empty, then by the definition of $X$, $\tau^-:=(G^+-X,G-X,\omega|_{V(G-X)})$ is a $\calG$-admissible triple with $|V(G-X)|\le n$.   Let $(\mu^-,\kappa^-):=(\mu^-_{\tau^-},\kappa^-_{\tau^-})$.  Thus, $(\mu^-,\kappa^-)$ is a weighted $(g_1,g_2,g_3)$ mixed labelling of $(G^+-X, G-X,\omega|_{V(G-X)})$ for $(A,I)$.
  For each $v\in V(G-X)$, since $g_1$ is non-decreasing,
  \begin{align}
    \begin{split}
      |\mu^-(v)|
       & \le \log\omega(G-X) -\log\omega(v)+g_1(|V(G-X)|) \\
       & \le \log\omega(G)-\log\omega(v)+g_1(n) \enspace,
    \end{split}
    \label{eq:mu-}
  \end{align}
  and for each clique $K$ in $G^+-X$, since $g_3$ is non-decreasing,
  \begin{align}
    \label{eq:split}
    \begin{split}
      |\mu^-(K)|
       & \le \log\omega(G-X)-\log\min_{v\in K}\omega(v)+g_3(|V(G-X)|)  \\
       & \le \log\omega(G)-\log\min_{v\in K}\omega(v)+g_3(n) \enspace.
    \end{split}
  \end{align}
  If $G-X$ is empty, no labelling from $\calS$ is needed; take $\mu^-$ and $\kappa^-$ to have empty domains.
  Extend the domain of $\mu^-$ to include the vertices in $X$ by defining $\mu^-(v):=\varepsilon$ for each $v\in X$.

  \paragraph{\boldmath The apex identifier $\sigma(u)$:}
  Fix a linear ordering $u_1,\ldots,u_{b}$ of the vertices in $X$.  By definition, $b\le a$.
  For each $u\in V(G)$, define
  \[
    \sigma(u) :=
    \begin{cases}
      \textrm{$\lceil\log(b+1)\rceil$-bit binary representation of $i$} & \textrm{if $u\in X$ and $u=u_i$,} \\
      \textrm{$\lceil\log(b+1)\rceil$ zero bits}                        & \textrm{if $u\in V(G)\setminus X$.}
    \end{cases}
  \]

  \paragraph{\boldmath The vertex label $\mu(v)$:} For each $v\in V(G)$, let $c(v)$ be a $b$-bit string whose $i$-th bit indicates if $vu_i\in E(G)$.
  For each $v\in V(G)$, define
  \[
    \mu(v):=\langle \sigma(v),c(v),\mu^-(v) \rangle \enspace .
  \]
  For $v\notin X$, \eqref{eq:mu-} bounds $|\mu^-(v)|$; for $v\in X$, $|\mu^-(v)|=0$. Therefore, for each $v\in V(G)$,
  \begin{align*}
    |\sigma(v)|+b+|\mu^-(v)|
     & \le \lceil\log(b+1)\rceil + b +\log\omega(G) - \log\omega(v)+g_1(n) \\
     & \le\log\omega(G) - \log\omega(v)+g_1(n)+2a                                                                            \\
     & \le \log n + g_1(n)+2a+1
     & \textrm{since $\omega$ is nice}                                                                                       \\
     & \in\Oh(a+\log n)\enspace,
    \intertext{and therefore}
    |\mu(v)|
     & \leq |\sigma(v)|+b+|\mu^-(v)| +\Oh(\lgg |\sigma(v)|+\lgg b+\lgg |\mu^-(v)|) & \text{by \eqref{encoding_cost}} \\
     & \leq \log\omega(G) - \log\omega(v)+g_1(n) + \Oh(a+\log\log n) \enspace .
  \end{align*}

  \paragraph{Adjacency testing:}
  We now describe the adjacency tester $A'$.
  Given the labels $\mu(v)=\langle \sigma(v),c(v),\mu^-(v)\rangle$ and $\mu(w)=\langle \sigma(w),c(w),\mu^-(w)\rangle$ for two vertices $v$ and $w$ of $G$,
  $A'$ works as follows.
  First, by checking the values of $\sigma(v)$ and $\sigma(w)$, the tester checks if $v$ and/or $w$ lie in $X$.
  If $\sigma(v)$ encodes an index $i\in[b]$, we return the $i$-th bit of $c(w)$ (which indicates whether $wu_i = wv \in E(G)$).
  Otherwise, if $\sigma(w)$ encodes an index $j\in[b]$, we return the $j$-th bit of $c(v)$ (which indicates whether $vu_j = vw \in E(G)$).
  If neither $v$ nor $w$ is in $X$, then $A'(\mu(v),\mu(w)):=A(\mu^-(v),\mu^-(w))$.

  \paragraph{\boldmath The clique label $\mu(K)$:}
  For each clique $K$ in $G^+$, let $c(K)$ be a $b$-bit string whose $i$-th bit indicates if $u_i\in K$.
  For each clique $K$ in $G^+$, define
  \[
    \mu(K):=
    \begin{cases}
      \langle c(K) \rangle           & \text{if $K\subseteq X$,} \\
      \langle c(K),\mu^-(K\setminus X)\rangle & \text{otherwise.}
    \end{cases}
  \]
  To see that $\mu$ is injective when restricted to cliques, consider two distinct cliques $K$ and $L$ in $G^+$. If $K\cap X \neq L\cap X$, then $c(K)\neq c(L)$, so $\mu(K)\neq \mu(L)$ because their first parts differ. Otherwise, $K\setminus X \neq L\setminus X$. If one is empty, say $K\subseteq X$, then $\mu(K)$ has one part while $\mu(L)$ has two parts, so $\mu(K)\neq\mu(L)$. If both are non-empty, then $\mu^-(K\setminus X)\neq \mu^-(L\setminus X)$ since $\mu^-$ is injective on cliques of $G^+-X$, which implies $\mu(K)\neq \mu(L)$.

  Next, we bound the length of the clique labels. First consider a clique $K\not\subseteq X$. Since $K\setminus X \subseteq K$, we have $\min_{v\in K\setminus X}\omega(v) \ge \min_{v\in K}\omega(v)$, so by~\eqref{eq:split},
  \[
    |\mu^-(K\setminus X)| + b
    \leq \log\omega(G)-\log\min_{v\in K}\omega(v) + g_3(n-|X|)+a \enspace,
  \]
  and therefore
  \begin{align*}
    |\mu(K)|
     & \leq |\mu^-(K\setminus X)|+b + \Oh(\lgg|\mu^-(K\setminus X)|+\lgg b)
     & \text{by \eqref{encoding_cost}}                                                               \\
     & \leq \log\omega(G)-\log\min_{v\in K}\omega(v) + g_3(n) + \Oh(a+\log\log n) \\
     & \leq \log\omega(G)-\log\min_{v\in K}\omega(v) + g'_3(n) \enspace .
  \end{align*}
  If $K\subseteq X$, then $|\mu(K)| \le b + \Oh(\lgg b) \le a + \Oh(\log a) \le g'_3(n) \le \log\omega(G)-\log\min_{v\in K}\omega(v) + g'_3(n)$ since $\omega(G) \ge \min_{v\in K}\omega(v) \ge 1$.

  \paragraph{\boldmath The local identifier $\kappa(K,u)$:}
  For each clique $K$ in $G^+$ and each $u\in K$, define
  \[
    \kappa(K,u) :=
    \begin{cases}
      0 \circ \kappa^-(K\setminus X,u) & \text{if $u\not\in X$,} \\
      1 \circ \sigma(u)                & \text{if $u\in X$.}
    \end{cases}
  \]
  Note that $\kappa(K,u)$ is a concatenation, so it does not incur the encoding cost from \eqref{encoding_cost}. Thus,
  \[
    |\kappa(K,u)|\le 1+ \max\{g_2(n), \lceil\log(a+1)\rceil\} \le g_2(n)+\Oh(\log a) \le g_2'(n) \enspace.
  \]

  \paragraph{Identity testing:}
  We now describe the identity tester $I'$.
  Given $\mu(K)$, $\kappa(K,u)$, and $\mu(v)$ for some clique $K$ in $G^+$, $u\in K$, and $v\in V(G)$,
  we compute $I'(\mu(K),\kappa(K,u),\mu(v))$ as follows.
  First, by examining $\sigma(v)$ and the first bit of $\kappa(K,u)$, the tester determines whether $v\in X$ and whether $u\in X$.
  If $u,v\in X$, then $I'(\mu(K),\kappa(K,u),\mu(v))=1$ if $\kappa(K,u)=1\circ\sigma(v)$, and $I'(\mu(K),\kappa(K,u),\mu(v))=0$ otherwise.
  If $u,v\not\in X$, then since $u\in K$ we have that $K\not\subseteq X$, which means $\mu(K)$ has a second part $\mu^-(K\setminus X)$. The tester extracts this second part and sets
  $I'(\mu(K),\kappa(K,u),\mu(v)) := I(\mu^-(K\setminus X),\kappa^-(K\setminus X,u),\mu^-(v))$.
  Otherwise, one of $u$ or $v$ is in $X$ and the other is not, so $u\neq v$ and we set $I'(\mu(K),\kappa(K,u),\mu(v)) := 0$.

  This completes the description of the adjacency tester $A'$, the identity tester $I'$, and the weighted $(g_1',g_2',g_3')$ mixed labelling $(\mu_\tau,\kappa_\tau):=(\mu,\kappa)$ of $\tau=(G^+,G,\omega)$ for $(A',I')$. Since $\tau$ is an arbitrary nice $\calG^{+a}$-admissible triple, this completes the construction of the weak weighted mixed labelling scheme $(\calS',A',I')$.

  \paragraph{\boldmath $g'_0$-orientability:}
  We now show that $(\calS',A')$ is $g'_0$-orientable.  Fix $n\in\N^+$. Fix an orientation of $U_n$ with maximum outdegree at most $g_0(n)$.  We must show that the $n$-universal graph $U'_n$ defined by $(\calS',A')$ has an orientation in which each vertex has outdegree at most $g'_0(n)$.
  We now explain how to choose the orientation for an arbitrary edge $\{\beta,\beta'\}$ of $U'_n$.   By the definition of $U'_n$, there exists a nice triple $\tau:=(G^+,G,\omega) \in \calT'_{\calG^{+a},\le n}$ and an edge $vv'\in E(G)$ such that the mixed labelling $(\mu,\kappa):=(\mu_\tau,\kappa_\tau)$ of $(G^+,G,\omega)$ for $(A',I')$ has $\mu(v)=\beta$ and $\mu(v')=\beta'$.  Let $X:=X(G^+)$, so $G^+-X\in \cal{G}$. If $G-X$ is non-empty, its restricted weight function remains positive integer-valued, $|V(G-X)|\le n$, and $\omega(G-X)\le\omega(G)\le 2n$. Thus, $\tau^-\in\calT'_{\calG,\le n}$, even though $\tau^-$ need not be nice.  By the definition of the adjacency tester $A'$, exactly one of the following two cases applies:
  \begin{enumerate}
    \item If neither $v$ nor $v'$ is in $X$, then $A'(\beta,\beta')=A(\mu^-(v),\mu^-(v'))=1$.  The triple $\tau^-$ witnesses this edge in $U_n$. Use the fixed orientation of $\{\mu^-(v),\mu^-(v')\}$ in $U_n$. If $\{\mu^-(v),\mu^-(v')\}$ is oriented towards $\mu^-(v)$, then we orient $\{\beta,\beta'\}$ towards $\beta=\mu(v)$. If $\{\mu^-(v),\mu^-(v')\}$ is oriented towards $\mu^-(v')$, then we orient $\{\beta,\beta'\}$ towards $\beta'=\mu(v')$.

    Fix a vertex $\beta$ in $U'_n$.
    We now argue that this branch of the adjacency testing procedure orients at most $g_0(n)\cdot 2^a$ edges away from $\beta$, over all choices of $\beta'$.  The bitstring $\beta$ determines the label $\mu(v)=\langle \sigma(v),c(v),\mu^-(v)\rangle$. Importantly, $c(v)$ is a bit-vector of length exactly $|X|$, so $\beta$ strictly fixes $|X|$. Each edge oriented away from $\beta$ by this branch of the adjacency testing procedure corresponds to an edge oriented away from $\mu^-(v)$ in $U_n$.    Since $\beta$ includes $\mu^-(v)$, there are at most $g_0(n)$ choices of $\mu^-(v')$ that result in an edge being oriented away from $\beta$ in this way.  Since $v'\not\in X$, $\sigma(v')$ is a sequence of $0$s of length $\lceil \log (|X|+1)\rceil$, and $c(v')$ is a bit-vector of length $|X|$, leaving $2^{|X|}$ choices for $c(v')$.  Thus, the number of possible strings $\beta'=\langle \sigma(v'),c(v'),\mu^-(v')\rangle$ for which the edge $\{\beta,\beta'\}$ is oriented away from $\beta$ by this branch of the adjacency testing procedure is at most $g_0(n) 2^{|X|} \leq g_0(n) 2^{a}$.

    \item Otherwise, at least one of $v$ or $v'$ is in $X$.  If both $v,v'\in X$, then we orient $\{\beta,\beta'\}$ towards the lexicographically larger string.  If $\{v,v'\}\cap X=\{v\}$, then we orient $\{\beta,\beta'\}$ towards $\beta$. If $\{v,v'\}\cap X=\{v'\}$, then we orient $\{\beta,\beta'\}$ towards $\beta'$.

    Fix a vertex $\beta$ in $U'_n$.  We now argue that this branch of the adjacency testing procedure orients at most $a 2^{a}$ edges away from $\beta$. This corresponds to a case in which $v'\in X$, so $\beta'=\langle \sigma(v'),c(v'),\varepsilon\rangle$. There are at most $|X|$ choices for $\sigma(v')$ and $2^{|X|}$ choices for $c(v')$. Again, $|X|$ is fixed by the length of $c(v)$ in $\beta$, so the number of strings $\beta'$ oriented away from $\beta$ is at most $|X|\cdot 2^{|X|} \leq a 2^{a}$.
  \end{enumerate}
  Therefore, for any fixed $\beta\in V(U'_n)$, the outdegree of $\beta$ in this orientation of $U'_n$ is at most $g_0(n) 2^{a} + a 2^{a} = (g_0(n)+a)2^{a} \leq g_0'(n)$.

  \paragraph{\boldmath Strong $g'_0$-orientability:}  We now verify that $(\calS',A',I')$ is strongly $g'_0$-orientable.
  Let $\tau_0 := (G_0^+, G_0, \omega_0) \in \calT'_{\calG^{+a}, \le n}$ be a nice triple, and let $(\mu_{\tau_0},\kappa_{\tau_0})\in\calS'$ be the mixed labelling of $\tau_0$ for $(A',I')$. Let $K_0$ be a clique in $G_0^+$ and let $u_0\in K_0$. Let $\zeta:=\mu_{\tau_0}(K_0)$ and let $\gamma:=\kappa_{\tau_0}(K_0,u_0)$. We need to show that there are at most $g'_0(n)$ strings $\beta\in V(U'_n)$ such that there exists a nice triple $\tau:=(G^+,G,\omega) \in \calT'_{\calG^{+a},\le n}$, some clique $K$ in $G^+$, and some $u\in K$ such that $\mu_\tau(K)=\zeta$, $\kappa_\tau(K,u)=\gamma$, and $\mu_\tau(u)=\beta$.

  Let $\tau:=(G^+,G,\omega)$, $K$, $u$, and $\beta$ be as described above and let $X:=X(G^+)$.  Let $(\mu,\kappa):=(\mu_\tau,\kappa_\tau)$. Since $\zeta = \mu(K)$, either $\zeta = \langle c(K) \rangle$ or $\zeta = \langle c(K), \mu^-(K\setminus X)\rangle$. In both cases, $c(K)$ is a bit-vector of length $|X|$, so $\zeta$ fixes $|X|$.
  \begin{itemize}[nosep,nolistsep]
    \item If $u\in X$, then $\beta=\mu(u)=\langle\sigma(u),c(u),\varepsilon\rangle$. Since $\gamma = \kappa(K,u) = 1 \circ \sigma(u)$, the string $\gamma$ determines $\sigma(u)$. There are $2^{|X|}$ choices for $c(u)$. Because $|X|$ is fixed by $\zeta$, there are at most $2^{|X|} \le 2^a$ such strings $\beta$.
    \item If $u\in V(G-X)$, then $\beta=\langle \sigma(u),c(u),\mu_{\tau^-}^-(u)\rangle$, and $I(\mu_{\tau^-}^-(K\setminus X),\kappa_{\tau^-}^-(K\setminus X,u),\mu_{\tau^-}^-(u))=1$.  The string $\mu_{\tau^-}^-(K\setminus X)$ is extracted from $\zeta$ and the string $\kappa^-_{\tau^-}(K\setminus X,u)$ is extracted from $\gamma$. Here $\tau^-\in\calT'_{\calG,\le n}$, since $|V(G-X)|\le n$ and $\omega(G-X)\le\omega(G)\le 2n$. The triple $\tau^-$, the clique $K\setminus X$, and the vertex $u$ witness the occurrence of these three strings in $\calS$. Since $(\calS,A,I)$ is strongly $g_0$-orientable, there are at most $g_0(n)$ choices for $\mu_{\tau^-}^-(u)$ given $\zeta$ and $\gamma$.  Since $u\not\in X$, $\sigma(u)$ is a string of $\lceil\log(|X|+1)\rceil$ zeros.  The bit-vector $c(u)$ is an arbitrary binary string of length $|X|$. Because $|X|$ is fixed by $\zeta$, there are at most $g_0(n)\cdot 2^{|X|} \le g_0(n) 2^a$ choices for $\beta=\mu_\tau(u)$ arising from such witnesses.
  \end{itemize}
  Therefore, the number of possible choices for $\beta$ is at most $2^a+ g_0(n) 2^{a} \le (2g_0(n)+a)2^{a} = g'_0(n)$.

  Finally, extend $(\calS',A',I')$ to all $\calG^{+a}$-admissible triples by assigning to $(G^+,G,\omega)$ its labelling for the normalized weights $\widehat\omega(v):=\ceil{|V(G)|\omega(v)/\omega(G)}$. By \cref{lem:nice_weights}, the vertex and clique label-length bounds increase by at most $1$, which is absorbed into the stated $\Oh(a)$ terms. By \cref{obs:normalization_orientability}, the resulting full weighted scheme remains strongly $g'_0$-orientable. \qedhere
\end{proof}

\subsection{Tree-Decompositions: Combining via Clique-Sums}
\begin{lem}\label{lem:the_big_one}
  Let $\calG$ be a monotone class of graphs that admits a strongly $g_0$-orientable weighted $(g_1,g_2,g_3)$ mixed labelling scheme where $g_i:\N\to\R^{\ge 0}$ is a non-decreasing function for each $i\in\{0,1,2,3\}$ and $g_i(n)\in\Oh(\log n)$ for each $i\in\{1,2,3\}$. Let $b:\N\to\{2,3,\dots\}$ and $k:\N\to\N^+$ be non-decreasing functions such that $b(n)\le n/2$ and $k(n)\le \log n$ for all sufficiently large integers $n$.  Let $\calG'$ be the graph class that consists of, for each $n\in\N$, every $n$-vertex graph $G$ that has a rooted tree-decomposition of adhesion-width at most $k(n)$ whose torsos belong to $\calG$. 
  Then $\calG'$ admits an injective $g'_0$-orientable $f(n)$-bit adjacency labelling scheme with
  \begin{align*}
    g_0'(n) & = g_0(n)\cdot \exp(\Oh(k(n)(g_3(n)+g_2(n)+k(n)+\log b(n)+\log\log n))) \enspace, \\
    f(n) & = \log n + g_1(n) + \Oh((g_3(n)+k(n)+\log\log n)\cdot \log_{b(n)} n) \\
    & \quad + \Oh(k(n)\cdot (g_3(n)+g_2(n) + k(n) + \log b(n) + \log\log n )) \enspace .
  \end{align*}
\end{lem}

\begin{proof}
  Let $g_0,g_1,g_2,g_3$, $\calG$, $\calG'$, $k$, and $b$ be as given in the statement of the lemma. For any fixed constant $n_0$, graphs on at most $n_0$ vertices have labels of length $\Oh(1)$ and $U_n'$ has $\Oh(1)$ vertices and edges, so the bounds hold trivially for $n < n_0$. Thus, we may assume that $n$ is sufficiently large so that $b(n)\le n/2$ and $k(n)\le \log n$.
  We will define an adjacency tester $A'$ and an adjacency labelling $\mathdefin{\calA'}:=(\ell_{G_0}:G_0\in\calG')$ of $\calG'$ for $A'$.  Since it is impractical to describe the operation of $A'$ without first describing the format of the labels given as arguments to $A'$, we begin by describing the construction of an adjacency labelling $\ell_{G_0}:V(G_0)\to\{0,1\}^*$ for an arbitrary graph $G_0\in\calG'$.

  Let $G_0$ be an $n$-vertex graph in $\calG'$ and let $\calT_0$ be a rooted tree-decomposition of $G_0$ witnessing that $G_0\in\calG'$, so $\calT_0$ has adhesion-width at most $k:=k(n)$ and each torso of $\calT_0$ belongs to $\calG$. Let $b:=b(n)$.  Let $G$, $\calT:=(T,(B_x\mid x\in V(T)))$, and $Y\subseteq E(T)$ be the result of applying \cref{all_singing_all_dancing} to the tree-decomposition $\calT_0$ of $G_0$, so that $G$ has an isolated vertex $v_0$ with $G-v_0=G_0$, $B_{x_0}=\{v_0\}$, and $\calT$ is a tidy rooted tree-decomposition of $G$ with root $x_0$.  
  
  Define $A_{x_0}:=\emptyset$ and, for each $x\in V(T)\setminus\{x_0\}$, let $A_x$ be the parent adhesion of $x$ in $\calT$.  For each $x\in V(T)$, $\torso{G}{\calT,B_x}-A_x$ is either a $1$-vertex graph (if $x=x_0$), can be made into a subgraph of a torso of $\calT_0$ by deleting at most one adhesion of $\calT_0$ (so at most $k$ vertices), or its vertex set is contained in the union of at most $3$ adhesions of $\calT_0$ (so it has at most $3k$ vertices). Since $\calG$ is monotone, each subgraph of a graph in $\calG$ is in $\calG$. In each case, therefore, $\torso{G}{\calT,B_x}-A_x\in\calG^{+3k}$.

  For each $v\in V(G)$, let $\mathdefin{\omega(v)}:=1$ and note that $\omega:V(G)\to\{1\}$ is a nice weight function. For each $x\in V(T)$, let $\mathdefin{\omega(x)}:=\omega(B_x\setminus A_x)$. Then $\omega(T)=\sum_{x\in V(T)}\omega(x)=\omega(G)=n+1$, since $\{B_x\setminus A_x\mid x\in V(T)\}$ is a partition of $V(G)$.  For each $x\in V(T)$, let $\mathdefin{C_x}:=N_{T_x}(x)$ be the set of $T$-children of $x$. For each $x\in V(T)$, and each clique $K$ in $\torso{G}{\calT,B_x}-A_x$, let $\mathdefin{C_K}:=\{y\in C_x\mid A_y\setminus A_x=K\}$.  Since $\calT$ is tidy, $A_y\setminus A_x\neq\emptyset$ for each $y\in C_x$.  Therefore, each $y\in C_x$ is contained in $C_K$ for exactly one clique $K=A_y\setminus A_x$ in $\torso{G}{\calT,B_x}-A_x$. We say that $K$ is \defin{light} if $\sum_{y\in C_K}\omega(T_y) < \omega(T_x)/b$, and is \defin{non-light} otherwise.

  Let $F:=T-Y$.  Properties 
  \ref{req:root_bag}, 
  \ref{req:height}, 
  \ref{req:subtree_size}, 
  \ref{req:same_size}, 
  \ref{req:adhesion_width}, and
  \ref{req:subtree_torso_apex} 
  of \cref{all_singing_all_dancing} translate to the following properties of $\calT$, $Y$, and $F$:
  \begin{enumerate}[nosep,nolistsep,label=(P\arabic*)]
    \item\label{td_width} $\calT$ has adhesion-width at most $2k$;
    \item\label{torsos_in_G_plus_k} for each $x\in V(T)$, 
    $\torso{G}{\calT,B_x}-A_x\in \calG^{+3k}$;
    \item\label{component_size} each component of $F$ has at most $4b-1$ nodes;
    \item\label{no_fast_shrink} for any two nodes $x,y$ in the same component of $F$, $|\log\omega(T_x)-\log\omega(T_y)|\le\log b$;
    \item\label{light_clique_buster} for each $x\in V(F)$, each light clique $K$ in $\torso{G}{\calT,B_x}-A_x$, and each $y\in C_K$, $xy\in Y$; and
    \item\label{y_on_path} each root-to-leaf path $P$ in $T$ contains at most $2\log_b n +\Oh(1)$ edges in $Y$.
  \end{enumerate}
  (Property~\ref{light_clique_buster} holds because if $xy\notin Y$, then $x$ and $y$ lie in the same component of $F$, so by~\ref{req:same_size} we would have $\omega(T_y)\ge \omega(T_x)/b$, which contradicts the definition of a light clique since $\omega(T_y)\le \sum_{z\in C_K}\omega(T_z) < \omega(T_x)/b$.)

  By \cref{add_apexes}, $\calG^{+3k}$ admits a strongly $h_0$-orientable weighted $(h_1,h_2,h_3)$ mixed labelling scheme $\mathdefin{\calS}:=((\mu_\tau,\kappa_\tau):\tau\in\calT_{\calG^{+3k}})$ with 
  \begin{align}
    h_0(n) & = (2g_0(n) + 3k)2^{3k} \enspace , \nonumber \\
    h_1(n) & = g_1(n) + \Oh(k(n)+\log\log n) \enspace , \label{h1_td} \\ 
    h_2(n) & = g_2(n) + \Oh(\log k(n) + \log\log n) \enspace , \label{h2_td} \\ 
    h_3(n) & = g_3(n) + \Oh(k(n) + \log\log n) \enspace \label{h3_td} . 
  \end{align} 
  For each $n\in\N^+$, let $U_n$ be the $n$-universal graph defined by $(\calS,A)$.

  \paragraph{\boldmath The weight function $\delta$:} 

  We now define a weight function $\delta:V(G)\to\R^+$ and, for each $x\in V(T)$, we define $\delta(x)=\delta(B_x\setminus A_x)$.  This weight function is designed so that, for any node $x\in V(T)$, $\tfrac{1}{4b}\omega(T_x)\le \delta(x)\le 3k\omega(T_x)$.  In fact, the upper bound will have an even stronger property:  For each component $C$ of $F$, $\delta(C)\le 6k\omega(T_{r(C)})$.

  For each $x\in V(T)$, let $\mathdefin{C^-_x}:=\{y\in C_x\mid xy\in Y\}$ and let $\mathdefin{C^+_x}:=C_x\setminus C^-_x$.  For each $x\in V(T)$ and each clique $K$ in $\torso{G}{\calT,B_x}-A_x$, let $\mathdefin{C^-_K}:=C_K\cap C^-_x$. For each $x\in V(T)$ and each $v\in B_x\setminus A_x$, let $\mathdefin{C_v}:=\{y\in C_x\mid v\in A_y\}$, let $\mathdefin{C^-_v}:=C_v\cap C^-_x$ and let $\mathdefin{C^+_v}:=C_v\cap C^+_x$.  
  
  For each $x\in V(T)$ and each $v\in B_x\setminus A_x$, define 
  \[
    \mathdefin{\delta(v)} := \omega(v)+\sum_{y\in C^-_v} \omega(T_y) + \tfrac{1}{4b}\cdot\sum_{y\in C^+_v}\omega(T_y) \enspace .
  \]
  For each $x\in V(T)$, let 
  \[
    \mathdefin{\delta(x)}:=\delta(B_x\setminus A_x) \enspace . 
  \]
  The definition of $\delta(v)$ is chosen so that $\delta(x)$ approximates $\omega(T_x)$ for each $x\in V(T)$, in a manner we now make precise.
  For each $x\in V(T)$, $y\in C_x$, and $v\in B_x\setminus A_x$, $y\in C_v$ if and only if $v\in A_y\setminus A_x$.  By property~\ref{td_width} of $\calT$, we have $|A_y|\le 2k\le 3k$ for each $y\in V(T)$, so each $y\in V(T)$ belongs to at most $3k$ sets $C_v$. Thus, 
  \begin{equation}
    \delta(x) 
    \le  \sum_{v\in B_x\setminus A_x}\left(\omega(v)+ \sum_{y\in C_v} \omega(T_y) \right) 
    \le \omega(x) + \sum_{y\in C_x} 3k\cdot\omega(T_y) \le 3k\cdot\omega(T_x) \enspace .
    \label{delta_x_upper_bound}
  \end{equation}
  Since $\calT$ is a tidy tree-decomposition, $A_y\setminus A_x\neq\emptyset$ for each $y\in C_x$, and thus every $y \in C_x$ appears in $C_v$ for at least one $v \in B_x \setminus A_x$.  Therefore
  \begin{equation}
    \delta(x) = \sum_{v\in B_x\setminus A_x}\delta(v) \ge \omega(x) + \sum_{y\in C^-_x} \omega(T_y) + \tfrac{1}{4b}\sum_{y\in C^+_x} \omega(T_y)
    \ge \tfrac{1}{4b}\omega(T_x) \enspace \label{delta_x_lower_bound}
  \end{equation}
  for each $x\in V(T)$.

  Beyond this, the definition of $\delta$ is chosen so that $\delta(C)$ approximates $\omega(T_{r(C)})$ for each component $C$ of $F$:
  \begin{align}
    \delta(C) 
    & = \sum_{x\in V(C)}\sum_{v\in B_x\setminus A_x}\left(\omega(v)+\sum_{y\in C^-_v} \omega(T_y) 
      + \tfrac{1}{4b}\sum_{y\in C^+_v}\omega(T_y)\right) \nonumber \\
    & \le \sum_{x\in V(C)}\left(\omega(x)+\sum_{y\in C^-_x} 3k\cdot\omega(T_y) 
      + \tfrac{1}{4b}\cdot\sum_{y\in C^+_x}3k\cdot\omega(T_y)\right) 
      & \text{(since $|A_y|\le 3k$)} \nonumber\\
    & \le 3k\sum_{x\in V(C)}\left(\omega(x)+\sum_{y\in C^-_x} \omega(T_y)\right) 
      + \tfrac{3k}{4b}\sum_{x\in V(C)}\sum_{y\in C^+_x}\omega(T_y) 
      & \text{(since $k\ge 1$)}\nonumber\\
    & \le 3k\sum_{x\in V(C)}\left(\omega(x)+\sum_{y\in C^-_x} \omega(T_y)\right) 
      + \tfrac{3k}{4b}\sum_{x\in V(C)}\omega(T_{x}) 
      & \text{(since $\textstyle T_x\supset\bigcup_{y\in C^+_x}T_y$)}\nonumber\\
    & \le 3k\sum_{x\in V(C)}\left(\omega(x)+\sum_{y\in C^-_x} \omega(T_y)\right) 
      + \tfrac{3k}{4b}\sum_{x\in V(C)}\omega(T_{r(C)}) 
      & \text{(since $T_{r(C)}\supseteq T_x$)}\nonumber\\
    & \le 3k\sum_{x\in V(C)}\left(\omega(x)
      + \sum_{y\in C^-_x} \omega(T_y)\right) + 3k\cdot \omega(T_{r(C)}) 
      & \text{(since $|V(C)|<4b$, by \ref{component_size})} \nonumber \\
    & \le 6k\cdot\omega(T_{r(C)})
    \enspace . \label{delta_c_upper_bound}
  \end{align}
  where the last inequality holds because the sets in $\{V(C)\} \cup \{ V(T_y)\mid x\in V(C), y\in C^-_x\}$ are pairwise disjoint subsets of $V(T_{r(C)})$. The final inequality is most conveniently written as 
  \begin{equation}  
    \omega(T_{r(C)}) \ge \delta(C) / (6k) \enspace . \label{omega_trC_lower_bound}
  \end{equation}  

  \paragraph{\boldmath The inner-tree codes $\rho_C$ for components $C$ of $F$:}
  For each component $C$ of $F$, let $\mathdefin{\rho_C}:V(C)\to\{0,1\}^*$ be the prefix-free code obtained by applying \cref{shannon_fano_elias_code} to the set $S=V(C)$ with weight function $\delta|_S$ obtained by restricting $\delta$ to $S$.  Then, for each $x\in V(C)$:
  \begin{align}
      |\rho_C(x)| 
      & \le \log \delta(C) - \log \delta(x) + \Oh(1) \label{rho_C_bound_defn} \\
      & \le \log (6k\cdot\omega(T_{r(C)})) - \log \delta(x) + \Oh(1) 
        & \text{(by \eqref{omega_trC_lower_bound})} \nonumber \\
      & = \log\omega(T_{r(C)}) - \log\delta(x) + \Oh(\log\log n) 
        & \text{(since $k\in \Oh(\log n)$)} \nonumber \\
      & \le \log\omega(T_{r(C)}) - \log(\omega(T_x)/(4b)) + \Oh(\log\log n) 
        & \text{(by \eqref{delta_x_lower_bound})} \nonumber \\
      & \le \log\omega(T_{r(C)}) - \log(\omega(T_{r(C)})/(4b^2)) + \Oh(\log\log n) 
        & \text{(by property~\ref{no_fast_shrink} of $F$)} \nonumber \\
      & \le 2\log b + \Oh(\log\log n)  .  \label{rho_C_bound_2logb}
  \end{align}

  \paragraph{\boldmath The inner-node labelling $\mu_x$:} The next part of our labelling involves a mixed labelling $(\mu_x,\kappa_x)$ of $\tau_x:=(\torso{G}{\calT,B_x}-A_x,G[B_x]-A_x,\widehat\delta_x)$ for each $x\in V(T)$.  Here $\widehat\delta_x$ is a nice weight function on $B_x\setminus A_x$ that is closely related to $\delta|_{B_x\setminus A_x}$. The only reason for the (distracting) introduction of $\widehat\delta_x$ is that later, we will make use of the fact that $(\calS,A,I)$ is strongly $h_0$-orientable.  In particular, the $n$-universal graph $U_n$ defined by $(\calS,A)$ has an orientation of maximum outdegree $h_0(n)$. For this to be useful, we require that the label $\mu_x(v)$ assigned to each vertex $v\in B_x\setminus A_x$ correspond to a vertex in $U_n$, and this is only guaranteed if $\tau_x\in \calT'_{\calG^{+3k},\le n}$. This is why we require $\widehat\delta_x$ to be nice.
  
  For each $x\in V(T)$, let $m_x:=|B_x\setminus A_x|$ and define the weight function $\widehat\delta_x:B_x\setminus A_x\to\N^+$ by
  \[
    \widehat\delta_x(v):=\left\lceil\frac{m_x\delta(v)}{\delta(x)}\right\rceil\enspace.
  \]
  As in the proof of \cref{obs:nice}, $\widehat\delta_x$ is nice, $\widehat\delta_x(B_x\setminus A_x)<2m_x$, and
  \[
    \log\widehat\delta_x(B_x\setminus A_x)-\log\widehat\delta_x(v)
    \le\log\delta(x)-\log\delta(v)+1\enspace
  \]
  for each $v\in B_x\setminus A_x$.
  Let $\tau_x:=(\torso{G}{\calT,B_x}-A_x,G[B_x]-A_x,\widehat\delta_x)$ and let $\mathdefin{(\mu_x,\kappa_x)}:=(\mu_{\tau_x},\kappa_{\tau_x})$. The root bag has $m_{x_0}=1$, and every other set $B_x\setminus A_x$ excludes $v_0$, so $m_x\le n$ for each $x\in V(T)$. Thus, $\tau_x\in\calT'_{\calG^{+3k},\le n}$.  Increasing the $\Oh(\log\log n)$ terms in $h_1$ and $h_3$ by $1$ preserves \cref{h1_td,h3_td}, so we obtain, for each $x\in V(T)$,

  \begin{equation}
    |\mu_x(v)|\le \log\delta(x)-\log\delta(v) + h_1(n)
      \label{mu_x_v_upper_bound}
  \end{equation}
  for each $v\in B_x\setminus A_x$ and 
  \begin{equation}
    |\mu_x(K)|\le \log\delta(x)-\log(\min_{u\in K} \delta(u)) + h_3(n)
      \label{mu_x_K_upper_bound}
  \end{equation}
  for each clique $K$ in $\torso{G}{\calT,B_x}-A_x$.  
  
  The definition of $\delta$ has been carefully chosen so that every non-light clique receives a short label: If $K$ is a non-light clique then 
  \begin{align}
    \min_{u\in K} \delta(u) 
    & \ge \min_{u\in K}\left(\sum_{y\in C_u}\tfrac{1}{4b}\omega(T_y) \right) 
    & \text{(by definition of $\delta$)} \nonumber \\
    & \ge \tfrac{1}{4b}\sum_{y\in C_K} \omega(T_y)  
    & \text{(since $C_K\subseteq C_u$ for all $u\in K$)} \nonumber \\
    & \ge \omega(T_x)/(4b^2) 
    & \text{(since $K$ is not light).} \label{nonlight_clique_delta_lower_bound}
  \end{align}

  Therefore, for any non-light clique $K$ in $\torso{G}{\calT,B_x}-A_x$,  
  \begin{align}
    |\mu_x(K)|
    & \le h_3(n) + \log\delta(x)-\log(\min_{u\in K} \delta(u)) \nonumber \\
    & \le h_3(n) + \log\delta(x)-\log(\omega(T_x)/(4b^2)) & \text{(by \eqref{nonlight_clique_delta_lower_bound})} \nonumber \\
    & \le h_3(n) + \log(3k\cdot\omega(T_x))-\log(\omega(T_x)/(4b^2)) \nonumber \\
    & \le h_3(n) + 2\log b  + \Oh(\log\log n) \enspace .
      \label{mu_x_K_light_upper_bound}
  \end{align}

  \paragraph{\boldmath The inter-tree codes $\rho_K$:}

  For each $x\in V(T)$ and each clique $K$ in $\torso{G}{\calT,B_x}-A_x$ such that $C^-_K\neq\emptyset$, let $\mathdefin{\rho_K}:C^-_K\to\{0,1\}^*$ be obtained by applying \cref{shannon_fano_elias_code} to the set $S=C^-_K$ with weight defined by $y\mapsto\omega(T_y)$.  Then, for each $y\in C^-_K$, 
  \begin{align}
    |\rho_K(y)|
      & \le 
      \textstyle\log\left(\sum_{y'\in C^-_K}\omega(T_{y'})\right) - \log\omega(T_y) + \Oh(1) \nonumber \\
      & \le \log(\min_{u\in K}\delta(u)) - \log\omega(T_y) + \Oh(1) \enspace ,
      \label{rho_K_y_upper_bound}
  \end{align}
  since $y'\in C^-_u$ for each $u\in K$.

  \paragraph{\boldmath The outer-tree codes $\tau(C,D)$:}
  Let $C$ and $D$ be two components of $F$ such that some edge $xy\in Y$ joins them, with $x\in V(C)$, $y\in V(D)$ and $x=\parent_T(y)$. Thus $y=r(D)$ is the root of component $D$. Then we define
  \[
     \mathdefin{\tau(C,D)}:= \langle\rho_C(x), \mu_x(A_y\setminus A_x), \rho_{A_y\setminus A_x}(y)\rangle \enspace .
  \]
  Then
  \begin{align}
    |\tau(C,D)|
      & = |\rho_C(x)| + |\mu_x(A_y\setminus A_x)| + |\rho_{A_y\setminus A_x}(y)| + \Oh(\log\log n) \nonumber \\
    & \le \log\delta(C) - \log\delta(x) + |\mu_x(A_y\setminus A_x)| + |\rho_{A_y\setminus A_x}(y)| + \Oh(\log\log n) 
      & \text{(by \cref{rho_C_bound_defn})} \nonumber \\
    &  \le \log\delta(C) - \log\left(\min_{u\in A_y\setminus A_x}\delta(u)\right) + h_3(n) + |\rho_{A_y\setminus A_x}(y)| + \Oh(\log\log n) 
      & \text{(by \cref{mu_x_K_upper_bound})} \nonumber \\
    & \le \log\delta(C) - \log\omega(T_y) + h_3(n) + \Oh(\log\log n)   
      & \text{(by \cref{rho_K_y_upper_bound})} \nonumber \\
    & \le \log\delta(C) - \log(\delta(D)/(6k)) + h_3(n) + \Oh(\log\log n)   
      & \text{(by \cref{omega_trC_lower_bound})} \nonumber \\
    & \le \log\delta(C) - \log\delta(D) + h_3(n) + \Oh(\log\log n) 
      & \text{\hspace{-1in}(since $k\le\log n$)}. \label{tau_upper_bound}
    \end{align}
    For a fixed component $C$ of $F$, $\tau(C,\cdot)$ is a function that maps each component $D$ of $F$ whose root is adjacent to a node of $C$ onto a binary string.  For each component $C$ of $F$, the mapping $\tau(C,\cdot)$ is injective, because of the following: $\rho_C:V(C)\to\{0,1\}^*$ is injective by definition; $\mu_{x}$ is injective over the cliques in $\torso{G}{\calT,B_x}-A_x$ for each $x\in V(C)$; and $\rho_{K}:C^-_K\to\{0,1\}^*$ is injective for each clique $K$ in $\torso{G}{\calT,B_x}-A_x$.

  \paragraph{\boldmath The node id $\id(x)$:}

  Refer to \cref{hopper}. Let $x\in V(T)$ and let $x_1y_1,\ldots,x_sy_{s}$ be the sequence of edges in $Y$ encountered along the path from the root $x_0$ of $T$ to $x$. Let $x_{s+1}:=x$ and let $y_0:=x_0$ be the root of $T$. 
  \begin{figure}
    \centering
    \includegraphics[page=2,width=0.5\textwidth]{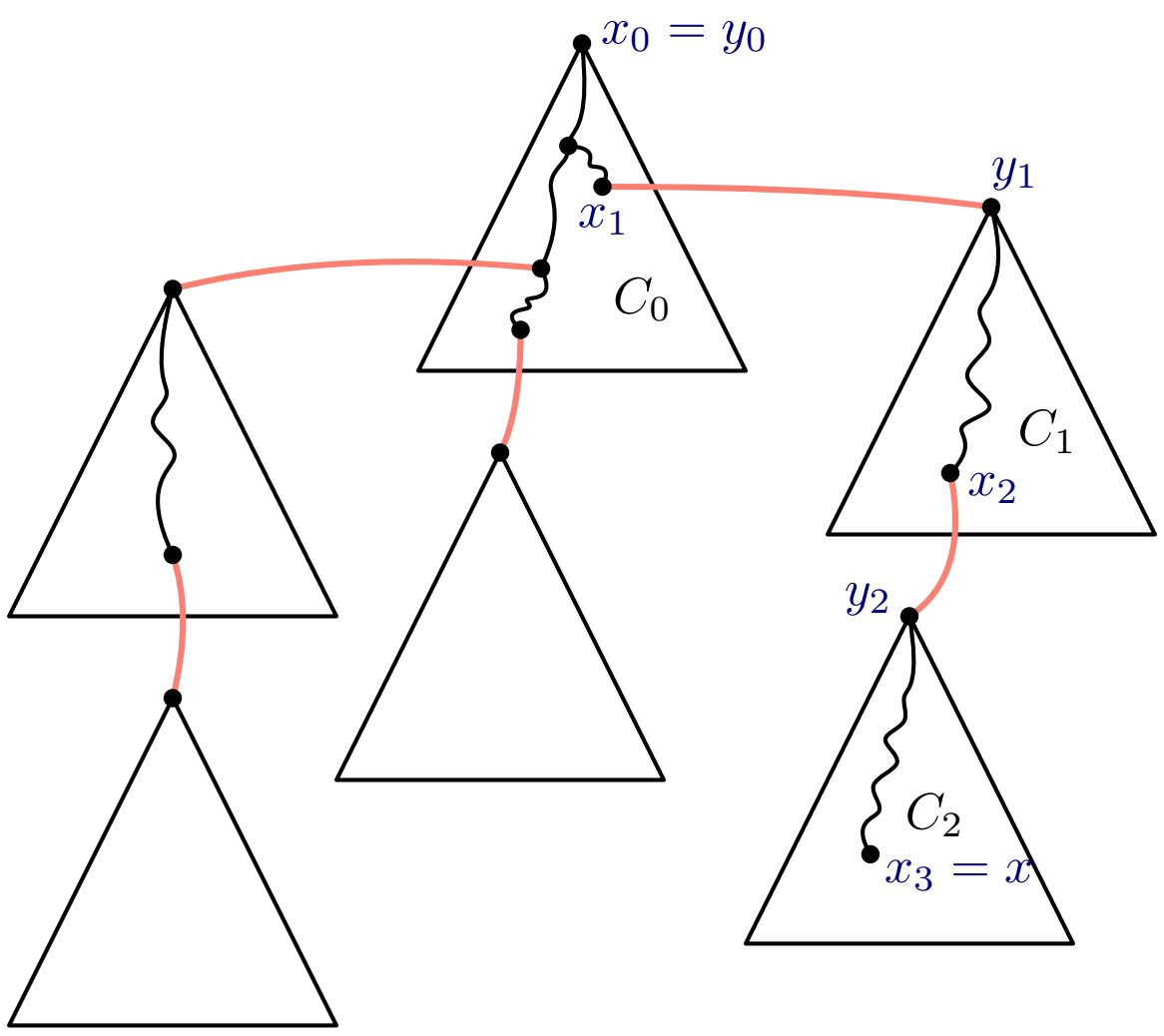}
    \caption{The tree $T$ with edges in $Y$ coloured and an example of the sequences $y_0, x_1,y_1,x_2,y_2,x_3$ and $C_0,C_1,C_2$ used to define $\id(x)=\langle\tau(C_0,C_1), \tau(C_1,C_2), \rho_{C_2}(x) \rangle$.}
    \label{hopper}
  \end{figure}
  Then $x_iy_i\in Y$ for each $i\in \{1,\ldots,s\}$, and $y_i$ is the root of some component $C_i$ of $F$ for each $i\in \{0,\ldots,s\}$.  Define
  \begin{align*}
    \mathdefin{\id(C_s)} & := \langle\tau(C_0,C_1),\tau(C_1,C_2),\ldots,\tau(C_{s-1},C_s)\rangle \enspace , \\
    \mathdefin{\id(x)} & := \langle\tau(C_0,C_1),\tau(C_1,C_2),\ldots,\tau(C_{s-1},C_s),\rho_{C_s}(x)\rangle \enspace .
  \end{align*}
  Then 
  \begin{align}
    |\id(C_s)| 
    & \le \sum_{i=0}^{s-1}(|\tau(C_i,C_{i+1})| + \Oh(\log\log n)) \nonumber \\
    & \le \sum_{i=0}^{s-1} \left( \log\delta(C_i)-\log\delta(C_{i+1}) \right) + s\cdot h_3(n)+\Oh((s+1)\log\log n) 
      & \text{(by \cref{tau_upper_bound})} \nonumber \\
    & = \log\delta(C_0) - \log\delta(C_s) + s\cdot h_3(n)+\Oh((s+1)\log\log n) \nonumber \\
    & = \log\delta(C_0) - \log\delta(C_s) + \Oh((h_3(n)+\log\log n)\cdot\log_b n) & \text{(by \ref{y_on_path})} \nonumber \\
    & \le \log\omega(G) - \log\delta(C_s) + \Oh((h_3(n)+\log\log n)\cdot\log_b n) & \text{(by \cref{delta_c_upper_bound})}  \label{id_C_upper_bound}
  \end{align}
  and
  \begin{align}
    |\id(x)| & = |\id(C_s)| + |\rho_{C_s}(x)| + \Oh(\log\log n) \nonumber \\
    & \le \log\omega(G) - \log\delta(x) + \Oh((h_3(n)+\log\log n)\cdot\log_b n)  
      & \text{(by \cref{id_C_upper_bound} and \cref{rho_C_bound_defn})} \enspace . \label{id_x_bound}
  \end{align}
  It is straightforward to verify (by induction on $s$) that the function $\id:V(T)\to\{0,1\}^*$ is injective. (With a slight abuse of notation, for $i\le s$ we say that $\id(C_i)$ is a \emph{prefix} of $\id(C_s)$ and $\id(x)$, and write $\id(C_i)=\id(C_{i-1})\circ\tau(C_{i-1},C_i)$ and $\id(x)=\id(C_s)\circ\rho_{C_s}(x)$ to denote the multipart encodings obtained by appending one part.)

  \paragraph{\boldmath The vertex label $\ell(v)$:}

  For each $x\in V(T)$ and each $v\in B_x\setminus A_x$, define
  \[
    \mathdefin{\ell(v)} := \langle \id(x), \mu_x(v),\alpha(v)\rangle \enspace .
  \]
  The parts $\id(x)$ and $\mu_x(v)$ are defined above. We now describe $\alpha(v)$, which contains an entry $\alpha(v,u)$ for each $u\in N_G(v)\cap A_x$.

  Let $u\in N_G(v)\cap A_x$.  Then $u\in B_{x'}\setminus A_{x'}$ for some $T$-ancestor $x'\neq x$ of $x$. Let $y'$ be the child of $x'$ on the path from $x'$ to $x$. The idea is that $\alpha(v,u)$ (in combination with $\id(x)$) should store enough information to reconstruct $\id(x')$, $\mu_{x'}(A_{y'}\setminus A_{x'})$, and $\kappa_{x'}(A_{y'}\setminus A_{x'},u)$.  Let $C'$ be the component of $F$ that contains $x'$ and let $D'$ be the component of $F$ that contains $y'$.  Observe that some prefix of $\id(x)$ contains $\id(C')$.  Furthermore, we can reconstruct $\id(x')$ by appending $\rho_{C'}(x')$ to $\id(C')$.  Let \defin{$\idx(C')$} be the number of edges in $Y$ on the path from the root of $T$ to $r(C')$. There are two cases to consider for the edge $x'y'$.   
  \begin{itemize}
    \item If $x'y'\in Y$, then $C'\neq D'$.  Then $\id(x)$ contains $\id(D')$ as a prefix, and $\id(D')$ contains $\tau(C',D')=\langle \rho_{C'}(x'), \mu_{x'}(A_{y'}\setminus A_{x'}), \rho_{A_{y'}\setminus A_{x'}}(y') \rangle$.  In this case, the only extra information needed in $\alpha(v,u)$ is $\idx(C')$ and $\kappa_{x'}(A_{y'}\setminus A_{x'}, u)$.  

    \item If $x'y'\not\in Y$, then $C'=D'$. In this case, $\id(x)$ contains $\id(C')$ but does not contain $\id(x')$, since it does not contain $\rho_{C'}(x')$.  In this case, $\alpha(v,u)$ will store $\rho_{C'}(x')$ and $\mu_{x'}(A_{y'}\setminus A_{x'})$ in addition to $\idx(C')$ and $\kappa_{x'}(A_{y'}\setminus A_{x'}, u)$.  We now explain why these strings are sufficiently short.
    
    By \cref{rho_C_bound_2logb}, $|\rho_{C'}(x')|\in \Oh(\log b + \log\log n)$.  More importantly, since $C'=D'$, $x'y'\not\in Y$, so $A_{y'}\setminus A_{x'}$ is a non-light clique in $\torso{G}{\calT,B_{x'}}-A_{x'}$, by \ref{light_clique_buster}. Therefore, $|\mu_{x'}(A_{y'}\setminus A_{x'})|\le h_3(n)+ \Oh(\log b +\log\log n)$, by \cref{mu_x_K_light_upper_bound}.  This allows us to store both $\rho_{C'}(x')$ and $\mu_{x'}(A_{y'}\setminus A_{x'})$ in $\alpha(v,u)$.  
  \end{itemize} 
  Summarizing,
  \[
    \mathdefin{\alpha(v,u)} :=
    \begin{cases}
      \langle 0, \bin(\idx(C')), \kappa_{x'}(A_{y'}\setminus A_{x'}, u)\rangle&\text{if $x'y'\in Y$,} \\
      \langle 1, \bin(\idx(C')), \rho_{C'}(x'), \mu_{x'}(A_{y'}\setminus A_{x'}), \kappa_{x'}(A_{y'}\setminus A_{x'}, u)\rangle&\text{if  $x'y'\notin Y$.} 
    \end{cases}
  \]
  Then, in the first case,
  \[
    |\alpha(v,u)| \le h_2(n) + \Oh(\log\log n) \enspace .
  \] 
  In the second case,
  \[
    |\alpha(v,u)| \le h_3(n) + h_2(n)+ \Oh(\log b+\log\log n) \enspace .
  \]
  In either case, $|\alpha(v,u)|\le h_3(n)+h_2(n)+ \Oh(\log b+\log\log n)$.  
  Define
  \[
    \mathdefin{\alpha(v)} := \langle \alpha(v,u):u\in N_G(v)\cap A_x \rangle \enspace .
  \] 
  Then,
  \begin{align*}
    |\alpha(v)| & \le |N_G(v)\cap A_x|\cdot (h_3(n)+h_2(n)+\Oh(\log b+\log\log n)) \\
    & \le 3k\cdot(h_3(n) + h_2(n) + \Oh(\log b+\log\log n) ) \enspace ,
    & \text{(by property~\ref{td_width} of $\calT$)}
  \end{align*}
  and, by \cref{id_x_bound,mu_x_v_upper_bound} and since $\delta(v)\ge\omega(v)$,
  \begin{align*}
    |\ell(v)| & \le |\id(x)| + |\mu_x(v)| + |\alpha(v)| + \Oh(\log\log n) \\
     & \le \log \omega(G) - \log \omega(v) + h_1(n) + \Oh((h_3(n)+\log\log n)\log_b n) + |\alpha(v)| \\
     & \le \log \omega(G) - \log \omega(v) + h_1(n) \\
     & \quad + \Oh((h_3(n)+\log\log n)\log_b n) \\
     & \quad + \Oh(k\cdot(h_3(n) + h_2(n) + \log b + \log\log n) ) \enspace ,
  \end{align*}
  which satisfies the bound on $f(n)$.

  We first verify that $\ell_{G_0}$ is injective. For distinct $v,w\in V(G_0)$, let $x$ and $y$ be their respective home nodes in $\calT$. If $x\neq y$, then $\id(x)\neq\id(y)$ because $\id:V(T)\to\{0,1\}^*$ is injective, so $\ell_{G_0}(v)=\langle\id(x),\dots\rangle \neq \langle\id(y),\dots\rangle = \ell_{G_0}(w)$. If $x=y$, then $v,w\in B_x\setminus A_x$ are distinct vertices in the torso at $x$. By \cref{injective_vertex_labels}, the vertex labelling $\mu_x$ is injective, so $\mu_x(v)\neq\mu_x(w)$, and therefore $\ell_{G_0}(v)\neq\ell_{G_0}(w)$.

  \paragraph{\boldmath The adjacency tester $A'$:} We now describe the adjacency tester $A':\{0,1\}^*\times\{0,1\}^*\to\{0,1\}$ which takes two labels $\ell(v)$ and $\ell(w)$ and returns $1$ if and only if $vw\in E(G)$.  Suppose $v\in B_x\setminus A_x$ and $w\in B_y\setminus A_y$ for some $x,y\in V(T)$.  Then $\ell(v)=\langle \id(x), \mu_x(v), \alpha(v)\rangle$ and $\ell(w)=\langle \id(y), \mu_y(w), \alpha(w)\rangle$.  The adjacency tester begins by comparing $\id(x)$ and $\id(y)$.  If $\id(x)=\id(y)$ then $x=y$, since $\id:V(T)\to\{0,1\}^*$ is injective.  In this case $v,w\in B_x\setminus A_x$,  and $vw\in E(G)$ if and only if $vw\in E(G[B_x\setminus A_x])$.  In this case, the tester returns $A'(\ell(v),\ell(w)):=A(\mu_x(v),\mu_y(w))$.  
  
  Now suppose that $\id(x)\neq\id(y)$, so $x\neq y$. Since $x\neq y$, the only way in which $v$ and $w$ can be adjacent in $G$ is if $v\in N_G(w)\cap A_y$ or $w\in N_G(v)\cap A_x$.  The former case requires that $x$ be a $T$-ancestor of $y$, and the latter case requires that $y$ be a $T$-ancestor of $x$.  Let $C$ be the component of $F$ that contains $x$ and let $D$ be the component of $F$ that contains $y$.  Removing the last part of $\id(x)$ gives $\id(C)$ and removing the last part of $\id(y)$ gives $\id(D)$, so the adjacency tester has access to $\id(C)$ and $\id(D)$.  There are several cases to consider:
  \begin{enumerate}
    \item If $\id(C)\neq\id(D)$ but $\id(C)$ is a prefix of $\id(D)$, then the only way in which $v$ and $w$ can be adjacent is if $x$ is a $T$-ancestor of $y$ and $v\in N_G(w)\cap A_y$.  To check if this is the case, the adjacency tester 
    examines each entry $\alpha(w,u)$ in $\alpha(w)$ and tests if $u=v$.  If it finds such an entry, then $vw=uw\in E(G)$, so $A'(\ell(v),\ell(w)):=1$. If no such entry is found, $A'(\ell(v),\ell(w)):=0$. All that remains is to explain how $A'$ tests if $v=u$ for some entry $\alpha(w,u)$ in $\alpha(w)$. An entry $\alpha(w,u)$ can have one of two forms, which are distinguished by their first part.
    \begin{itemize}
      \item If $\alpha(w,u)=\langle 0,\bin(\idx(C')), \kappa_{x'}(A_{y'}\setminus A_{x'}, u)\rangle$, then $u\in B_{x'}\setminus A_{x'}$ and $C'$ is the component of $F$ that contains $x'$. Furthermore $x'y'\in Y$, so $y'$ is in a different component $D'\neq C'$ of $F$ than $x'$. Then, the first $\idx(C')+1$ parts of $\id(D)$ are $\id(D')=\id(C')\circ\tau(C',D')$. The last part of $\id(D')$ is $\tau(C',D')=\langle\rho_{C'}(x'),\mu_{x'}(A_{y'}\setminus A_{x'}), \rho_{A_{y'}\setminus A_{x'}}(y')\rangle$.  Therefore, the tester can extract $\id(C')$ and $\id(x')=\id(C')\circ\rho_{C'}(x')$ using $\id(D)$ and $\idx(C')$.  After doing this, the tester checks if $\id(x')=\id(x)$.  If not, then $x'\neq x$ and $u\neq v$ since $u\in B_{x'}\setminus A_{x'}$ and $v\in B_x\setminus A_x$. If $\id(x')=\id(x)$ then $x'=x$, since $\id:V(T)\to\{0,1\}^*$ is injective.  In this case, the tester can check if $u=v$ by checking if $I(\mu_{x'}(A_{y'}\setminus A_{x'}),\kappa_{x'}(A_{y'}\setminus A_{x'}, u), \mu_{x}(v))=1$ (the first argument to $I(\cdot,\cdot,\cdot)$ comes from $\tau(C',D')$, the second from $\alpha(w,u)$ and the third from $\ell(v)=\langle\id(x),\mu_x(v),\alpha(v)\rangle$).  If so, then $v=u\in N_G(w)$ so the tester outputs $A'(\ell(v),\ell(w)):=1$.  If not, the tester continues with the next entry (if any) in $\alpha(w)$.
      
      \item If $\alpha(w,u)=\langle 1, \bin(\idx(C')), \rho_{C'}(x'), \mu_{x'}(A_{y'}\setminus A_{x'}), \kappa_{x'}(A_{y'}\setminus A_{x'}, u)\rangle$, then $u\in B_{x'}\setminus A_{x'}$ and $C'$ is the component of $F$ that contains $x'$.  In this case, the first $\idx(C')$ parts of $\id(y)$ are $\id(C')$.  The tester computes $\id(x')=\id(C')\circ\rho_{C'}(x')$.  After doing this, the tester checks if $\id(x')=\id(x)$.  If not, then $x'\neq x$ and $u\neq v$ since $u\in B_{x'}\setminus A_{x'}$ and $v\in B_x\setminus A_x$. If $\id(x')=\id(x)$ then $x'=x$ since $\id$ is injective.  In this case, the tester can check if $u=v$ by checking if $I(\mu_{x'}(A_{y'}\setminus A_{x'}),\kappa_{x'}(A_{y'}\setminus A_{x'}, u), \mu_{x}(v))=1$ (the first and second arguments to $I(\cdot,\cdot,\cdot)$ come from $\alpha(w,u)$ and the third comes from $\ell(v)=\langle\id(x),\mu_x(v),\alpha(v)\rangle$).  If so, then $v=u\in N_G(w)$ so the tester outputs $A'(\ell(v),\ell(w)):=1$.  If not, the tester continues with the next entry (if any) in $\alpha(w)$. 
    \end{itemize}
    \item If $\id(C)\neq\id(D)$ but $\id(D)$ is a prefix of $\id(C)$, then the tester proceeds exactly as in the previous case, but swapping the roles of $v$ and $w$. 
    \item If $\id(C)=\id(D)$ then $C=D$ and $x$ and $y$ are in the same component of $F$.  In this case, there are two ways in which $v$ and $w$ can be adjacent:  Either $x$ is a $T$-ancestor of $y$ and $v\in N_G(w)\cap A_y$ or $y$ is a $T$-ancestor of $x$ and $w\in N_G(v)\cap A_x$.  Case 1 above explains how to test the former and Case 2 explains how to test the latter.
    The fact that $\id(C)=\id(D)$ simplifies some parts of the test, but also introduces a new complication.  Specifically, the tester does not immediately know if $x$ is an ancestor of $y$, $y$ is an ancestor of $x$, or if neither of $x$ nor $y$ is an ancestor of the other.  Therefore, the tester must check if $u=v$ for some $u\in N_G(w)\cap A_y$ and must also check if $u=w$ for some $u\in N_G(v)\cap A_x$. We only explain how to test if $u=v$ for some $u\in N_G(w)\cap A_y$, since testing if $u=w$ for some $u\in N_G(v)\cap A_x$ can be done in a symmetric manner.

    The tester tests each entry $\alpha(w,u)$ in $\alpha(w)$.  There are two possibilities, that are distinguished by the first part of $\alpha(w,u)$.

    \begin{itemize}
      \item If $\alpha(w,u)=\langle 0,\bin(\idx(C')), \kappa_{x'}(A_{y'}\setminus A_{x'}, u)\rangle$, then $u\in B_{x'}\setminus A_{x'}$ and $C'$ is the component of $F$ that contains $x'$.  In this case, the edge $x'y'\in Y$, so $x'$ and $y'$ are in different components of $F$.  In particular, $C'\neq C$, so $x'\neq x$.  Thus $u\neq v$ since $u\in B_{x'}\setminus A_{x'}$ and $v\in B_x\setminus A_x$.

      \item If $\alpha(w,u)=\langle 1, \bin(\idx(C')), \rho_{C'}(x'), \mu_{x'}(A_{y'}\setminus A_{x'}), \kappa_{x'}(A_{y'}\setminus A_{x'}, u)\rangle$, then $u\in B_{x'}\setminus A_{x'}$ and $C'$ is the component of $F$ that contains $x'$.  In this case, the first $\idx(C')$ parts of $\id(y)$ give $\id(C')$.  Then $\id(x')=\id(C')\circ\rho_{C'}(x')$.  The tester first checks that $\id(x')=\id(x)$.  If not, then $x'\neq x$ and $u\neq v$ since $u\in B_{x'}\setminus A_{x'}$ and $v\in B_x\setminus A_x$. If $\id(x')=\id(x)$ then $x'=x$ since $\id$ is injective.  In this case, the tester can check if $u=v$ by checking if $I(\mu_{x'}(A_{y'}\setminus A_{x'}),\kappa_{x'}(A_{y'}\setminus A_{x'}, u), \mu_{x}(v))=1$.  If so, then $v=u\in N_G(w)$ so the tester outputs $A'(\ell(v),\ell(w)):=1$.  If not, the tester continues with the next entry (if any) in $\alpha(w)$.
    \end{itemize}

    \item If neither $\id(C)$ nor $\id(D)$ is a prefix of the other then no node in $C$ is a $T$-ancestor of a node in $D$ and no node in $D$ is a $T$-ancestor of a node in $C$.  Thus $x$ is not a $T$-ancestor of $y$ and $y$ is not a $T$-ancestor of $x$.  Thus $vw\not\in E(G)$ so $A'(\ell(v),\ell(w)):=0$.
  \end{enumerate}

  \paragraph{\boldmath $g'_0$-orientability:}

  All that remains is to argue that the adjacency labelling scheme we have just described is $g_0'$-orientable.  For each $G_0\in\calG'$, let $\ell_{G_0}:V(G_0)\to\{0,1\}^*$ be the adjacency labelling of $G_0$ for $A'$ defined above (by restricting $\ell$ to $V(G_0)=V(G)\setminus\{v_0\}$). Let $\calA':=(\ell_{G_0}: G_0\in\calG')$.  We must show that, for each $n\in\N^+$, the $n$-universal graph $U'_n$ defined by $(\calA',A')$ has an orientation in which each vertex has outdegree at most $g'_0(n)$.

  For each $n\in\N^+$ and each edge $\beta\beta'$ in $U'_n$ there exists
  a graph $G_0\in\calG'$ with at most $n$ vertices and an edge $vv'\in E(G_0)$ such that $\beta=\ell_{G_0}(v)$ and $\beta'=\ell_{G_0}(v')$. Let $\ell:=\ell_G$ where $G$ is the graph obtained from $G_0$ as above.  We now explain how to choose the orientation of the edge $\beta\beta'$ by examining the operation of the adjacency tester on the input pair $(\ell_{G_0}(v),\ell_{G_0}(v'))$.\footnote{Although this doesn't happen, it would not matter if the orientation of $\beta\beta'$ depended on the choice of $G_0$, since we will argue that the number of edges oriented away from a fixed $\beta$ by each branch of the adjacency testing algorithm is at most $g_0'(n)$.} Let $\calT:=(T,(B_x\mid x\in V(T)))$ be the tree-decomposition of $G$ used to define $\ell_G$ and let $((\mu_x,\kappa_x)\mid x\in V(T))$, $Y\subseteq E(T)$ and $F=T-Y$ be defined as above.  Let $x$ be the home node of $v$ in $\calT$, let $x'$ be the home node of $v'$ in $\calT$, let $C$ be the component of $F$ that contains $x$ and let $C'$ be the component of $F$ that contains $x'$.  By the definition of the adjacency tester $A'$, exactly one of the following cases applies:
  \begin{enumerate}\setcounter{enumi}{-1}
    \item $\id(x)=\id(x')$. In this case $x=x'$ since $\id:V(T)\to\{0,1\}^*$ is injective.  Thus, $v,v'\in B_x\setminus A_x$.  Then $1=A'(\ell(v),\ell(v'))=A(\mu_x(v),\mu_x(v'))$.    The triple $\tau_x=(\torso{G}{\calT,B_x}-A_x,G[B_x]-A_x,\widehat\delta_x)$ belongs to $\calT'_{\calG^{+3k},\le n}$: its graph has $m_x\le n$ vertices and its positive integer weights sum to less than $2m_x\le 2n$. Since $vv'\in E(G[B_x]-A_x)$, this triple witnesses that $\mu_x(v)\mu_x(v')$ is an edge of the universal graph $U_n$ for $(\calS,A)$.  Since the mixed labelling scheme for $(A,I)$ is $h_0$-orientable, there is an orientation of $U_{n}$ of outdegree at most $h_0(n)$.  Without loss of generality, suppose that this orientation orients the edge $\mu_x(v)\mu_x(v')$ towards $\mu_x(v')$.  In this case, we orient the edge $\beta\beta'$ of $U'_n$ towards $\beta'$.
    
    We now bound the maximum contribution to the outdegree of $\beta$ from this case. When $\id(x)=\id(x')$ and $U_n$ orients $\mu_x(v)\mu_x(v')$ towards $\mu_x(v')$, the number of possible values for $\mu_{x'}(v')$ is at most $h_0(n)$.   Since $|\alpha(v')|\in \Oh(k(h_3(n)+h_2(n)+\log b+\log\log n))$, the number of possible values for $\alpha(v')$ is at most  $\exp(\Oh(k(h_3(n)+h_2(n)+\log b+\log\log n)))$.  Thus, the number of edges of $U'_n$ incident to $\beta$ that are oriented away from $\beta$ by this branch of the adjacency testing procedure is at most $h_0(n)\cdot\exp(\Oh(k(h_3(n)+h_2(n)+\log b+\log\log n)))$.
    
    \item $\id(x)\neq\id(x')$.  In this case, $x\neq x'$.  Since $vv'\in E(G)$, one of $x$ or $x'$ is a $T$-ancestor of the other.  We orient the edge $\beta\beta'$ towards the endpoint whose home node is an ancestor of the other's home node; that is, if $x'$ is a $T$-ancestor of $x$, we orient the edge $\beta\beta'$ towards $\beta'$. Consequently, an edge is oriented away from $\beta$ towards $\beta'$ only when the home node $x'$ of $\beta'$ is an ancestor of $x$.  Since $x'$ is a $T$-ancestor of $x$, $\id(C')$ is a prefix of $\id(C)$ and is completely determined by $\id(C)$ and $\idx(C')$.  Furthermore, $\idx(C')$ is contained in $\alpha(v,v')$, which is one of the at most $3k$ entries in $\alpha(v)$.  Thus, there are at most $3k$ choices for $\id(C')$.  Then $\id(x')$ is obtained by appending $\rho_{C'}(x')$ to $\id(C')$. Since $|\rho_{C'}(x')|\le 2\log b + \Oh(\log\log n)$, there are at most $\Oh(b^2\log n)$ choices for $\rho_{C'}(x')$.  Thus, there are at most $\Oh(k b^2\log n)$ choices for $\id(x')$. 
    
    Let $y'$ be the $T$-child of $x'$ that is on the path in $T$ from $x'$ to $x$.  There are now two cases to consider.
    \begin{itemize}
      \item If $\alpha(v,v')=\langle 0,\bin(\idx(C')),\kappa_{x'}(A_{y'}\setminus A_{x'},v')\rangle$ then the edge $x'y'$ is in $Y$.  Let $D'$ be the component of $T-Y$ that contains $y'$.  Then $\id(D')$ is a prefix of $\id(C)$ and the last part of $\id(D')$ is $\tau(C',D')=\langle \rho_{C'}(x'),\mu_{x'}(A_{y'}\setminus A_{x'}),\rho_{A_{y'}\setminus A_{x'}}(y')\rangle$.  Since $vv'\in E(G)$, and $x'$ is the home node of $v'$ in $\calT$, $v'\in A_{y'}\setminus A_{x'}$. By definition,  $I(\mu_{x'}(A_{y'}\setminus A_{x'}),\kappa_{x'}(A_{y'}\setminus A_{x'},v'),\mu_{x'}(v'))=1$. 
      The triple $\tau_{x'}\in\calT'_{\calG^{+3k},\le n}$, together with the clique $A_{y'}\setminus A_{x'}$ and vertex $v'$, witnesses the occurrence of this clique label, local identifier, and vertex label. Since $(\calS,A,I)$ is strongly $h_0$-orientable, there are at most $h_0(n)$ choices for $\mu_{x'}(v')$ for each fixed clique label and local identifier.

      \item If $\alpha(v,v')=\langle 1, \bin(\idx(C')), \rho_{C'}(x'), \mu_{x'}(A_{y'}\setminus A_{x'}), \kappa_{x'}(A_{y'}\setminus A_{x'}, v')\rangle$, then $x'y'\notin Y$.  Since $v'\in A_{y'}\setminus A_{x'}$, we have $I(\mu_{x'}(A_{y'}\setminus A_{x'}),\kappa_{x'}(A_{y'}\setminus A_{x'},v'),\mu_{x'}(v'))=1$ by definition.  The triple $\tau_{x'}\in\calT'_{\calG^{+3k},\le n}$, together with the clique $A_{y'}\setminus A_{x'}$ and vertex $v'$, witnesses the occurrence of this clique label, local identifier, and vertex label. Since $(\calS,A,I)$ is strongly $h_0$-orientable, there are at most $h_0(n)$ choices for $\mu_{x'}(v')$ for each fixed clique label and local identifier.
    \end{itemize}
    In either case, there are at most $h_0(n)$ choices for $\mu_{x'}(v')$. 

    Summarizing, there are at most $\Oh(k b^2\log n)\le\exp(\Oh(\log b+\log\log n))$ choices for $\id(x')$.  There are at most $h_0(n) = (2g_0(n)+3k)2^{3k} \le g_0(n) \exp(\Oh(k))$ choices for $\mu_{x'}(v')$.  There are at most $\exp(\Oh(k(h_3(n)+h_2(n)+\log b+\log\log n)))$ choices for $\alpha(v')$.  Thus, the number of edges of $U'_n$ incident to $\beta$ that are oriented away from $\beta$ by this branch of the adjacency testing procedure is at most $h_0(n)\cdot \exp(\Oh(k(h_3(n)+h_2(n)+\log b+\log\log n)))$. 
  \end{enumerate}
  Thus, $U'_n$ has an orientation of outdegree at most $g'_0(n)\le g_0(n)\exp(\Oh(k(h_3(n)+h_2(n)+\log b+\log\log n)))$.  
\end{proof}

\section{Proof of \texorpdfstring{\cref{main_result,main_result_induced_universal,main_result_induced_universal_sparse,product_version}}{Theorems \ref*{main_result}, \ref*{main_result_induced_universal_sparse}, and \ref*{product_version}, and Corollary \ref*{main_result_induced_universal}} with Explicit Bounds}
\label{proof_main_result}

We now prove our main results, \cref{main_result,main_result_induced_universal,main_result_induced_universal_sparse,product_version}, with explicit lower order terms.  We rely on the following lemma, which we prove in \cref{weighted_product_proof}.
\begin{lem}\label{gk_labelling_precise}
  For any fixed integer $k\ge 0$, the graph class $\calR_{k}$ admits a strongly $g_0$-orientable weighted $(g_1,g_2,g_3)$ mixed labelling scheme, where
  \begin{align*}
  g_0(n)&\le \exp(\Oh(\sqrt{\log n})) \enspace ,\\ 
  g_1(n)&\in \Oh(\sqrt{\log n}) \enspace ,\\
  g_2(n)&\in \Oh(1) \enspace ,\\
  g_3(n)&\in \Oh(\sqrt{\log n}) \enspace .
\end{align*}
(More generally, for any non-decreasing function $k:\N\to\N$ with $k(n)\in o(\log n)$, the same scheme achieves $g_0(n)\le \exp(\Oh(k(n)\sqrt{\log n}))$, $g_1(n)\in \Oh(k(n)\sqrt{\log n})$, $g_3(n)\in \Oh(k(n)\sqrt{\log n})$, and $g_2(n)\le \lceil\log(2k(n)+2)\rceil\in \Oh(\log(k(n)+2))$.)
\end{lem}

\begin{lem}
  \label{main_result_detailed}
  Every proper minor-closed graph class $\calC$ admits an injective $\exp(\Oh((\log n)^{3/4}))$-orientable $(\log n + \Oh\big((\log n)^{3/4}\big))$-bit adjacency labelling scheme.
\end{lem}

\begin{proof}
  Let $k,a\geq 0$ be integers given by \cref{gmpst} for the class $\calC$.  
  Since $\calR_k$ is monotone, $\calR_k^{+a}$ is also monotone.
  By \cref{gk_labelling_precise} (if $a=0$) or by \cref{add_apexes,gk_labelling_precise} (if $a\ge 1$), $\calR_k^{+a}$ admits a strongly $\exp(\Oh(\sqrt{\log n}))$-orientable weighted $(g'_1,g'_2,g'_3)$ mixed labelling scheme, where $g'_1(n),g'_3(n)\in \Oh(\sqrt{\log n})$ and $g'_2(n)\in \Oh(1)$.

  Let $\calQ$ be the class of graphs that have a tree-decomposition whose torsos are in $\calR_k^{+a}$. The graphs in $\calR_k^{+a}$ have maximum clique size at most $2k+a+2$. Since each adhesion of a tree-decomposition is a clique in some torso, each graph in $\calQ$ has a rooted tree-decomposition of adhesion-width at most $k':=2k+a+2$ whose torsos are in $\calR_k^{+a}$.

  By \cref{lem:the_big_one} with $b(n)=\lceil 2^{(\log n)^{3/4}} \rceil$, $\calQ$ admits an injective $\exp(\Oh((\log n)^{3/4}))$-orientable $f(n)$-bit adjacency labelling scheme with
  \begin{align*}
    f(n) 
      & \le \log n + g'_1(n) + \Oh\big((g'_3(n) + k' + \log\log n)\log_{b(n)} n \\ & \quad {} + k'(g'_3(n) + g'_2(n) + k' + \log b(n) + \log\log n)\big) \\
      & = \log n + \Oh((\log n)^{3/4}) \enspace.
  \end{align*}
  By the choice of $k$ and $a$, $\calC\subseteq\calQ$.  Therefore, $\calC$ admits an injective $\exp(\Oh((\log n)^{3/4}))$-orientable $f(n)$-bit adjacency labelling scheme with $f(n)=\log n + \Oh((\log n)^{3/4})$.
\end{proof}

Together with \cref{sec:orientable}, \cref{main_result_detailed} immediately yields our main results:

\begin{proof}[Proof of \cref{main_result,main_result_induced_universal,main_result_induced_universal_sparse}]
  Let $\calC$ be a proper minor-closed graph class. By \cref{main_result_detailed}, $\calC$ admits an $f(n)$-bit adjacency labelling scheme with $f(n)=\log n + \Oh((\log n)^{3/4}) = (1+o(1))\log n$. This proves \cref{main_result}.

  Furthermore, this scheme is $g_0'(n)$-orientable with $g_0'(n) \le \exp(\Oh((\log n)^{3/4})) = n^{o(1)}$. For each $n\in\N^+$, the $n$-universal graph $U_n$ defined by this scheme has $|V(U_n)| \le 2^{f(n)+1} = n\cdot 2^{\Oh((\log n)^{3/4})} = n^{1+o(1)}$ vertices. By orientability, $U_n$ has an orientation with maximum outdegree at most $g_0'(n)$, so $|E(U_n)| \le |V(U_n)|\cdot g_0'(n) \le n^{1+o(1)}\cdot n^{o(1)} = n^{1+o(1)}$ edges. Since every $n$-vertex graph in $\calC$ is isomorphic to an induced subgraph of $U_n$, this completes the proof of \cref{main_result_induced_universal,main_result_induced_universal_sparse}.
\end{proof}

\begin{proof}[Proof of \cref{product_version}]
  Let $k\ge 1$ be an integer. By \cref{gk_labelling_precise}, the class $\calR_k$ admits a strongly $g_0$-orientable weighted $(g_1,g_2,g_3)$ mixed labelling scheme with $g_1(n)\in \Oh(k\sqrt{\log n})$ and $g_0(n)\le \exp(\Oh(k\sqrt{\log n}))$. By \cref{adj_labelling_from_mixed_labelling}, $\calR_k$ admits an injective $g_0$-orientable $f(n)$-bit adjacency labelling scheme with $f(n) = \log n + g_1(n) = \log n + \Oh(k\sqrt{\log n})$. For each $n\in\N^+$, the $n$-universal graph $U^k_n$ defined by this scheme has $|V(U^k_n)| \le 2^{f(n)+1} = n\cdot \exp(\Oh(k\sqrt{\log n}))$ vertices and, by orientability, at most $|V(U^k_n)|\cdot g_0(n) \le n\cdot\exp(\Oh(k\sqrt{\log n}))$ edges. By injectivity of the vertex labels, every $n$-vertex graph in $\calR_k$ is isomorphic to an induced subgraph of $U^k_n$.
\end{proof}

\section{Discussion}\label{conclusion}

We conclude with a brief discussion about constants, computational complexity, and extensions beyond minor-closed classes.  The constant $a$ in \cref{gmpst} is surprisingly small. Suppose that some graph $X$ is not in $\calG$, where $X-A$ is planar for some non-empty set $A\subseteq V(X)$. Then Theorem~42 in \citep{dujmovic.joret.ea:planar} (which employs the structure theorem of \citet{DvoTho}) shows that \cref{gmpst} holds with $a=|A|-1$. For example, if $K_t$ is not in $\calG$, then \cref{gmpst} holds with $a=\max\{t-5,0\}$. Moreover, using the recent polynomial bounds in the Graph Minor Structure Theorem by \citet{gorsky2025polynomial}, in \cref{gmpst} one can obtain polynomial bounds on $k$ and $a$ (as a function of the excluded minor). However, it is open whether one can simultaneously get the above optimal bound on $a$ and a polynomial bound on $k$. \citet[Conjecture~18.9]{gorsky2025polynomial} conjecture this is possible (also see the discussion in \citep{GSW26}).

We now explain why the proof of \cref{main_result} is constructive and yields a polynomial-time labelling algorithm.  The only difficult part of this process is finding the decomposition described by the Graph Minor Product Structure Theorem (\cref{gmpst}).  \citet{dujmovic.joret.ea:planar} show how to obtain this decomposition in polynomial time, given the structural decomposition in the original Robertson--Seymour Graph Minor Structure Theorem. Finding the Robertson--Seymour decomposition is complicated, but polynomial-time algorithms exist \cite{kawarabayashi2018flat,kawarabayashi2020new,gorsky2025polynomial,demaine2005algorithmic}.  Together, these give the parameters $k$ and $a$, the tree-decomposition $\calT_0$ of $G_0$, and the embedding of each torso of $\calT_0$ into a graph in $\calR_k^{+a}$. The adjacency labelling scheme in \cite{dujmovic.esperet.ea:adjacency} and its modification to weighted mixed labellings in \cref{weighted_product_proof,add_apexes} are easily implemented in polynomial time directly from their definitions.  The transformation of $\calT_0$ into the tidy rooted tree-decomposition $\calT=(T,(B_x\mid x\in V(T)))$ and edge set $Y\subseteq E(T)$ in \cref{all_singing_all_dancing} (via \cref{skinny_edges,bag_separator,log_b_td}) also runs in polynomial time.  Finally, the construction in \cref{lem:the_big_one} requires only computing the weight function $\delta$, alphabetic prefix-free codes (\cref{shannon_fano_elias_tree,shannon_fano_elias_code}), and weighted mixed labellings on the modified torsos $\torso{G}{\calT,B_x}-A_x\in\calR_k^{+a+3k'}$ (where $k':=2k+a+2$), and then assembling each vertex label $\ell(v)=\langle\id(x),\mu_x(v),\alpha(v)\rangle$ from the root-to-$x$ path in $T$.

Finally, we note that \cref{lem:the_big_one} applies beyond proper minor-closed classes. Recall that, for fixed $k$, $\calR_k$ consists of the graphs isomorphic to subgraphs of $H\boxtimes P$ for some graph $H$ of treewidth at most $k$ and some path $P$. This class is monotone and has an efficient strongly orientable weighted mixed labelling scheme (\cref{gk_labelling_precise}), so \cref{lem:the_big_one} applies to it. Moreover, for suitable $k=k(p)$, $\calR_k$ contains every $p$-planar graph (graphs that can be drawn in the plane with at most $p$ crossings per edge)~\citep{dujmovic.morin.ea:graph}, and more generally every $p$-matching-planar graph~\citep{DBLP:journals/corr/abs-2507-22395} (graphs that have a drawing in the plane such that, for every edge $e$, every matching among the edges crossing $e$ has size at most $p$). Note that these classes are not minor-closed. Consequently, for all fixed $p$ and $q$, the class of graphs having a tree-decomposition of adhesion-width at most $q$ in which every torso is $p$-matching-planar (again a class that is not minor-closed) has a $(1+o(1))\log n$-bit adjacency labelling scheme and $n$-induced-universal graphs with $n^{1+o(1)}$ vertices and edges.

\section*{Acknowledgement} 
This work was initiated during the Ninth Annual Workshop on Geometry and Graphs (v2), held January 21--28, 2022, at the Bellairs Research Institute of McGill University.  The authors are thankful to the other organizers and participants for providing a stimulating research environment.

\section*{Research History and Statement on AI Use}

The results in this paper are the product of years of work by the authors, all of whom have been thinking about adjacency-labelling and sparse universal graphs since at least February 2020, when the work on planar graphs resulting in publication \cite{dujmovic.esperet.ea:adjacency} began at the Eighth Workshop on Geometry and Graphs.  The work on adjacency labelling for proper minor-closed graph families began at the 2022 workshop mentioned above and several of the ideas in the paper (using weights to deal with an unbalanced tree-decomposition and using clique labels as a means of navigating a tree-decomposition) come from that week.  The final ideas (decomposition into skinny subtrees and a mixed labelling scheme for skinny subtrees) needed to complete the adjacency labelling scheme for proper minor-closed graph classes that resulted in publication \cite{dujmovic.gavoille.ea:adjacency_minor} were obtained during a visit to Ottawa in February 2026.  The strongly $g_0$-orientable version of this scheme that yields sparse universal graphs was developed during a visit to Kraków in June 2026.  The revision history is visible at \url{https://github.com/patmorin/mflabelling}.

Piotr~Micek gave a talk about adjacency labelling for proper minor-closed graph classes at the MATRIX Workshop on Global Structure and Geometry of Graphs in Creswick, Australia in April 2026. Pat~Morin gave a similar talk at the Jagiellonian TCS Seminar in June 2026.  Vida~Dujmović gave a talk about sparse induced-universal graphs for proper minor-closed graph classes at the 11th Polish Combinatorial Conference in September 2026.

Every line of \LaTeX\ in this paper was written by the authors. One of the authors, Pat~Morin, used an AI-enabled text editor (Antigravity IDE with Gemini) to proofread the manuscript and suggest corrections for minor errors.  The final version was checked and approved by all authors.

\setlength{\bibsep}{0.3ex plus 0.1ex minus 0.1ex}
\bibliographystyle{plainurlnat}
\bibliography{mflabelling}

\newpage
\appendix
\crefalias{section}{appendix}
\section{A Mixed Labelling Scheme for \texorpdfstring{$\calR_k$}{Rk}}
\label{weighted_product_proof}

In this appendix, we modify the labelling scheme of \citet{dujmovic.esperet.ea:adjacency} and its improvements by \citet{gawrychowski.janczewski:simpler} in order to show the existence of a strongly $g_0$-orientable efficient weighted mixed labelling scheme for $\calR_k$ mentioned in \cref{weighted_mixed_section}.  In addition to mixing results from \cite{dujmovic.esperet.ea:adjacency,gawrychowski.janczewski:simpler}, we obtain a further simplification using logarithmic height tree-decompositions (see \cref{low_height_tree_decomp}, below).  This not only simplifies the labelling schemes from \cite{dujmovic.esperet.ea:adjacency,gawrychowski.janczewski:simpler}, but also produces a sparse induced-universal graph without requiring any of the additional machinery used by \citet{EJM23}.  

A rooted tree-decomposition $(T,(B_x\mid x\in V(T)))$ is a \defin{binary tree-decomposition} if $T$ is a binary tree. Variants of the following result appear in several contexts, with different bounds for the width of the tree-decomposition and the height of the tree $T$. For example, \citet{bodlaender.hagerup:parallel} use a similar result in the context of parallel algorithms and \citet{jain.tewari:space} use a similar result in the context of space-efficient algorithms. The result cited here is due to \citet{chatterjee.ibsen-jensen.ea:optimal}.

\begin{lem}[\protect{\cite[Lemma~7]{chatterjee.ibsen-jensen.ea:optimal}}]\label{low_height_tree_decomp}
  Let $H$ be an $n$-vertex graph of treewidth at most $k$. Then $H$ has a rooted binary tree-decomposition $\calT:=(T,(B_x\mid x\in V(T)))$ of width at most $6k+5$ in which $\height(T)\in \Oh(\log n)$.
\end{lem}

\begin{proof}[Proof of \cref{gk_labelling_precise}]
  \cref{gk_labelling_precise} is establishing several results at once, using combinations of techniques from \cite{dujmovic.esperet.ea:adjacency,gawrychowski.janczewski:simpler,EJM23} and \cref{low_height_tree_decomp}.  Thus, we present the proof in stages.
  We begin by describing a $(\log n+O(\sqrt{\log n}))$-bit adjacency labelling scheme for $\calR_k$ that uses an adjacency tester $A$.  Next we describe clique labels, local identifiers, and an identity tester $I$ that, along with the adjacency labelling for $A$, give an unweighted strongly $g_0$-orientable $(g_1,g_2,g_3)$ mixed labelling scheme.  Finally, we show how to make this into a weighted strongly $g_0$-orientable $(g_1,g_2,g_3)$ mixed labelling scheme using the trick of \emph{simulating weight by multiplicity}.

  Since $\calR_0\subseteq\calR_1$, we may assume without loss of generality that $k\ge 1$; we may also assume that $n\ge 2$, since for $n=1$ the bounds hold trivially.
  As usual, we first describe the format of the labels before describing the operation of the adjacency tester $A$. Let $\tau=(G^+,G,\omega)\in\calT_{\calR_k}$ and let $n=|V(G)|$. (For now we assume $\omega(v)=1$ for all $v\in V(G)$.) Thus, $G^+$ and $G$ are subgraphs of $H\boxtimes P$, where $H$ is a graph of treewidth at most $k$ and $P$ is a path. Without loss of generality, assume that $P=(1,2,\dots,h)$, with $h\le n$.  Each vertex of $G$ is a pair $(v,i)$ with an \defin{$H$-coordinate} $v\in V(H)$ and a \defin{$P$-coordinate} $i\in V(P)$. We may assume, without loss of generality, that $G$ contains at least one vertex with $P$-coordinate $i$ for each $i\in[h]$ (and therefore $h\le n$) and that each vertex of $H$ appears as the $H$-coordinate of at least one vertex of $G$ (and therefore $|V(H)|\le n$). The label of each vertex $(v,i)$ of $G$ has two parts. The first part, $\lambda_1(v,i)$, is determined by the $P$-coordinate, $i$. The second part, $\lambda_2(v,i)$, is obtained from a labelling of $H[V_i]$ (where $V_i$ is defined below), so it represents the $H$-coordinate, $v$.

  By \cref{low_height_tree_decomp}, $H$ has a rooted binary tree-decomposition $\calT:=(T,(B_x\mid x\in V(T)))$ of width at most $6k+5$ in which $\height(T)\in \Oh(\log n)$.   For each $x\in V(T)$, add edges to $H$ so that $B_x$ becomes a clique. (Note that for any clique $K$ in $G^+$, its projection onto $V(H)$ is a clique in the original treewidth-$k$ graph $H$ before adding these edges, and therefore has size at most $k+1$.)  Define a strict partial order $<$ on $V(H)$ such that $v<w$ if the home of $v$ in $\calT$ is a strict $T$-ancestor of the home of $w$ in $\calT$ (resolving ties arbitrarily for vertices with the same home node). For each $w\in V(H)$, let $N^+_H(w)=\{v\in N_H(w)\mid v < w\}$.  Then $|N^+_H(w)|\le 6k+5$ since $v\in N^+_H(w)$ implies that $v\in B_x\setminus\{w\}$ where $x$ is the home of $w$ in $\calT$.  Note that, for each clique $K$ in $H$, $K\subseteq N^+_H(w)\cup\{w\}$ for some $w\in K$. 

  \paragraph{\boldmath Mapping $V(H)$ onto intervals:}
  The first step in \cite{dujmovic.esperet.ea:adjacency} is an idea used (at least implicitly) by \citet{gavoille.labourel:shorter} in their labelling scheme for bounded-treewidth graphs.  Here we refine this idea using \cref{low_height_tree_decomp}.
  We map each vertex $v$ of $H$ to a closed interval \defin{$f(v)$} with integer endpoints in $\{1,\ldots,2|V(H)|\}$ and such that each integer is an endpoint of at most $1$ interval.  This mapping has the following properties (for some fixed $c>0$):
  \begin{enumerate}[nosep,label=(I\arabic*)]
    \item\label{nested_or_disjoint} For any two vertices $v,w\in V(H)$, $f(v)\subset f(w)$, $f(w)\subset f(v)$ or $f(v)\cap f(w)=\emptyset$.
    \item\label{intersection_subgraph} For each $v\in V(H)$ and each $w\in N^+_H(v)$, $f(v)\subset f(w)$.
    \item\label{thin} For each $p\in \R$, $|\{v\in V(H):p\in f(v)\}|\le ck\log n$. 
  \end{enumerate}

  The existence of such a mapping $f$ is easily established by a recursive procedure that takes a node $x$ of $T$ and an interval $[a,b]$ where $a$ and $b$ are integers with $b-a+1=2|\bigcup_{y\in V(T_x)} (B_y\setminus A_y)|$, where $A_y$ denotes the parent adhesion of $y$ in $\calT$. Let $v_1,\ldots,v_r$ be the vertices in $B_x\setminus A_x$ ordered so that $v_1< \cdots< v_r$.  Then $f(v_i)=[a+i-1,b-i+1]$ for each $i\in\{1,\ldots,r\}$.  If $x$ has a left child $x_1$ then let $n_1:=|\bigcup_{y\in V(T_{x_1})} (B_y\setminus A_y)|$. If $x$ has no left child, then $n_1:=0$.  The procedure is then applied recursively on the left child $x_1$ of $x$ (if any) with the interval $[a+r,m]$ and the procedure is then applied recursively to the right child $x_2$ of $x$ (if any) with the interval $[m+1,b-r]$, where $m:=a+r+2n_1-1$.  Initially, the procedure is called with $[a,b]=[1,2|V(H)|]$ and $x$ as the root of $T$.

  Let $G_f$ be the interval intersection graph with vertex-set $V(G_f):=V(H)$ and edge-set $E(G_f):=\{vw\in\binom{V(H)}{2}:f(v)\cap f(w)\neq\emptyset\}$.
  Property~\ref{intersection_subgraph} implies that $H$ is a subgraph of $G_f$. Property~\ref{thin} implies that $G_f$ has no clique larger than $ck\log n$, for some fixed constant $c$.
  The graph $G_f$ has a proper colouring $\mathdefin{\varphi}:V(H)\to[\lfloor ck\log n\rfloor]$ \cite{gilmore.hoffman:characterization}.
  The label for each vertex $(v,i)$ of $G$ will include the $\Oh(\log (k\log n))$-bit integer $\varphi(v)$.

  Since every clique $K$ in $H$ satisfies $K\subseteq N^+_H(w)\cup\{w\}$ for some $w\in K$, Property~\ref{intersection_subgraph} implies:
  \begin{enumerate}[nosep,label=(I\arabic*)]\setcounter{enumi}{3}
    \item\label{clique_property} Let $K$ be a clique in $H$.  Then there exists $w\in K$ such that $K\setminus\{w\}\subseteq N^+_H(w)$ and $f(w)\subset f(v)$ for each $v\in K\setminus\{w\}$.  
  \end{enumerate}
  Property \ref{clique_property} is a crucial difference between the scheme we describe here and the labelling schemes of \citet{dujmovic.esperet.ea:adjacency} and that of \citet{gavoille.labourel:shorter} on which it is based.  In essence, it ensures that the ancestor-descendant relationships in $T$ are consistent with containment relationships in the interval representation.  This simultaneously simplifies the labelling scheme and makes it $g_0$-orientable without requiring additional tricks used by \citet{EJM23}.\footnote{Specifically, Property~\ref{clique_property} completely eliminates the need for the extended signature $\sigma^+(v,i)$ of a vertex $(v,i)$ of $G$.  This allows for a number of simplifications.  Furthermore, the difference in length between $\sigma^+(v,i)$ and the signature $\sigma(v,i)$ is the main reason that the adjacency labelling scheme in \cite{dujmovic.esperet.ea:adjacency} is not $n^{o(1)}$-orientable, so this modification almost immediately yields a universal graph with $n^{1+o(1)}$ edges.}

  \paragraph{\boldmath Treating $G$ as a sequence of graphs and mapping $V(H)$ onto points:}
  The next idea in \cite{dujmovic.esperet.ea:adjacency} is to treat $V(G)$ as a sequence of vertex subsets of $H$.
  For each $i\in[h]$, let $V_i:=\{v\in V(H): \{(v,i-1),(v,i),(v,i+1)\}\cap V(G)\neq\emptyset \}$.    Each vertex $(v,i)$ of $G$ contributes a vertex to at most three sets $V_{i-1}$, $V_{i}$ and $V_{i+1}$, so $\sum_{i=1}^h |V_i|\le 3n$.  For each vertex $(v,i)$ of $G$, define the set $\mathdefin{X_{v,i}}:= \{v\}\cup (V_i\cap N^+_H(v))$.  In our labelling scheme, the label of $(v,i)$ will include enough information to identify the vertices in $N_G((v,i))\cap (X_{v,i}\times\{i-1,i,i+1\})$.  The following property shows that this is sufficient to cover every edge of $G$:
  \begin{enumerate}[label=(Q\arabic*),nosep,nolistsep]
    \item \label{obs:x_v_i_sufficient}
      For each edge $(v,i)(w,j)$ of $G$ with $w\in N^+_H(v)$, we have that $w\in X_{v,i}$.
  \end{enumerate}

  For each vertex $v$ of $H$, let \defin{$x_f(v)$} be the right endpoint of the interval $f(v)$.
  For each $i\in[h]$, let $\mathdefin{S_i}:=\{x_f(v):v\in V_i\}$.
  Note that $S_i$ is a set of integers, for each $i\in[h]$.
  It is helpful for consecutive sets $S_i$ and $S_{i+1}$ to satisfy certain properties, for each $i\in[h-1]$. Therefore, \citet{dujmovic.esperet.ea:adjacency} and \citet{gawrychowski.janczewski:simpler} define a sequence of supersets \defin{$S^+_1,\ldots,S^+_h$}. The exact definition of $S^+_1,\ldots,S^+_h$ is not necessary for the current discussion.
  For our purposes it is sufficient to know that for any fixed integer $q\ge 1$ the sets $S^+_1,\ldots,S^+_h$ can be chosen so that  
  \begin{enumerate}[nosep,label=(\alph*)]
    \item\label{plus_superset} $S^+_i\supseteq S_i$ for each $i\in[h]$;
    \item\label{plus_size} there is a constant $c>0$ such that $\sum_{i=1}^h |S^+_i|\le cq^2n$;
    \item\label{plus_chunking_down} for each $i\in\{1,\ldots,h-1\}$, for each $t\in[q]$, and for each set of $t$ consecutive elements of  $S^+_i$, at least $t-1$ of these elements are in $S^+_{i+1}$; and 
    \item\label{plus_chunking_up} for each $i\in\{2,\ldots,h\}$, for each $t\in[q]$, and for each set of $t$ consecutive elements of $S^+_i$, at least $t-1$ of these elements are in $S^+_{i-1}$.     
  \end{enumerate}
  We remark that the contents of the sets $S^+_1,\ldots,S^+_h$ are entirely determined by the contents of the sets $S_1,\ldots,S_h$.
  (Formally, the sequence $S^+_1,\ldots,S^+_h$ is the result of applying a function $Q:(2^\R)^*\to (2^\R)^*$ whose input is a sequence of finite sets of real numbers and whose output is a sequence of finite sets of real numbers. The function $Q$ is universal; it does not depend on $H$, $P$, or $G$.)

  \paragraph{\boldmath The row label $\lambda_1$:}
  With the sets $S^+_1,\ldots,S^+_h$ described, we can now discuss $\lambda_1$.
  The first part of $\lambda_1(v,i)$ is obtained by applying \cref{shannon_fano_elias_code} to the set $[h]$ using the weight function defined by $i \mapsto |S^+_i|$ to obtain a prefix-free code $\rho:[h]\to\{0,1\}^*$ and functions $\overleftarrow{\pi}:[h]\to\{0,1\}^*$ and $\overrightarrow{\pi}:[h]\to\{0,1\}^*$. By \cref{shannon_fano_elias_code} and Property~\ref{plus_size} of $S^+_1,\ldots,S^+_h$, $|\rho(i)|\le \log(\sum_{i=1}^h |S^+_i|) - \log|S^+_i| + \Oh(1)\le \log n - \log|S^+_i| + \Oh(1)$ for each $i\in[h]$.

  Then
  \[
    \mathdefin{\lambda_1(v,i)}:=\langle \rho(i),\overleftarrow{\pi}(i),\overrightarrow{\pi}(i)\rangle
  \]
  and
  \begin{equation}
    |\lambda_1(v,i)|\le \log n - \log |S^+_i| + \Oh(\log\log n) \enspace ,
    \label{lambda_one_length}
  \end{equation}
  for each vertex $(v,i)$ of $G$.

  \paragraph{B-trees:}

  Until this point, the labelling scheme we have been describing comes from \cite{dujmovic.esperet.ea:adjacency}, except for the modifications to the interval mapping function $f$ using \cref{low_height_tree_decomp} that guarantee \ref{clique_property}.  At this point, we continue with the B-tree based labelling scheme from \cite{gawrychowski.janczewski:simpler} rather than the binary-tree based scheme in \cite{dujmovic.esperet.ea:adjacency}.

  A \defin{B-tree} stores a set of real numbers in the leaves of a tree $B$ in which all leaves have the same depth.
  A \defin{weight-balanced B-tree} is parameterized by a parameter $a\ge 6$.  A weight-balanced B-tree $B$ is \defin{semi-balanced} if, for each node $x$ of $B$, the number of leaves in the subtree $B_x$ is at most $6a^{\height(B_x)}$ and for each \emph{non-root} node $x$ of $B$, the number of leaves in $B_x$ is at least $\tfrac{1}{2}a^{\height(B_x)}$.
  Let $B$ be a weight-balanced B-tree with $r\in O(n)$ leaves.  The weight-balance conditions guarantee that $\height(B)\le \log_a r + 1$.
  The weight-balance conditions also ensure that the number of children of any node $x$ of $B$ is at most $12a$.
  This implies that any root-to-leaf path in a B-tree can be encoded by a sequence of $\log_a r+1$ integers, each requiring $\lceil\log(12a)\rceil\le 5+\log a$ bits to represent, for a total of at most $\lceil\log(12a)\rceil\cdot(\log_a r+1)=\log r + \Oh(\log_a r + \log a)$ bits.
  The value of $a$ used in \cite{gawrychowski.janczewski:simpler} and here is approximately $2^{\sqrt{\log n}}$, so that $\log a\in O(\sqrt{\log n})$ and $\log_a n\in O(\sqrt{\log n})$.

  We use a sequence of weight-balanced B-trees, $B_1,\ldots,B_h$, where the leaves of $B_i$ store the values in $S^+_i$.  The implementation details of these B-trees are not important here. What is important is that the transformation from $B_i$ to $B_{i+1}$ is performed by a rebalancing operation that is applied only to nodes of depth $d_i$.  Rebalancing a node $x$ of depth $d_i$ involves replacing the children $y_1,\ldots,y_p$ of $x$ with children $y_1',\ldots,y_m'$ and redistributing the set of children of $y_1,\ldots,y_p$ so that they become children of $y_1',\ldots,y'_m$.  Therefore, the set of grandchildren of $x$ in $B_i$ is the same as the set of grandchildren of $x$ in $B_{i+1}$ (but the path in $B_i$ from $x$ to a grandchild $z$ is different from the path in $B_{i+1}$ from $x$ to $z$).  This has the following important consequence:
  \begin{enumerate}[nosep,label=(Q\arabic*)]\setcounter{enumi}{1}
    \item\label{path_splice} For each $i\in[h-1]$ and each $\ell\in S^+_{i}\cap S^+_{i+1}$, the root-to-$\ell$ path $P_{B_{i}}(\ell)$ in $B_{i}$ and the root-to-$\ell$ path $P_{B_{i+1}}(\ell)$ in $B_{i+1}$ differ by at most a single vertex.  More precisely, a subpath $xyz$ in $P_{B_{i}}(\ell)$ can be replaced by a subpath $xy'z$ to give $P_{B_{i+1}}(\ell)$.%
    \footnote{There are two exceptions. Occasionally, $P_{B_{i+1}}(\ell)$ is obtained by removing the first vertex of $P_{B_{i}}(\ell)$ or by adding a vertex at the beginning of $P_{B_{i}}(\ell)$. Since these exceptions are distracting and easily handled, we ignore them here.}
  \end{enumerate}

  \paragraph{\boldmath The vertex signature $\sigma(v,i)$:}

  For each $i\in[h]$ and each $x\in V(B_i)$, define the \defin{signature} \defin{$\sigma_{B_i}(x)$} as the encoding of $P_{B_i}(x)$ as a sequence of $\depth_{B_i}(x)$ integers in $\{0,\ldots,12a-1\}$.  Then 
  \[
    |\sigma_{B_i}(x)| \le \lceil\log (12a)\rceil\cdot\depth_{B_i}(x) 
    \le (\log a)\cdot\depth_{B_i}(x) + O(\sqrt{\log n})   \enspace .
  \]

  For each $i\in [h]$ and each $S\subseteq S^+_i$, define \defin{$\lca'_{B_i}(S)$} as the maximum depth non-leaf node in $B_i$ that contains all elements of $S$ in its subtree. (Equivalently, $\lca'_{B_i}(S)=\lca_{B_i}(S)$ if $\lca_{B_i}(S)$ is not a leaf of $B_i$, otherwise $\lca'_{B_i}(S)=\parent_{B_i}(\lca_{B_i}(S))$.) For each $i\in[h]$ and each $v\in V_i$, define $\mathdefin{x'_{B_i}(v)}:=\lca'_{B_i}(f(v)\cap S^+_i)$ and define \defin{$x_{B_i}(v)$} as the child of $x'_{B_i}(v)$ that has $x_f(v)$ as a leaf-descendant.   
  (The node $x'_{B_i}(v)$ exists because $x_f(v)\in S_i\subseteq S^+_i$ and $x_f(v)\in f(v)$ and $x_{B_i}(v)$ exists because $x'_{B_i}(v)$ is a non-leaf node.)  We make use of the following property of this definition:
  \begin{enumerate}[nosep,label=(Q\arabic*)]\setcounter{enumi}{2}
    \item\label{smallest_in} For each $i\in[h]$ and each $v\in V_i$, the smallest leaf-descendant of $x_{B_i}(v)$ is contained in $f(v)$.
  \end{enumerate}
  For each $i\in[h]$ and each $v\in V_i$, define the \defin{signature} $\mathdefin{\sigma(v,i)}:=\sigma_{B_i}(x_{B_i}(v))$.
  Then
  \begin{align*}
    |\sigma(v,i)| 
    & = |\sigma_{B_i}(x_{B_i}(v))| \\
    & \le  (\log a)\cdot \depth_{B_i}(\lca_{B_i}(f(v)\cap S^+_i)) + \Oh(\log_a n) \\
    & \le (\log a)\cdot\height(B_i) + \Oh(\log_a n) \\
    & \le \log |S^+_i| + \Oh(\log_a n + \log a) \\
    & = \log |S^+_i|+\Oh(\sqrt{\log n}) \enspace .
  \end{align*}
  The following is the critical property of the vertex signatures:
  \begin{enumerate}[nosep,label=(Q\arabic*)]\setcounter{enumi}{3}
   \item For any two vertices $(v,i)$ and $(w,j)$ of $G$, $(\rho(i),\sigma(v,i),\varphi(v))=(\rho(j),\sigma(w,j),\varphi(w))$ if and only if $(v,i)=(w,j)$.
  \end{enumerate}  
  Indeed, since $\rho$ is an injective function $\rho(i)=\rho(j)$ if and only if $i=j$.  If $\sigma(v,i)=\sigma(w,i)$ then $y:=x_{B_i}(v)=x_{B_i}(w)$. By \ref{smallest_in}, $f(v)$ and $f(w)$ each contain the smallest leaf-descendant of $y$, so $f(v)\cap f(w)\neq\emptyset$.   Since $\varphi$ is a proper colouring of the intersection graph $G_f$, $\varphi(v)=\varphi(w)$ and $f(v)\cap f(w)\neq\emptyset$ if and only if $v=w$.

  We make use of the fact that $\sigma(v,i)$ is defined for many pairs $(v,i)$ that are not vertices of $G$ (but not necessarily all vertices of $H\boxtimes P$).  We make use of the fact that, for each $i\in[h]$ and each $v\in V_i$, $\sigma(v,i)$ is defined, which implies that $\sigma(v,i)$ is defined if any of $(v,i-1)$, $(v,i)$, or $(v,i+1)$ is a vertex of $G$.

  Let $(v,i)$ and $(w,j)$ be two vertices of $G$ such that $w\in \{v\}\cup N^+_H(v)$ and $|j-i|\le 1$.  Since $|j-i|\le 1$, $x'_{B_j}(v)$, $x_{B_j}(v)$, and $\sigma(v,j)$ are all defined.  By Property~\ref{intersection_subgraph},  $f(v)\subseteq f(w)$.  Therefore  $x'_{B_j}(w)$ is a $B_j$-ancestor of $x'_{B_j}(v)$ and of its child $x_{B_j}(v)$.
  Therefore, $\sigma_{B_j}(x'_{B_j}(w))$ is a prefix of $\sigma(v,j)=\sigma_{B_j}(x_{B_j}(v))$.  Since $x_{B_j}(w)$ is a child of $x_{B_j}'(w)$, this implies that $\sigma(w,j)$ consists of a prefix of $\sigma(v,j)$ followed by the $\lceil\log(12a)\rceil$-bit representation of $c_{w,j}$, where $c_{w,j}\in\{0,\ldots,12a-1\}$ is the index of $x_{B_j}(w)$ in the list of children of $x'_{B_j}(w)$.  To summarize, it is possible to construct $\sigma(w,j)$ given $\sigma(v,j)$, $|\sigma(w,j)|$, and $c_{w,j}$.

  \paragraph{\boldmath The adjacency list $\alpha(v,i)$:}
  For each $(v,i)\in V(G)$, define
  \begin{align*}
    \mathdefin{\alpha^{-}(v,i)}
     & := \langle\; \langle |\sigma(w,i-1)|, c_{w,i-1}, \varphi(w)\rangle \mid w\in \{v\}\cup N^+_H(v),\; (w,i-1)\in N_G((v,i)) \;\rangle   \\
    \mathdefin{\alpha^{0}(v,i)}
     & := \langle\; \langle |\sigma(w,i)|, c_{w,i}, \varphi(w)\rangle \mid w\in N^+_H(v),\; (w,i)\in N_G((v,i)) \;\rangle               \\
    \mathdefin{\alpha^+(v,i)}
     & := \langle\; \langle |\sigma(w,i+1)|, c_{w,i+1}, \varphi(w)\rangle \mid w\in \{v\}\cup N^+_H(v),\; (w,i+1)\in N_G((v,i)) \;\rangle \\
    \mathdefin{\alpha(v,i)}
     & := \langle\; \alpha^-(v,i),\alpha^0(v,i),\alpha^+(v,i) \;\rangle        \enspace .
  \end{align*}
  Since $|N^+_H(v)|\le 6k+5$ for each $v\in V(H)$, $|\alpha(v,i)|\in\Oh(k\sqrt{\log n})$.

  \paragraph{\boldmath The transition codes $\overleftarrow{\tau}(v,i)$ and $\overrightarrow{\tau}(v,i)$:}
  The last (and critical) piece of the labelling schemes in \cite{dujmovic.esperet.ea:adjacency,gawrychowski.janczewski:simpler} is the notion of a \defin{transition code} that makes it possible to derive the signature $\sigma(v,i+1)$ (or $\sigma(v,i-1)$) from $\sigma(v,i)$ and a bitstring of length $\Oh(\sqrt{\log n})$. Formally, there is a universal function $J:(\{0,1\}^*)^2\to\{0,1\}^*$ such that, for each $i\in[h-1]$ and each $v\in V_i\cap V_{i+1}$ (which includes every $(v,i)\in V(G)$ with $i< h$), there exists a bitstring \defin{$\overrightarrow{\tau}(v,i)$} of length $\Oh(\sqrt{\log n})$ such that $J(\sigma(v,i),\overrightarrow{\tau}(v,i))=\sigma(v,i+1)$ (for $i=h$, let $\overrightarrow{\tau}(v,h):=\varepsilon$). Symmetrically, for each $i\in\{2,\ldots,h\}$ and each $v\in V_i\cap V_{i-1}$ (which includes every $(v,i)\in V(G)$ with $i\ge 2$), there exists a bitstring \defin{$\overleftarrow{\tau}(v,i)$} of length $\Oh(\sqrt{\log n})$ such that $J(\sigma(v,i),\overleftarrow{\tau}(v,i))=\sigma(v,i-1)$ (for $i=1$, let $\overleftarrow{\tau}(v,1):=\varepsilon$).

  We now describe the format of the bitstring $\overrightarrow{\tau}(v,i)$ and how it is used by the transition function $J$.  The bitstring $\overleftarrow{\tau}(v,i)$ is defined in a completely symmetric manner.  The transition from $\sigma(v,i)$ to $\sigma(v,i+1)$ proceeds in three steps.

  In the first step, we take $\sigma(v,i)$ onto $\sigma_{B_i}(\ell)$ for some leaf $\ell\in f(v)\cap S^+_{i}\cap S^+_{i+1}$. At least one such leaf $\ell$ exists because $v\in V_i\cap V_{i+1}$, and therefore $x_f(v)\in S_i\cap S_{i+1}\subseteq S^+_i\cap S^+_{i+1}$ is a leaf in $B_i$ and $B_{i+1}$. Let $y=x_{B_i}(v)$.  If $y$ is a leaf of $B_i$ then $y= x_f(v)$ because $x_f(v)\in S_i\subseteq S^+_i$ and $x_f(v)\in f(v)$.  In this case we can simply take $\ell:=y$ and we already have the signature $\sigma_{B_i}(\ell)=\sigma(v,i)$.  
  
  Now suppose that $y:=x_{B_i}(v)$ is not a leaf of $B_i$. By \ref{smallest_in}, the smallest leaf-descendant $\ell'$ of $y$ is in $f(v)$.   The signature $\sigma_{B_i}(\ell')$ can be obtained from $\sigma_{B_i}(y)=\sigma(v,i)$ by appending a sequence of $\lceil\log(12a)\rceil(\height(B_i)-\depth_{B_i}(y))$ zeros. If $\ell'\in S^+_{i+1}$ then we take $\ell:=\ell'$ and we are done.  If $\ell'\not\in S^+_{i+1}$ then $\ell'<x_f(v)$ and, by property~\ref{plus_chunking_down} of $S^+_1,\ldots,S^+_h$, the successor $\ell''$ of $\ell'$ in $S^+_i$ is in $S^+_{i+1}$.  In this case, we take $\ell:=\ell''$ and the signature $\sigma_{B_i}(\ell'')$ can be obtained from $\sigma_{B_i}(y)$ by appending a sequence of $\lceil\log(12a)\rceil(\height(B_i)-\depth_{B_i}(y))-1$ zeros followed by a $1$.  In each case, the information needed to obtain $\sigma_{B_i}(\ell)$ can be recorded in $\Oh(\log{\log n})$ bits.  Summarizing, there exists a universal function $J_1$ such that, for each $i< h$ and each $v\in V_i\cap V_{i+1}$, there exists $\ell\in f(v)\cap S^+_i\cap S^+_{i+1}$ and a bitstring $\overrightarrow{\tau}_1(v,\ell,i)$ of length $\Oh(\log{\log n})$ such that $J_1(\sigma(v,i),\overrightarrow{\tau}_1(v,\ell,i))=\sigma_{B_i}(\ell)$.

  In the second step we take $\sigma_{B_i}(\ell)$ onto $\sigma_{B_{i+1}}(\ell)$.  By property~\ref{path_splice} of $B_1,\ldots,B_h$, the path $P_{B_{i+1}}(\ell)$ can be obtained by replacing a subpath $xyz$ in $P_{B_{i}}(\ell)$ with a subpath $xy'z$.  In terms of the signature function $\sigma$, this difference between $P_{B_i}(\ell)$ and $P_{B_{i+1}}(\ell)$ can be described by an $O(\log\log n)$-bit encoding of $\depth_{B_i}(x)$ and two $\lceil\log (12a)\rceil$-bit integers that give the index of $y'$ in the list of children of $x$ and the index of $z$ in the list of children of $y'$.  Although this is the only change to the root-to-$\ell$ path $P_{B_i}(\ell)$, there is an additional change to $\sigma_{B_i}(\ell)$ that needs to be recorded.  Let $p:=\parent_{B_i}(\ell)$ be the parent of $\ell$.  Then the set of children of $p$ in $B_i$ and the set of children of $p$ in $B_{i+1}$ are different.  There are $B_i$-children of $p$ in $S^+_{i}\setminus S^+_{i+1}$ that are not in $B_{i+1}$ and there are $B_{i+1}$-children of $p$ in $S^+_{i+1}\setminus S^+_i$ that are not in $B_i$.  Therefore, the last integer in $\sigma_{B_i}(\ell)$ may be different from the last integer in $\sigma_{B_{i+1}}(\ell)$ even though each of these describes the edge $p\ell$. When this occurs, it can also be described by a $\lceil \log(12a)\rceil$-bit integer, namely the last integer in $\sigma_{B_{i+1}}(\ell)$.  To summarize, there exists a universal function $J_2:(\{0,1\}^*)^2\to \{0,1\}^*$ such that, for each $i\in[h-1]$ and each $\ell\in S^+_{i}\cap S^+_{i+1}$, there exists a bitstring $\overrightarrow{\tau}_2(\ell,i)$ of length $\Oh(\sqrt{\log n})$ such that $\sigma_{B_{i+1}}(\ell)=J_2(\sigma_{B_i}(\ell),\overrightarrow{\tau}_2(\ell,i))$.  

  In the third and final step, we take $\sigma_{B_{i+1}}(\ell)$ onto $\sigma(v,i+1)$.  Since $\ell\in f(v)\cap S^+_{i+1}$, the node $x'_{B_{i+1}}(v)$ is a $B_{i+1}$-ancestor of $\ell$.  Therefore, the signature $\sigma(v,i+1)$ can be obtained from $\sigma_{B_{i+1}}(\ell)$, $\depth_{B_{i+1}}(x_{B_{i+1}}(v))$, and the $\lceil\log(12a)\rceil$-bit integer $c_{v,i+1}$ that gives the index of $x_{B_{i+1}}(v)$ in the list of children of $x'_{B_{i+1}}(v)$.   Formally, there exists a universal function $J_3$ such that, for each $i< h$ and each $v\in V_i\cap V_{i+1}$, there exists a bitstring $\overrightarrow{\tau}_3(v,i)$ of length $\Oh(\sqrt{\log n})$ such that $J_3(\sigma_{B_{i+1}}(\ell),\overrightarrow{\tau}_3(v,i))=\sigma(v,i+1)$.
  
  Summarizing, the forward transition code
  \[
     \overrightarrow{\tau}(v,i):=\langle \overrightarrow{\tau}_1(v,\ell,i),\; \overrightarrow{\tau}_2(\ell, i),\;\overrightarrow{\tau}_3(v,i)\rangle 
  \]
  has length $\Oh(\sqrt{\log n})$, and the transition function $J$ on $\overrightarrow{\tau}(v,i)$ is defined by
  \[
     J(\sigma(v,i),\overrightarrow{\tau}(v,i)) := J_3(J_2(J_1(\sigma(v,i), \overrightarrow{\tau}_1(v,\ell,i)), \overrightarrow{\tau}_2(\ell, i)), \overrightarrow{\tau}_3(v,i)) = \sigma(v,i+1) \enspace .
  \]
  As mentioned above, the backward transition code $\overleftarrow{\tau}(v,i)$, also of length $O(\sqrt{\log n})$, is defined symmetrically (using property~\ref{plus_chunking_up} in place of~\ref{plus_chunking_down}) so that $J(\sigma(v,i),\overleftarrow{\tau}(v,i)) = \sigma(v,i-1)$.

  \paragraph{\boldmath The second part $\lambda_2(v,i)$:}
  This gives all the pieces needed for the second part of the label for a vertex $(v,i)$ of $G$:
  \[
    \mathdefin{\lambda_2(v,i)} := \langle \sigma(v,i),\; \varphi(v),\; \alpha(v,i),\; \overleftarrow{\tau}(v,i),\; \overrightarrow{\tau}(v,i) \rangle  \enspace .
  \]
  The complete vertex label for a vertex $(v,i)$ of $G$ is
  \[
    \mathdefin{\mu(v,i)} := \langle \lambda_1(v,i),\; \lambda_2(v,i) \rangle
  \]
  which has length $\log n + \Oh(\sqrt{\log n})$, for fixed $k$.  (For $k$ that grows with $n$, $|\mu(v,i)|=\log n + \Oh(k\sqrt{\log n})$.)

  \paragraph{\boldmath Adjacency testing:}

  We now describe a procedure that defines the output of the adjacency tester $A$.
  The operation of this procedure makes use of the fact that $\lambda_2(v,i)$ contains all the information needed to evaluate $\sigma(v,i)$, $\sigma(v,i-1)=J(\sigma(v,i),\overleftarrow{\tau}(v,i))$, and $\sigma(v,i+1)=J(\sigma(v,i),\overrightarrow{\tau}(v,i))$.

  Given the labels $\mu(v,i)$ and $\mu(w,j)$ for two vertices $(v,i)$ and $(w,j)$ of $G$, the procedure first uses $\lambda_1(v,i)$ and $\lambda_1(w,j)$ (along with the function $F$) to test if $i=j$, $i=j+1$ or $i=j-1$.  If none of these is the case, then the procedure can immediately conclude that $(v,i)(w,j)$ is not an edge of $G$, since $i\neq j$ and $ij$ is not an edge of $P$.  Since $G$ is a subgraph of $H\boxtimes P$, this implies that $A(\mu(v,i),\mu(w,j)):=0$.  Otherwise, the procedure deals with one of the following cases:
  \begin{enumerate}
    \item If $i=j$, then the procedure uses $\sigma(v,i)$ and $\alpha^0(v,i)$ to compute $\sigma(w',i)$ and $\varphi(w')$ for each $w'\in N^+_H(v)$ such that $(w',i)\in N_G((v,i))$.  For each of these the procedure checks if $(\sigma(w',i),\varphi(w'))=(\sigma(w,i),\varphi(w))$. If so, then $(w,i)=(w',i)\in N_G((v,i))$ so $A(\mu(v,i),\mu(w,j)):=1$.

    If this test fails for each of the entries in $\alpha^0(v,i)$, then the procedure concludes that $w\not\in X_{v,i}$ or $(w,i)\not\in N_G((v,i))$.  This still leaves the possibility that $v\in X_{w,i}$ and $(v,i)\in N_G((w,j))$.  The procedure checks this in exactly the same manner, but with the roles of $\mu(v,i)$ and $\mu(w,j)$ reversed.  If this test also fails to return a result, then $A(\mu(v,i),\mu(w,j)):=0$.

    \item If $i=j+1$ (so $j=i-1$), then the procedure computes $\sigma(v,i-1)=J(\sigma(v,i),\overleftarrow{\tau}(v,i))$. Then using the entries from $\alpha^-(v,i)$, it searches for an entry such that $(\sigma(w',i-1),\varphi(w'))=(\sigma(w,j),\varphi(w))$.  If such an entry is found, then $A(\mu(v,i),\mu(w,j)):=1$.

    If this search fails, then the procedure computes $\sigma(w,j+1)=J(\sigma(w,j),\overrightarrow{\tau}(w,j))$.  Then using the entries from $\alpha^+(w,j)$, it searches for an entry such that $(\sigma(v',j+1),\varphi(v'))=(\sigma(v,i),\varphi(v))$.  If such an entry is found, then $A(\mu(v,i),\mu(w,j)):=1$.  If this search also fails, then $A(\mu(v,i),\mu(w,j)):=0$.

    \item If $i=j-1$ (so $j=i+1$), then the procedure computes $\sigma(v,i+1)=J(\sigma(v,i),\overrightarrow{\tau}(v,i))$. Then using the entries from $\alpha^+(v,i)$, it searches for an entry such that $(\sigma(w',i+1),\varphi(w'))=(\sigma(w,j),\varphi(w))$.  If such an entry is found, then $A(\mu(v,i),\mu(w,j)):=1$.

    If this search fails, then the procedure computes $\sigma(w,j-1)=J(\sigma(w,j),\overleftarrow{\tau}(w,j))$. Then using the entries in $\alpha^-(w,j)$, it searches for an entry such that $(\sigma(v',j-1),\varphi(v'))=(\sigma(v,i),\varphi(v))$. If such an entry is found, then $A(\mu(v,i),\mu(w,j)):=1$. If this search also fails, then $A(\mu(v,i),\mu(w,j)):=0$.
  \end{enumerate}
  This concludes the description of the adjacency testing procedure $A$ as well as the $(\log n+O(\sqrt{\log n}))$-bit adjacency labelling scheme for $\calR_k$.

  \paragraph{Clique labelling and identity testing:}
  Thus far, we have described a $g_0$-orientable $(\log n+O(\sqrt{\log n}))$-bit adjacency labelling scheme for $\calR_k$ (for any fixed $k$).  Next we describe clique labels, the local identifiers, and the identity testing procedure that makes this into a strongly $g_0$-orientable unweighted $(g_1,g_2,g_3)$ mixed labelling scheme.

  Let $K$ be a clique in $G^+$ and let $V_K:=\{v\mid (v,j)\in K\}$. Then $K\subseteq V_K\times\{i_K-1,i_K\}$ and $V_K\subseteq V_{i_K}$ for some $i_K\in[h]$ (with $i_K\ge 2$ unless $h=1$).  Define 
  \[
    \lambda_1(K):=\langle \rho(i_K), \overleftarrow{\pi}(i_K), \overrightarrow{\pi}(i_K)\rangle \enspace.
  \]    
  Since $K$ is a clique in $G^+\subseteq H\boxtimes P$, $V_K$ is a clique in the original treewidth-$k$ graph $H$ (before adding edges within bags of $\calT$).  Therefore $|V_K|\le k+1$ and there exists $v_K\in V_K$ such that $V_K\subseteq \{v_K\}\cup N^+_H(v_K)$.  Since $K\subseteq V_K\times\{i_K-1,i_K\}$, $|K|\le 2(k+1) = 2k+2$.  Furthermore, at least one of $(v_K,i_K)$ or $(v_K,i_K-1)$ is in $K$. Therefore $v_K\in V_{i_K}$ (and $v_K\in V_{i_K-1}$ if $i_K\ge 2$), so $\sigma(v_K, i_K)$ and $\overleftarrow{\tau}(v_K,i_K)$ are defined.  We can now mimic the adjacency labelling.  Define 
  \begin{align*}
    \alpha^-(K)
      & := \langle\; \langle |\sigma(w,i_K-1)|, c_{w,i_K-1}, \varphi(w)\rangle \mid (w,i_K-1)\in K \;\rangle \\
    \alpha^0(K)
      & := \langle\; \langle |\sigma(w,i_K)|, c_{w,i_K}, \varphi(w)\rangle \mid (w,i_K)\in K \;\rangle   \\
    \alpha(K)
      & := \langle\; \alpha^-(K),\; \alpha^0(K)  \;\rangle \enspace .
  \end{align*}
  Since $|K|\le 2k+2$, $|\alpha(K)|\in\Oh(k\sqrt{\log n})$.  By Property~\ref{clique_property}, $f(v_K) \subseteq f(w)$ for each $w \in V_K$. As discussed above, this implies that $\sigma(w,i_K)$ is completely determined by $\sigma(v_K,i_K)$, $|\sigma(w,i_K)|$, and the integer $c_{w,i_K}$ (specifically, by the prefix of $\sigma(v_K,i_K)$ of length $|\sigma(w,i_K)|-\lceil\log(12a)\rceil$ together with $c_{w,i_K}$). Define
  \[
    \mu(K) := \langle\; \lambda_1(K),\; \sigma(v_K,i_K),\; \overleftarrow{\tau}(v_K,i_K),\; \alpha(K) \;\rangle \enspace .
  \]
  Then,
  \begin{align*}
    |\mu(K)|
      & \le |\lambda_1(K)| + |\sigma(v_K,i_K)| + |\overleftarrow{\tau}(v_K,i_K)| + \Oh(k(\sqrt{\log n} + \log\log n)) \\
      & \le \log n - \log|S^+_{i_K}| + (\log a)\depth_{B_{i_K}}(x_{B_{i_K}}(v_K)) + \Oh(\log_a n + k\sqrt{\log n}) \\
      & \le \log n - \log|S^+_{i_K}| + (\log a)\height(B_{i_K}) + \Oh(k\sqrt{\log n}) \\
      & \le \log n + \Oh(\log a+ \log_a n + k\sqrt{\log n}) \\
      & \le \log n + \Oh(k\sqrt{\log n}) \enspace .
  \end{align*}

  To see that $\mu$ is injective when restricted to cliques in $G^+$, suppose $\mu(K)=\mu(L)$ for two cliques $K$ and $L$ in $G^+$.
  The bitstring $\lambda_1(K)$ completely determines $\rho(i_K)$ and $\rho(i_K-1)=F(\rho(i_K),\overleftarrow{\pi}(i_K))$.
  From $\sigma(v_K,i_K)$ and $\overleftarrow{\tau}(v_K,i_K)$, we obtain the signature $\sigma(v_K,i_K-1)=J(\sigma(v_K,i_K), \overleftarrow{\tau}(v_K,i_K))$. The bitstrings $\sigma(v_K,i_K-1)$ and $\alpha^-(K)$ determine $\{(\sigma(w,i_K-1),\varphi(w)):(w,i_K-1)\in K\}$ because each $\sigma(w,i_K-1)$ is uniquely reconstructed from $\sigma(v_K,i_K-1)$, $|\sigma(w,i_K-1)|$, and $c_{w,i_K-1}$.  The bitstrings $\sigma(v_K,i_K)$ and $\alpha^0(K)$ determine $\{(\sigma(w,i_K),\varphi(w)):(w,i_K)\in K\}$ similarly.  Thus, $\mu(K)$ completely determines $\{(\rho(j),\sigma(w,j),\varphi(w)):(w,j)\in K\}$.  Since $\mu(K)=\mu(L)$, we have
  \[
    \{(\rho(j),\sigma(w,j),\varphi(w)):(w,j)\in K\} = \{(\rho(j),\sigma(w,j),\varphi(w)):(w,j)\in L\} \enspace .
  \]
  Since the triple $(\rho(j),\sigma(w,j),\varphi(w))$ is unique for each $(w,j)\in V(G)$, this implies that $K=L$.

  \paragraph{Local identifiers:}

  For each vertex $(u,j)\in K$, the local identifier for $(u,j)$ is a single integer $\kappa(K,(u,j))\in[2k+2]$ that indicates the position of the entry for $(u,j)$ in $\alpha(K)$.  This integer also implicitly encodes whether $j=i_K-1$  or $j=i_K$.  Thus $|\kappa(K,(u,j))| \le \lceil\log(2k+2)\rceil \in \Oh(\log(k+2))$.  For fixed $k$, this satisfies the requirement $g_2(n)\in\Oh(1)$.

  \paragraph{Identity testing:}

  We now describe the identity tester $I$.  Let $K$ be a clique in $G^+$, let $(u,j)\in K$, and let $(v,i)$ be a vertex of $G$.  Let $(v_K,i_K)$ be as in the description of $\mu(K)$, so that $K\subseteq V_K\times\{i_K-1,i_K\}$ and $f(v_K)\subseteq f(w)$ for each $(w,j)\in K$.   
  
  Given the labels $\mu(K)=\langle\lambda_1(K),\sigma(v_K,i_K),\overleftarrow{\tau}(v_K,i_K),\alpha(K)\rangle$ and $\mu(v,i)=\langle \lambda_1(v,i),\lambda_2(v,i)\rangle$, along with the identifier $\kappa(K,(u,j))$, the identity tester proceeds as follows. First, the tester examines $\lambda_1(v,i)$ and $\lambda_1(K)$ to determine if $i=i_K$ or $i=i_K-1$.  If neither of these is the case, then $j \neq i$ so $(u,j)\neq (v,i)$ and the tester outputs $I(\mu(K),\kappa(K,(u,j)),\mu(v,i)):=0$.

  If $i\in\{i_K-1,i_K\}$, then the tester locates the entry $\langle |\sigma(u,j)|, c_{u,j}, \varphi(u)\rangle$ in $\alpha(K)$ using $\kappa(K,(u,j))$, which encodes the position of this entry in $\alpha(K)$. The location of this entry also determines whether $j=i_K-1$ or $j=i_K$.  Since the tester now knows the values of $j$ and $i$, it can immediately determine that $I(\mu(K),\kappa(K,(u,j)),\mu(v,i)):=0$ if $j\neq i$.  Otherwise, $j=i$ and there are two cases to consider:
  \begin{itemize}
    \item If $i=i_K$ then the tester extracts $\sigma(u,i)$ from $\sigma(v_K,i_K)$ using $|\sigma(u,i)|$ and $c_{u,i}$.  If $(\sigma(u,i), \varphi(u)) = (\sigma(v,i), \varphi(v))$, then $(u,j)=(v,i)$ and the tester outputs $1$; otherwise it outputs $0$.
    \item If $i=i_K-1$ then the tester computes $\sigma(v_K,i)=\sigma(v_K,i_K-1)=J(\sigma(v_K,i_K),\overleftarrow{\tau}(v_K,i_K))$.  The tester extracts $\sigma(u,i)$ from $\sigma(v_K,i)$ using $|\sigma(u,i)|$ and $c_{u,i}$.  If $(\sigma(u,i), \varphi(u)) = (\sigma(v,i), \varphi(v))$, then $(u,j)=(v,i)$ and the tester outputs $1$; otherwise it outputs $0$.
  \end{itemize}

  \paragraph{Simulating weights by multiplicity:} Finally, we explain how to make the unweighted $(g_1,g_2,g_3)$ mixed labelling scheme described so far into a strongly $g_0$-orientable weighted mixed labelling scheme. By \cref{lem:nice_weights,obs:normalization_orientability}, it suffices to construct a strongly $g_0$-orientable weak weighted mixed labelling scheme; the additive constants in the vertex and clique label lengths are absorbed into $g_1$ and $g_3$. Let $\tau=(G^+,G,\omega)\in\calT_{\calR_k}$ be a nice triple and let $n:=|V(G)|$. Thus, $\omega:V(G)\to\N^+$ and $n\le\omega(G)\le 2n$, so $\log\omega(G)\le\log n+1$.

  The only modification needed is in the definition of the sets $S_1,\ldots,S_h$. Define $\omega(v,i):=0$ for each $(v,i)\in (V(H)\times\{0,\ldots,h+1\})\setminus V(G)$.  For each $i\in[h]$ and each $v\in V_i$, define $\delta(v,i) := \omega(v,i-1)+\omega(v,i)+\omega(v,i+1)$.  For each $i\in[h]$ and each $v\in V_i$, at least one of $(v,i-1)$, $(v,i)$, or $(v,i+1)$ is in $V(G)$, so $\delta(v,i) \ge 1$. Furthermore, $\sum_{i\in[h]}\sum_{v\in V_i}\delta(v,i) \le 3\omega(G)$ since each vertex $(v,i)$ of $G$ contributes $\omega(v,i)$ to $\delta(v,j)$ for at most three values of $j$. Fix $\epsilon \in (0,1/(6\omega(G))]$. For each $i\in\{1,\ldots,h\}$, define 
  \[
     S_i := \{ x_f(v)-\epsilon j: v\in V_i, j\in\{0,\ldots,\delta(v,i)-1\}\} \enspace .
  \]
  With this new definition of $S_1,\ldots,S_h$, everything proceeds exactly as described above. We define $(S^+_1,\ldots,S^+_h):= Q(S_1,\ldots,S_h)$ so that $S^+_1,\ldots,S^+_h$ satisfy properties \ref{plus_superset}--\ref{plus_chunking_up}. Note that $\sum_{i=1}^h |S_i| = \sum_{i\in[h]}\sum_{v\in V_i}\delta(v,i) \le 3\omega(G)$ so, by property~\ref{plus_size}, $\sum_{i=1}^h |S^+_i| \in \Oh(\omega(G))$.

  Before going further, we give an intuitive explanation of why this works.
  Replacing each leaf $x_f(v)$ in $B_i$ with a cluster of $\delta(v,i)$ leaves ensures that the subtree of $B_i$ rooted at $x'_{B_i}(v)$ contains at least $\delta(v,i)\ge \omega(v,i)$ leaves.  This places a lower bound on the height of the subtree of $B_i$ rooted at $x'_{B_i}(v)$.  Since $x_{B_i}(v)$ is a child of $x'_{B_i}(v)$, this places a lower bound on the height of the subtree of $B_i$ rooted at $x_{B_i}(v)$.  Since the height of $B_i$ is $\log_a |S^+_i| + \Oh(1)$, the lower bound on the height of the subtree rooted at $x_{B_i}(v)$ upper bounds the depth of $x_{B_i}(v)$ in $B_i$.  We now make this precise.

  For each $i\in[h]$ and each $v\in V_i$, let $S_i(v) := \{ x_f(v)-\epsilon j: j\in\{0,\ldots,\delta(v,i)-1\}\}$. Since $S_i(v) \subseteq f(v)\cap S_i \subseteq f(v)\cap S^+_i$, the node $x'_{B_i}(v) = \lca'_{B_i}(f(v)\cap S^+_i)$ is a $B_i$-ancestor of $\lca'_{B_i}(S_i(v))$. The node $\lca'_{B_i}(S_i(v))$ is an internal node of $B_i$ whose subtree contains all $\delta(v,i)$ leaves of $S_i(v)$. Because a subtree of height $s$ in a weight-balanced B-tree of parameter $a$ has at most $6a^s$ leaves, we have $6a^s \ge \delta(v,i)$, so $s \ge \log_a \delta(v,i) - O(1)$. Therefore,
  \begin{align*}
    \depth_{B_i}(x_{B_i}(v))
      & = \depth_{B_i}(x'_{B_i}(v)) + 1 \\
      & \le \depth_{B_i}(\lca'_{B_i}(S_i(v))) + 1 \\
      & \le \height(B_i) - \log_a(\delta(v,i)/6) + 1 \\
      & \le \log_a|S^+_i| - \log_a\delta(v,i) + \Oh(1) \enspace .
  \end{align*}
  Since $\delta(v,i) \ge \omega(v,i)$, this immediately implies that 
  \[
      \depth_{B_i}(x_{B_i}(v)) \le  \log_a|S^+_i| - \log_a\omega(v,i) + \Oh(1) \enspace .
  \]
  It follows that the length of $\sigma(v,i)$ satisfies
  \begin{align*}
    |\sigma(v,i)| 
      & = |\sigma_{B_i}(x_{B_i}(v))| \\
      & \le \lceil\log(12a)\rceil \depth_{B_i}(x_{B_i}(v)) \\
      & \le (5+\log a) (\log_a|S^+_i| - \log_a\omega(v,i) + \Oh(1)) \\
      & \le \log|S^+_i| - \log\omega(v,i) + \Oh(\log_a |S^+_i| + \log a) \\
      & \le \log|S^+_i| - \log\omega(v,i) + \Oh(\sqrt{\log n}) \enspace .
  \end{align*}
  Then, as in the unweighted case, the length of the label $\lambda_2(v,i)$ satisfies
  \begin{align*}
    |\lambda_2(v,i)|
      & \le |\sigma(v,i)| + \Oh(k\sqrt{\log n}) \\
      & \le \log |S^+_i|-\log \omega(v,i) + \Oh(k\sqrt{\log n}) \enspace .
  \end{align*}
  The code $\rho:[h]\to\{0,1\}^*$ is constructed from the weight function $\psi(i):=|S_i^+|$ using \cref{shannon_fano_elias_code}.  Therefore, 
  $|\rho(i)| \le \log\psi([h]) - \log|S^+_i| + \Oh(1) \le \log\omega(G) - \log|S^+_i| + \Oh(1)$.  Thus, 
  $|\lambda_1(v,i)|\le \log\omega(G) - \log |S^+_i| + \Oh(\log\log n)$.  Therefore,
  \begin{align*}
    |\mu(v,i)|
      & = |\lambda_1(v,i)|+|\lambda_2(v,i)| + \Oh(\log\log n) \\
      & \le \log\omega(G) - \log |S^+_i| + \log |S^+_i|-\log \omega(v,i) + \Oh(k\sqrt{\log n}) \\
      & = \log\omega(G) - \log\omega(v,i) + \Oh(k\sqrt{\log n}) \enspace ,
  \end{align*}
  which satisfies the requirement $g_1(n)\in O(k\sqrt{\log n})$.

  For the clique label $\mu(K)$, since at least one of $(v_K,i_K-1)$ or $(v_K,i_K)$ belongs to $K$, we have $\delta(v_K,i_K) \ge \min_{u\in K}\omega(u)$, and the same depth bound gives $|\sigma(v_K,i_K)| \le \log|S^+_{i_K}| - \log\delta(v_K,i_K) + \Oh(\sqrt{\log n})$. Therefore, we similarly obtain
  \begin{align*}
    |\mu(K)|
     & \le |\lambda_1(K)| + |\sigma(v_K,i_K)| + \Oh(k\sqrt{\log n}) \\
     & \le \log\omega(G) - \log|S^+_{i_K}| + \log|S^+_{i_K}| - \log\delta(v_K,i_K) + \Oh(k\sqrt{\log n}) \\
     & \le \log\omega(G) - \log\min_{u\in K}\omega(u) + \Oh(k\sqrt{\log n}) \enspace ,
  \end{align*}
  which satisfies the requirement $g_3(n)\in\Oh(k\sqrt{\log n})$.

  \paragraph{\boldmath $g_0$-orientability:}  
  For each nice triple $\tau=(G^+,G,\omega)\in\calT_{\calR_k}$, let $(\mu_\tau,\kappa_\tau)$ be the mixed labelling of $\tau$ for $(A,I)$ described above, and let $\calS$ be the indexed collection of these labellings over all nice triples in $\calT_{\calR_k}$. Then $(\calS,A,I)$ is a weak weighted mixed labelling scheme for $\calR_k$. Fix $n\in\N^+$, and let $U_n$ be the $n$-universal graph defined by $(\calS,A)$, using only nice triples in $\calT'_{\calR_k,\le n}$. For every such triple, $|V(G)|\le n$ and $\omega(G)\le 2n$ and $\sum_i|S_i^+|\in\Oh(\omega(G))\subseteq\Oh(n)$. Consequently, the bounds on the lengths of the auxiliary label components established above are uniform over these triples, with $n$ now denoting the universal-graph parameter rather than the order of a particular graph. We now describe an orientation of $U_n$ in which each vertex has outdegree at most $g_0(n)$.

  Let $\beta\beta'$ be an edge of $U_n$. Then there exists a nice triple $\tau=(G^+,G,\omega)\in\calT'_{\calR_k,\le n}$ and an edge $e=(v,i)(w,j)\in E(G)$ such that $\beta=\mu_\tau(v,i)$ and $\beta'=\mu_\tau(w,j)$. Since $G^+\in\calR_k$, $G^+$ is a subgraph of $H\boxtimes P$ for some graph $H$ of treewidth at most $k$ and path $P=(1,\ldots,h)$. Since $(v,i)(w,j)\in E(G)$, we have $|i-j|\le 1$. In $H$, the ancestor ordering $<$ of $\calT$ is a strict partial order. Exactly one of the following holds for the pair of vertices $(v,i)$ and $(w,j)$:
  \begin{itemize}[nosep,nolistsep]
    \item $v\neq w$: then either $w < v$ (so $w\in N^+_H(v)$) or $v < w$ (so $v\in N^+_H(w)$). We orient the edge $\beta\beta'$ towards the endpoint whose $H$-coordinate is smaller in $<$. That is, if $w\in N^+_H(v)$, we orient $\beta\beta'$ towards $\beta'$; if $v\in N^+_H(w)$, we orient $\beta\beta'$ towards $\beta$.
    \item $v = w$: then $i\neq j$ and $|i-j|=1$. We orient the edge $\beta\beta'$ towards the endpoint with smaller $P$-coordinate. That is, if $j=i-1$, we orient $\beta\beta'$ towards $\beta'$; if $i=j-1$, we orient $\beta\beta'$ towards $\beta$.
  \end{itemize}

  We now bound the maximum contribution to the outdegree of $\beta$ from this orientation rule. An edge incident to $\beta=\mu_\tau(v,i)$ is oriented away from $\beta$ towards $\beta'=\mu_\tau(w,j)$ only when $w\in N^+_H(v)$ (with $j\in\{i-1,i,i+1\}$) or when $w=v$ and $j=i-1$. In every such case, $w\in \{v\}\cup N^+_H(v)$, so $(w,j)$ corresponds to an entry in $\alpha(v,i)$. We consider the three cases corresponding to the operation of the adjacency tester $A(\beta,\beta')$:
  \begin{enumerate}
    \item $i=j$: In this case, the edge is oriented away from $\beta$ only if $w\in N^+_H(v)$. The tester examines the at most $6k+5$ entries in $\alpha^0(v,i)$ and finds the entry $\langle |\sigma(w,i)|, c_{w,i}, \varphi(w)\rangle$. Using $\sigma(v,i)$ and this entry, $\sigma(w,i)$ and $\varphi(w)$ are uniquely determined. Since $\beta$ fixes $\lambda_1(v,i)=\lambda_1(w,j)$, the prefix $\langle \lambda_1(w,j), \sigma(w,j), \varphi(w)\rangle$ of $\beta'$ is fixed by $\beta$ and the entry. There are at most $6k+5$ such entries.
    
    \item $j=i-1$: In this case, the edge is oriented away from $\beta$ if $w\in \{v\}\cup N^+_H(v)$. The tester computes $\sigma(v,i-1)=J(\sigma(v,i),\overleftarrow{\tau}(v,i))$ and searches in $\alpha^-(v,i)$. The matching entry $\langle |\sigma(w,i-1)|, c_{w,i-1}, \varphi(w)\rangle$ uniquely determines $\sigma(w,i-1)$ and $\varphi(w)$. Since $\rho(j)=F(\rho(i),\overleftarrow{\pi}(i))$ is uniquely determined by $\lambda_1(v,i)$, the prefix $\langle \rho(j), \sigma(w,j), \varphi(w)\rangle$ of $\beta'$ is fixed by $\beta$ and the entry. There are at most $6k+6$ such entries.

    \item $j=i+1$: In this case, the edge is oriented away from $\beta$ only if $w\in N^+_H(v)$. The tester computes $\sigma(v,i+1)=J(\sigma(v,i),\overrightarrow{\tau}(v,i))$ and searches in $\alpha^+(v,i)$. The matching entry uniquely determines $\sigma(w,i+1)$ and $\varphi(w)$, and $\rho(j)=F(\rho(i),\overrightarrow{\pi}(i))$ is determined by $\lambda_1(v,i)$. There are at most $6k+5$ such entries.
  \end{enumerate}
  Across all three cases, there are at most $(6k+5)+(6k+6)+(6k+5) = 18k+16$ choices for the prefix $\langle \rho(j), \sigma(w,j), \varphi(w)\rangle$ of $\beta'$. The remaining components of $\beta'$---namely $\overleftarrow{\pi}(j)$, $\overrightarrow{\pi}(j)$, $\alpha(w,j)$, $\overleftarrow{\tau}(w,j)$, and $\overrightarrow{\tau}(w,j)$---have total length $\Oh(k\sqrt{\log n})$, so there are at most $\exp(\Oh(k\sqrt{\log n}))$ choices for these components. Thus, the number of edges oriented away from $\beta$ is at most $(18k+16)\cdot \exp(\Oh(k\sqrt{\log n})) = \exp(\Oh(k\sqrt{\log n})) = g_0(n)$.

  \paragraph{\boldmath Strong $g_0$-orientability:}  
  All that remains is to show that $(\calS,A,I)$ is strongly $g_0$-orientable. Fix $n\in\N^+$. Let $\tau_0=(G^+_0,G_0,\omega_0)\in\calT'_{\calR_k,\le n}$ be a nice triple, let $(\mu_0,\kappa_0):=(\mu_{\tau_0}, \kappa_{\tau_0})$, let $K_0$ be a clique in $G^+_0$, and let $(u_0,j_0)\in K_0$. Define $\zeta:=\mu_0(K_0)$ and $\gamma:=\kappa_0(K_0,(u_0,j_0))$. We must upper bound the number of strings $\beta$ for which there exists a nice triple $\tau=(G^+,G,\omega)\in\calT'_{\calR_k,\le n}$, a clique $K$ in $G^+$, and a vertex $(u,j)\in K$ such that $\mu_\tau(K)=\zeta$, $\kappa_\tau(K,(u,j))=\gamma$, $\mu_\tau(u,j)=\beta$, and $I(\zeta,\gamma,\beta)=1$. Let $\beta$ be such a string. The local identifier $\gamma$ indexes an entry in $\alpha(K) = \langle \alpha^-(K), \alpha^0(K)\rangle$, which gives $|\sigma(u,j)|$, $c_{u,j}$, $\varphi(u)$, and whether $j=i_K$ or $j=i_K-1$:
  \begin{itemize}[nosep,nolistsep] 
    \item If $j=i_K$, then $\lambda_1(u,j)=\lambda_1(K)$ is completely determined by $\zeta$. The string $\sigma(v_K,i_K)$ extracted from $\zeta$ together with $|\sigma(u,i_K)|$ and $c_{u,i_K}$ uniquely determines $\sigma(u,i_K)=\sigma(u,j)$. Along with $\varphi(u)$, the prefix $\langle \lambda_1(u,j), \sigma(u,j), \varphi(u)\rangle$ of $\beta$ is fixed. The remaining components $\alpha(u,j), \overleftarrow{\tau}(u,j), \overrightarrow{\tau}(u,j)$ of $\beta$ have total length $\Oh(k\sqrt{\log n})$. Thus there are at most $\exp(\Oh(k\sqrt{\log n}))$ such strings $\beta$.

    \item If $j=i_K-1$, then $\rho(j)=F(\rho(i_K),\overleftarrow{\pi}(i_K))$ is determined by $\lambda_1(K)$. The tester computes $\sigma(v_K,i_K-1)=J(\sigma(v_K,i_K),\overleftarrow{\tau}(v_K,i_K))$ using the strings $\sigma(v_K,i_K)$ and $\overleftarrow{\tau}(v_K,i_K)$ extracted from $\zeta$. Together with $|\sigma(u,j)|$ and $c_{u,j}$ from $\gamma$, this uniquely determines $\sigma(u,j)$. Along with $\varphi(u)$, the prefix $\langle \rho(j), \sigma(u,j), \varphi(u)\rangle$ of $\beta$ is fixed. The remaining components of $\beta$, namely $\overleftarrow{\pi}(j)$, $\overrightarrow{\pi}(j)$, $\alpha(u,j)$, $\overleftarrow{\tau}(u,j)$, and $\overrightarrow{\tau}(u,j)$, have total length $\Oh(k\sqrt{\log n})$. Thus there are at most $\exp(\Oh(k\sqrt{\log n}))$ such strings $\beta$.
  \end{itemize}
  In either case, the number of choices for $\beta$ is at most $\exp(\Oh(k\sqrt{\log n})) = g_0(n)$. Therefore, $(\calS,A,I)$ is strongly $g_0$-orientable.

  To obtain the full weighted scheme, for each $\tau=(G^+,G,\omega)\in\calT_{\calR_k}$, assign the labelling in $\calS$ of $(G^+,G,\widehat\omega)$, where
  \[
    \widehat\omega(u):=\ceil{\frac{|V(G)|\omega(u)}{\omega(G)}}\enspace.
  \]
  By \cref{lem:nice_weights}, this increases the vertex and clique label-length bounds by at most $1$. If $\tau\in\calT'_{\calR_k,\le n}$, then $\widehat\omega(G)\le\omega(G)\le 2n$ by \cref{obs:nice}. Hence, by \cref{obs:normalization_orientability}, the full weighted scheme is strongly $g_0$-orientable with the same function $g_0$. This proves \cref{gk_labelling_precise}. \qedhere
\end{proof}

Let $\calH_k$ be the class of graphs of treewidth at most $k$.  \citet{EJM23} show that, for each fixed $k\in\N^+$ and each $n\in\N^+$, $\calH_k$ has an $n$-induced-universal graph with $n(\log n)^{\Oh(k^2)}$ vertices and edges.
We finish this appendix by remarking that this can be improved using the ideas above.  In particular, the use of \cref{low_height_tree_decomp} eliminates a step used in \cite{EJM23} that leads to the appearance of $k^2$ in the exponent.

\begin{lem}\label{treewidth_result}
  There exists a constant $c>0$ such that, for each fixed $k\in\N^+$, the class of graphs $\calH_k$ admits an injective $g_0$-orientable $(\log n + \Oh(k\log\log n))$-bit adjacency labelling scheme with $g_0(n)\in \Oh((\log n)^{ck})$.
\end{lem}

\begin{proof}[Proof Sketch]
  Let $G$ be an $n$-vertex graph of treewidth at most $k$.  By \cref{low_height_tree_decomp}, $G$ has a tree-decomposition $\calT=(T,(B_x:x\in V(T)))$ of width at most $6k+5$ where $T$ is a binary tree of height $\Oh(\log n)$. As in the proof of \cref{gk_labelling_precise}, define a strict partial order $<$ on $V(G)$ such that $v<w$ if the home node of $v$ in $\calT$ is a strict $T$-ancestor of the home node of $w$ (resolving ties arbitrarily for vertices with the same home node), and for each $v\in V(G)$ define $N^+_G(v):=\{u\in N_G(v):u<v\}$. Note that $|N^+_G(v)|\le 6k+5$.
  Map each vertex $v$ of $G$ to a closed interval $f(v)$ with integer endpoints in $\{1,\dots,2n\}$ satisfying Properties~\ref{nested_or_disjoint}--\ref{thin} as described in the proof of \cref{gk_labelling_precise}.
  For each $v\in V(G)$, let $x_f(v)$ be the right endpoint of $f(v)$. Since each integer is an endpoint of at most one interval, the points in $S:=\{x_f(v):v\in V(G)\}$ are distinct, and $x_f(v)\in f(v)\cap S$ for every $v\in V(G)$.

  Let $B$ be a balanced binary tree whose leaves are the elements of $S$, ordered as in $\R$. Then $B$ has height $\lceil \log_2 n\rceil \le \log n + 1$.
  For each $v\in V(G)$, define $x_B(v):=\lca_B(f(v)\cap S)$. Since $x_f(v)\in f(v)\cap S$, the set $f(v)\cap S$ is non-empty, so $x_B(v)$ is well-defined. For any $u\in N^+_G(v)$, Property~\ref{intersection_subgraph} implies $f(v)\subset f(u)$, so $\emptyset\neq (f(v)\cap S)\subseteq (f(u)\cap S)$. Consequently, $x_B(u)$ is an ancestor of $x_B(v)$ in $B$.

  For each node $x\in V(B)$, we upper bound the number of vertices $v\in V(G)$ with $x_B(v)=x$:
  \begin{itemize}[nosep,nolistsep]
    \item If $x$ is an internal node of $B$ with left child $x_1$ and right child $x_2$, choose a real number $p_x$ strictly between the maximum leaf in the subtree of $x_1$ and the minimum leaf in the subtree of $x_2$. If $x_B(v)=x$, then $f(v)\cap S$ contains a leaf in the subtree of $x_1$ and a leaf in the subtree of $x_2$, so the interval $f(v)$ contains $p_x$.
    \item If $x$ is a leaf of $B$, then $x$ is a point $p_x\in S$, and $x_B(v)=x$ implies $p_x\in f(v)$.
  \end{itemize}
  In either case, every vertex $v\in V(G)$ with $x_B(v)=x$ satisfies $p_x\in f(v)$. By Property~\ref{thin}, at most $ck\log n$ vertices are mapped to $x$. We can therefore assign each vertex $v\in V(G)$ a colour $\varphi(v)\in [\lfloor ck\log n\rfloor]$ such that any two distinct vertices mapped to the same node of $B$ receive different colours (for instance, using the proper colouring of the interval intersection graph $G_f$).

  For each $v\in V(G)$, define the signature $\sigma(v)$ as the binary string describing the path $P_B(x_B(v))$ from the root of $B$ to $x_B(v)$. Then $|\sigma(v)|\le \log n + 1$. Observe that, for any $v,w\in V(G)$, $(\sigma(v),\varphi(v)) = (\sigma(w),\varphi(w))$ if and only if $v=w$.
  For each $v\in V(G)$ and each $u\in N^+_G(v)$, since $x_B(u)$ is an ancestor of $x_B(v)$, $\sigma(u)$ is the prefix of $\sigma(v)$ of length $|\sigma(u)|$. Define $\alpha(v,u):=\langle |\sigma(u)|,\varphi(u)\rangle$. Then $|\alpha(v,u)|\le \Oh(\log\log n + \log k)$. For each $v\in V(G)$, define
  \[
    \alpha(v):=\langle \alpha(v,u):u\in N^+_G(v)\rangle \enspace .
  \]
  Since $|N^+_G(v)|\le 6k+5$, $|\alpha(v)|\le \Oh(k\log\log n)$. The label for a vertex $v\in V(G)$ is
  \[
    \ell_G(v) := \langle \sigma(v),\varphi(v),\alpha(v)\rangle \enspace ,
  \]
  which has length $|\ell_G(v)| \le \log n + \Oh(k\log\log n)$. Since $(\sigma(v),\varphi(v))$ uniquely identifies $v$, the labelling is injective.

  Given the labels $\ell_G(v)$ and $\ell_G(w)$ of two vertices $v,w\in V(G)$, the adjacency tester $A$ operates as follows:
  \begin{itemize}[nosep,nolistsep]
    \item If $\sigma(w)$ is a prefix of $\sigma(v)$ and $\alpha(v)$ contains an entry $\langle|\sigma(w)|,\varphi(w)\rangle$, then $A$ returns $1$.
    \item Else, if $\sigma(v)$ is a prefix of $\sigma(w)$ and $\alpha(w)$ contains an entry $\langle|\sigma(v)|,\varphi(v)\rangle$, then $A$ returns $1$.
    \item Otherwise, $A$ returns $0$.
  \end{itemize}
  To verify correctness, suppose the first condition holds. The entry $\langle|\sigma(w)|,\varphi(w)\rangle$ in $\alpha(v)$ certifies the existence of some $u\in N^+_G(v)$ such that $|\sigma(u)|=|\sigma(w)|$ and $\varphi(u)=\varphi(w)$. Since both $\sigma(u)$ and $\sigma(w)$ are prefixes of $\sigma(v)$ of the same length, $\sigma(u)=\sigma(w)$. Together with $\varphi(u)=\varphi(w)$, this implies $u=w$, so $w\in N^+_G(v)$ and $vw\in E(G)$. The second case is symmetric. Conversely, if $vw\in E(G)$, then either $w\in N^+_G(v)$ (satisfying the first condition) or $v\in N^+_G(w)$ (satisfying the second condition).

  To see that the scheme is $g_0$-orientable, let $U_n$ be the $n$-universal graph defined by this scheme. For each edge $e=\{\beta,\beta'\}\in E(U_n)$, there exists some graph $G\in\calH_k$ with at most $n$ vertices and an edge $vw\in E(G)$ such that $\beta=\ell_G(v)$ and $\beta'=\ell_G(w)$. If $w\in N^+_G(v)$, orient $e$ towards $\beta'$; if $v\in N^+_G(w)$, orient $e$ towards $\beta$ (breaking ties arbitrarily if different graphs witness opposite orientations).
  Now consider any vertex $\beta=\langle \sigma,\varphi,\alpha\rangle\in V(U_n)$. An edge is oriented away from $\beta$ towards $\beta'$ only if $\beta'$ corresponds to an in-neighbour $w\in N^+_G(v)$ for some witness graph. For each such neighbour, $\alpha$ contains an entry $\langle s,c\rangle=\langle|\sigma(w)|,\varphi(w)\rangle$. This entry completely determines the prefix $\langle \sigma(w),\varphi(w)\rangle$ of $\beta'$, where $\sigma(w)$ is the prefix of $\sigma$ of length $s$ and $\varphi(w)=c$. The only remaining component of $\beta'$ is $\alpha(w)$, which has length at most $c'k\log\log n$ for some constant $c'>0$. Thus, there are at most $2^{c'k\log\log n} = (\log n)^{c'k}$ choices for $\alpha(w)$. Since $\alpha$ contains at most $6k+5$ entries, the number of edges oriented away from $\beta$ in $U_n$ is at most $(6k+5)(\log n)^{c'k} \in \Oh((\log n)^{ck})$ for an absolute constant $c>0$. Therefore, the scheme is $g_0$-orientable with $g_0(n)\in\Oh((\log n)^{ck})$.
\end{proof}

\begin{cor}\label{cor:treewidth_universal}
  For each fixed $k\in\N^+$ and each $n\in\N^+$, the class $\calH_k$ of graphs of treewidth at most $k$ has an $n$-induced-universal graph with $n(\log n)^{\Oh(k)}$ vertices and edges.
\end{cor}
\begin{proof}
  By \cref{treewidth_result}, $\calH_k$ admits an injective $g_0$-orientable $f(n)$-bit adjacency labelling scheme with $f(n)=\log n + \Oh(k\log\log n)$ and $g_0(n)\in\Oh((\log n)^{ck})$. For each $n\in\N^+$, the $n$-universal graph $U_n$ defined by this scheme has $|V(U_n)|\le 2^{f(n)+1} = n(\log n)^{\Oh(k)}$ vertices and, by orientability, at most $|V(U_n)|\cdot g_0(n) \le n(\log n)^{\Oh(k)}$ edges. By injectivity, every $n$-vertex graph in $\calH_k$ is isomorphic to an induced subgraph of $U_n$.
\end{proof}

\end{document}